


\ifx\shlhetal\undefinedcontrolsequence\let\shlhetal\relax\fi

\input amstex
\expandafter\ifx\csname mathdefs.tex\endcsname\relax
  \expandafter\gdef\csname mathdefs.tex\endcsname{}
\else \message{Hey!  Apparently you were trying to
  \string\input{mathdefs.tex} twice.   This does not make sense.} 
\errmessage{Please edit your file (probably \jobname.tex) and remove
any duplicate ``\string\input'' lines}\endinput\fi




\catcode`\X=12\catcode`\@=11

\def\n@wcount{\alloc@0\count\countdef\insc@unt}
\def\n@wwrite{\alloc@7\write\chardef\sixt@@n}
\def\n@wread{\alloc@6\read\chardef\sixt@@n}
\def\r@s@t{\relax}\def\v@idline{\par}\def\@mputate#1/{#1}
\def\l@c@l#1X{\firstpart.#1}\def\gl@b@l#1X{#1}\def\t@d@l#1X{{}}

\def\crossrefs#1{\ifx\all#1\let\tr@ce=\all\else\def\tr@ce{#1,}\fi
   \n@wwrite\cit@tionsout\openout\cit@tionsout=\jobname.cit 
   \write\cit@tionsout{\tr@ce}\expandafter\setfl@gs\tr@ce,}
\def\setfl@gs#1,{\def\@{#1}\ifx\@\empty\let\next=\relax
   \else\let\next=\setfl@gs\expandafter\xdef
   \csname#1tr@cetrue\endcsname{}\fi\next}
\def\m@ketag#1#2{\expandafter\n@wcount\csname#2tagno\endcsname
     \csname#2tagno\endcsname=0\let\tail=\all\xdef\all{\tail#2,}
   \ifx#1\l@c@l\let\tail=\r@s@t\xdef\r@s@t{\csname#2tagno\endcsname=0\tail}\fi
   \expandafter\gdef\csname#2cite\endcsname##1{\expandafter
     \ifx\csname#2tag##1\endcsname\relax?\else\csname#2tag##1\endcsname\fi
     \expandafter\ifx\csname#2tr@cetrue\endcsname\relax\else
     \write\cit@tionsout{#2tag ##1 cited on page \folio.}\fi}
   \expandafter\gdef\csname#2page\endcsname##1{\expandafter
     \ifx\csname#2page##1\endcsname\relax?\else\csname#2page##1\endcsname\fi
     \expandafter\ifx\csname#2tr@cetrue\endcsname\relax\else
     \write\cit@tionsout{#2tag ##1 cited on page \folio.}\fi}
   \expandafter\gdef\csname#2tag\endcsname##1{\expandafter
      \ifx\csname#2check##1\endcsname\relax
      \expandafter\xdef\csname#2check##1\endcsname{}%
      \else\immediate\write16{Warning: #2tag ##1 used more than once.}\fi
      \multit@g{#1}{#2}##1/X%
      \write\t@gsout{#2tag ##1 assigned number \csname#2tag##1\endcsname\space
      on page \number\count0.}%
   \csname#2tag##1\endcsname}}

\def\multit@g#1#2#3/#4X{\def\t@mp{#4}\ifx\t@mp\empty%
      \global\advance\csname#2tagno\endcsname by 1 
      \expandafter\xdef\csname#2tag#3\endcsname
      {#1\number\csname#2tagno\endcsnameX}%
   \else\expandafter\ifx\csname#2last#3\endcsname\relax
      \expandafter\n@wcount\csname#2last#3\endcsname
      \global\advance\csname#2tagno\endcsname by 1 
      \expandafter\xdef\csname#2tag#3\endcsname
      {#1\number\csname#2tagno\endcsnameX}
      \write\t@gsout{#2tag #3 assigned number \csname#2tag#3\endcsname\space
      on page \number\count0.}\fi
   \global\advance\csname#2last#3\endcsname by 1
   \def\t@mp{\expandafter\xdef\csname#2tag#3/}%
   \expandafter\t@mp\@mputate#4\endcsname
   {\csname#2tag#3\endcsname\lastpart{\csname#2last#3\endcsname}}\fi}
\def\t@gs#1{\def\all{}\m@ketag#1e\m@ketag#1s\m@ketag\t@d@l p
\let\realscite\scite
\let\realstag\stag
   \m@ketag\gl@b@l r \n@wread\t@gsin
   \openin\t@gsin=\jobname.tgs \re@der \closein\t@gsin
   \n@wwrite\t@gsout\openout\t@gsout=\jobname.tgs }
\outer\def\localtags{\t@gs\l@c@l}
\outer\def\globaltags{\t@gs\gl@b@l}
\outer\def\newlocaltag#1{\m@ketag\l@c@l{#1}}
\outer\def\newglobaltag#1{\m@ketag\gl@b@l{#1}}

\newif\ifpr@ 
\def\m@kecs #1tag #2 assigned number #3 on page #4.%
   {\expandafter\gdef\csname#1tag#2\endcsname{#3}
   \expandafter\gdef\csname#1page#2\endcsname{#4}
   \ifpr@\expandafter\xdef\csname#1check#2\endcsname{}\fi}
\def\re@der{\ifeof\t@gsin\let\next=\relax\else
   \read\t@gsin to\t@gline\ifx\t@gline\v@idline\else
   \expandafter\m@kecs \t@gline\fi\let \next=\re@der\fi\next}
\def\pretags#1{\pr@true\pret@gs#1,,}
\def\pret@gs#1,{\def\@{#1}\ifx\@\empty\let\n@xtfile=\relax
   \else\let\n@xtfile=\pret@gs \openin\t@gsin=#1.tgs \message{#1} \re@der 
   \closein\t@gsin\fi \n@xtfile}

\newcount\sectno\sectno=0\newcount\subsectno\subsectno=0
\newif\ifultr@local \def\ultralocal{\ultr@localtrue}
\def\firstpart{\number\sectno}
\def\lastpart#1{\ifcase#1 \or a\or b\or c\or d\or e\or f\or g\or h\or 
   i\or k\or l\or m\or n\or o\or p\or q\or r\or s\or t\or u\or v\or w\or 
   x\or y\or z \fi}

\def\resetall{\global\advance\sectno by 1\subsectno=0
   \gdef\firstpart{\number\sectno}\r@s@t}
\def\resetsub{\global\advance\subsectno by 1
   \gdef\firstpart{\number\sectno.\number\subsectno}\r@s@t}
\def\newsection#1\par{\resetall\vskip0pt plus.3\vsize\penalty-250
   \vskip0pt plus-.3\vsize\bigskip\bigskip
   \message{#1}\leftline{\bf#1}\nobreak\bigskip}
\def\subsection#1\par{\ifultr@local\resetsub\fi
   \vskip0pt plus.2\vsize\penalty-250\vskip0pt plus-.2\vsize
   \bigskip\smallskip\message{#1}\leftline{\bf#1}\nobreak\medskip}


\newdimen\marginshift

\newdimen\margindelta
\newdimen\marginmax
\newdimen\marginmin

\def\margininit{       
\marginmax=3 true cm                  
				      
\margindelta=0.1 true cm              
\marginmin=0.1true cm                 
\marginshift=\marginmin
}    

\def\t@gsjj#1,{\def\@{#1}\ifx\@\empty\let\next=\relax\else\let\next=\t@gsjj
   \def\@@{p}\ifx\@\@@\else
   \expandafter\gdef\csname#1cite\endcsname##1{\citejj{##1}}
   \expandafter\gdef\csname#1page\endcsname##1{?}
   \expandafter\gdef\csname#1tag\endcsname##1{\tagjj{##1}}\fi\fi\next}
\newif\ifshowstuffinmargin
\showstuffinmarginfalse
\def\jjtags{\showstuffinmargintrue
\ifx\all\relax\else\expandafter\t@gsjj\all,\fi}

\def\tagjj#1{\realstag{#1}\mginpar{\zeigen{#1}}}
\def\citejj#1{\zeigen{#1}\mginpar{\rechnen{#1}}}

\def\rechnen#1{\expandafter\ifx\csname stag#1\endcsname\relax ??\else
                           \csname stag#1\endcsname\fi}

\newdimen\theight

\def\marginfont{\sevenrm}

\def\trymarginbox#1{\setbox0=\hbox{\marginfont\hskip\marginshift #1}%
		\global\marginshift\wd0 
		\global\advance\marginshift\margindelta}

\def \mginpar#1{%
\ifvmode\setbox0\hbox to \hsize{\hfill\rlap{\marginfont\quad#1}}%
\ht0 0cm
\dp0 0cm
\box0\vskip-\baselineskip
\else 
             \vadjust{\trymarginbox{#1}%
		\ifdim\marginshift>\marginmax \global\marginshift\marginmin
			\trymarginbox{#1}%
                \fi
             \theight=\ht0
             \advance\theight by \dp0    \advance\theight by \lineskip
             \kern -\theight \vbox to \theight{\rightline{\rlap{\box0}}%
\vss}}\fi}


\def\t@gsoff#1,{\def\@{#1}\ifx\@\empty\let\next=\relax\else\let\next=\t@gsoff
   \def\@@{p}\ifx\@\@@\else
   \expandafter\gdef\csname#1cite\endcsname##1{\zeigen{##1}}
   \expandafter\gdef\csname#1page\endcsname##1{?}
   \expandafter\gdef\csname#1tag\endcsname##1{\zeigen{##1}}\fi\fi\next}
\def\verbatimtags{\showstuffinmarginfalse
\ifx\all\relax\else\expandafter\t@gsoff\all,\fi}
\def\zeigen#1{\hbox{$\langle$}#1\hbox{$\rangle$}}
\def\margincite#1{\ifshowstuffinmargin\mginpar{\rechnen{#1}}\fi}

\def\(#1){\edef\dot@g{\ifmmode\ifinner(\hbox{\noexpand\etag{#1}})
   \else\noexpand\eqno(\hbox{\noexpand\etag{#1}})\fi
   \else(\noexpand\ecite{#1})\fi}\dot@g}

\newif\ifbr@ck
\def\eat#1{}
\def\[#1]{\br@cktrue[\br@cket#1'X]}
\def\br@cket#1'#2X{\def\temp{#2}\ifx\temp\empty\let\next\eat
   \else\let\next\br@cket\fi
   \ifbr@ck\br@ckfalse\br@ck@t#1,X\else\br@cktrue#1\fi\next#2X}
\def\br@ck@t#1,#2X{\def\temp{#2}\ifx\temp\empty\let\neext\eat
   \else\let\neext\br@ck@t\def\temp{,}\fi
   \def\teemp{#1}\ifx\teemp\empty\else\rcite{#1}\fi\temp\neext#2X}
\def\resetbr@cket{\gdef\[##1]{[\rtag{##1}]}}
\def\references{\resetbr@cket\newsection References\par}

\newtoks\symb@ls\newtoks\s@mb@ls\newtoks\p@gelist\n@wcount\ftn@mber
    \ftn@mber=1\newif\ifftn@mbers\ftn@mbersfalse\newif\ifbyp@ge\byp@gefalse
\def\defm@rk{\ifftn@mbers\n@mberm@rk\else\symb@lm@rk\fi}
\def\n@mberm@rk{\xdef\m@rk{{\the\ftn@mber}}%
    \global\advance\ftn@mber by 1 }
\def\rot@te#1{\let\temp=#1\global#1=\expandafter\r@t@te\the\temp,X}
\def\r@t@te#1,#2X{{#2#1}\xdef\m@rk{{#1}}}
\def\b@@st#1{{$^{#1}$}}\def\str@p#1{#1}
\def\symb@lm@rk{\ifbyp@ge\rot@te\p@gelist\ifnum\expandafter\str@p\m@rk=1 
    \s@mb@ls=\symb@ls\fi\write\f@nsout{\number\count0}\fi \rot@te\s@mb@ls}
\def\byp@ge{\byp@getrue\n@wwrite\f@nsin\openin\f@nsin=\jobname.fns 
    \n@wcount\currentp@ge\currentp@ge=0\p@gelist={0}
    \re@dfns\closein\f@nsin\rot@te\p@gelist
    \n@wread\f@nsout\openout\f@nsout=\jobname.fns }
\def\m@kelist#1X#2{{#1,#2}}
\def\re@dfns{\ifeof\f@nsin\let\next=\relax\else\read\f@nsin to \f@nline
    \ifx\f@nline\v@idline\else\let\t@mplist=\p@gelist
    \ifnum\currentp@ge=\f@nline
    \global\p@gelist=\expandafter\m@kelist\the\t@mplistX0
    \else\currentp@ge=\f@nline
    \global\p@gelist=\expandafter\m@kelist\the\t@mplistX1\fi\fi
    \let\next=\re@dfns\fi\next}
\def\symbols#1{\symb@ls={#1}\s@mb@ls=\symb@ls} 
\def\bigsymbol{\textstyle}
\symbols{\bigsymbol\ast,\dagger,\ddagger,\sharp,\flat,\natural,\star}
\def\ftnumbers{\ftn@mberstrue} \def\ftsymbols{\ftn@mbersfalse}
\def\paginal{\byp@ge} \def\resetftnumbers{\ftn@mber=1}
\def\ftnote#1{\defm@rk\expandafter\expandafter\expandafter\footnote
    \expandafter\b@@st\m@rk{#1}}

\long\def\jump#1\endjump{}
\def\ssum{\mathop{\lower .1em\hbox{$\textstyle\Sigma$}}\nolimits}

\def\qed{\nobreak\kern 1em \vrule height .5em width .5em depth 0em}
\def\newneq{\hbox{\rlap{\hbox to 1\wd9{\hss$=$\hss}}\raise .1em 
   \hbox to 1\wd9{\hss$\scriptscriptstyle/$\hss}}}
\def\subsetne{\setbox9 = \hbox{$\subset$}\mathrel{\hbox{\rlap
   {\lower .4em \newneq}\raise .13em \hbox{$\subset$}}}}
\def\supsetne{\setbox9 = \hbox{$\subset$}\mathrel{\hbox{\rlap
   {\lower .4em \newneq}\raise .13em \hbox{$\supset$}}}}

\def\vbar{\mathchoice{\vrule height6.3ptdepth-.5ptwidth.8pt\kern-.8pt}
   {\vrule height6.3ptdepth-.5ptwidth.8pt\kern-.8pt}
   {\vrule height4.1ptdepth-.35ptwidth.6pt\kern-.6pt}
   {\vrule height3.1ptdepth-.25ptwidth.5pt\kern-.5pt}}
\def\f@dge{\mathchoice{}{}{\mkern.5mu}{\mkern.8mu}}
\def\b@c#1#2{{\rm \mkern#2mu\vbar\mkern-#2mu#1}}
\def\b@b#1{{\rm I\mkern-3.5mu #1}}
\def\b@a#1#2{{\rm #1\mkern-#2mu\f@dge #1}}
\def\bb#1{{\count4=`#1 \advance\count4by-64 \ifcase\count4\or\b@a A{11.5}\or
   \b@b B\or\b@c C{5}\or\b@b D\or\b@b E\or\b@b F \or\b@c G{5}\or\b@b H\or
   \b@b I\or\b@c J{3}\or\b@b K\or\b@b L \or\b@b M\or\b@b N\or\b@c O{5} \or
   \b@b P\or\b@c Q{5}\or\b@b R\or\b@a S{8}\or\b@a T{10.5}\or\b@c U{5}\or
   \b@a V{12}\or\b@a W{16.5}\or\b@a X{11}\or\b@a Y{11.7}\or\b@a Z{7.5}\fi}}

\catcode`\X=11 \catcode`\@=12


\expandafter\ifx\csname citeadd.tex\endcsname\relax
\expandafter\gdef\csname citeadd.tex\endcsname{}
\else \message{Hey!  Apparently you were trying to
\string\input{citeadd.tex} twice.   This does not make sense.} 
\errmessage{Please edit your file (probably \jobname.tex) and remove
any duplicate ``\string\input'' lines}\endinput\fi

\sectno=-1   
\localtags
\NoBlackBoxes
\define\mr{\medskip\roster}
\define\sn{\smallskip\noindent}
\define\mn{\medskip\noindent}
\define\bn{\bigskip\noindent}
\define\ub{\underbar}
\define\wilog{\text{without loss of generality}}
\define\ermn{\endroster\medskip\noindent}
\define\dbca{\dsize\bigcap}
\define\dbcu{\dsize\bigcup}
\define \nl{\newline}
\magnification=\magstep 1
\documentstyle{amsppt}

{    
\catcode`@11

\ifx\alicetwothousandloaded@\relax
  \endinput\else\global\let\alicetwothousandloaded@\relax\fi

\gdef\subjclass{\let\savedef@\subjclass
 \def\subjclass##1\endsubjclass{\let\subjclass\savedef@
   \toks@{\def\usualspace{{\rm\enspace}}\eightpoint}%
   \toks@@{##1\unskip.}%
   \edef\thesubjclass@{\the\toks@
     \frills@{{\noexpand\rm2000 {\noexpand\it Mathematics Subject
       Classification}.\noexpand\enspace}}%
     \the\toks@@}}%
  \nofrillscheck\subjclass}
} 

\pageheight{8.5truein}
\topmatter
\title{Are ${\frak{\lowercase{a}}}$ and ${\frak{\lowercase{d}}}$ your cup of
tea?}\endtitle
\author {Saharon Shelah \thanks {\null\newline I would like to thank 
Alice Leonhardt for the beautiful typing. \null\newline
The research was supported by The Israel Science Foundation founded by the
Israel Academy of Sciences and Humanities. Publication 700.} \endthanks}
\endauthor  
\affil{Institute of Mathematics\\
 The Hebrew University\\
 Jerusalem, Israel
 \medskip
 Rutgers University\\
 Mathematics Department\\
 New Brunswick, NJ  USA} \endaffil
\mn
\abstract  We show that consistently, every MAD family has cardinality
strictly bigger than the dominating number, that is ${\frak a} > {\frak d}$,
thus solving one of the oldest problems on cardinal invariants of the
continuum.  The method is a contribution to the theory of iterated
forcing for making the continuum large.   \endabstract
\endtopmatter
\document  

\expandafter\ifx\csname alice2jlem.tex\endcsname\relax
  \expandafter\xdef\csname alice2jlem.tex\endcsname{\the\catcode`@}
\else \message{Hey!  Apparently you were trying to
\string\input{alice2jlem.tex}  twice.   This does not make sense.}
\errmessage{Please edit your file (probably \jobname.tex) and remove
any duplicate ``\string\input'' lines}\endinput\fi

\expandafter\ifx\csname bib4plain.tex\endcsname\relax
  \expandafter\gdef\csname bib4plain.tex\endcsname{}
\else \message{Hey!  Apparently you were trying to \string\input
  bib4plain.tex twice.   This does not make sense.}
\errmessage{Please edit your file (probably \jobname.tex) and remove
any duplicate ``\string\input'' lines}\endinput\fi

\def\renewcommand{\newcommand}	       
\edef\cite{\the\catcode`@}%
\catcode`@ = 11
\let\@oldatcatcode = \cite
\chardef\@letter = 11
\chardef\@other = 12
%
%
%
%
\def\@innerdef#1#2{\edef#1{\expandafter\noexpand\csname #2\endcsname}}%
%
%
\@innerdef\@innernewcount{newcount}%
\@innerdef\@innernewdimen{newdimen}%
\@innerdef\@innernewif{newif}%
\@innerdef\@innernewwrite{newwrite}%
%
%
%
\def\@gobble#1{}%
%
%
%
\ifx\inputlineno\@undefined
   \let\@linenumber = \empty 
\else
   \def\@linenumber{\the\inputlineno:\space}%
\fi
%
%
%
\def\@futurenonspacelet#1{\def\cs{#1}%
   \afterassignment\@stepone\let\@nexttoken=
}%
\begingroup 
\def\\{\global\let\@stoken= }%
\\ 
\endgroup
\def\@stepone{\expandafter\futurelet\cs\@steptwo}%
\def\@steptwo{\expandafter\ifx\cs\@stoken\let\@@next=\@stepthree
   \else\let\@@next=\@nexttoken\fi \@@next}%
\def\@stepthree{\afterassignment\@stepone\let\@@next= }%
%
%
%
\def\@getoptionalarg#1{%
   \let\@optionaltemp = #1%
   \let\@optionalnext = \relax
   \@futurenonspacelet\@optionalnext\@bracketcheck
}%
%
%
\def\@bracketcheck{%
   \ifx [\@optionalnext
      \expandafter\@@getoptionalarg
   \else
      \let\@optionalarg = \empty
      \expandafter\@optionaltemp
   \fi
}%
\def\@@getoptionalarg[#1]{%
   \def\@optionalarg{#1}%
   \@optionaltemp
}%
%
%
%
\def\@nnil{\@nil}%
\def\@fornoop#1\@@#2#3{}%
\def\@for#1:=#2\do#3{%
   \edef\@fortmp{#2}%
   \ifx\@fortmp\empty \else
      \expandafter\@forloop#2,\@nil,\@nil\@@#1{#3}%
   \fi
}%
\def\@forloop#1,#2,#3\@@#4#5{\def#4{#1}\ifx #4\@nnil \else
       #5\def#4{#2}\ifx #4\@nnil \else#5\@iforloop #3\@@#4{#5}\fi\fi
}%
\def\@iforloop#1,#2\@@#3#4{\def#3{#1}\ifx #3\@nnil
       \let\@nextwhile=\@fornoop \else
      #4\relax\let\@nextwhile=\@iforloop\fi\@nextwhile#2\@@#3{#4}%
}%
%
%
%
\@innernewif\if@fileexists
\def\@testfileexistence{\@getoptionalarg\@finishtestfileexistence}%
\def\@finishtestfileexistence#1{%
   \begingroup
      \def\extension{#1}%
      \immediate\openin0 =
         \ifx\@optionalarg\empty\jobname\else\@optionalarg\fi
         \ifx\extension\empty \else .#1\fi
         \space
      \ifeof 0
         \global\@fileexistsfalse
      \else
         \global\@fileexiststrue
      \fi
      \immediate\closein0
   \endgroup
}%
%
%
%
%
\def\bibliographystyle#1{%
   \@readauxfile
   \@writeaux{\string\bibstyle{#1}}%
}%
\let\bibstyle = \@gobble
%
%
\let\bblfilebasename = \jobname
\def\bibliography#1{%
   \@readauxfile
   \@writeaux{\string\bibdata{#1}}%
   \@testfileexistence[\bblfilebasename]{bbl}%
   \if@fileexists
      \nobreak
      \@readbblfile
   \fi
}%
\let\bibdata = \@gobble
%
%
\def\nocite#1{%
   \@readauxfile
   \@writeaux{\string\citation{#1}}%
}%
\@innernewif\if@notfirstcitation
%
%
\def\cite{\@getoptionalarg\@cite}%
%
%
\def\@cite#1{%
   \let\@citenotetext = \@optionalarg
   \printcitestart
   \nocite{#1}%
   \@notfirstcitationfalse
   \@for \@citation :=#1\do
   {%
      \expandafter\@onecitation\@citation\@@
   }%
   \ifx\empty\@citenotetext\else
      \printcitenote{\@citenotetext}%
   \fi
   \printcitefinish
}%
\def\@onecitation#1\@@{%
   \if@notfirstcitation
      \printbetweencitations
   \fi
   \expandafter \ifx \csname\@citelabel{#1}\endcsname \relax
      \if@citewarning
         \message{\@linenumber Undefined citation `#1'.}%
      \fi
      \expandafter\gdef\csname\@citelabel{#1}\endcsname{%
\strut
\vadjust{\vskip-\dp\strutbox
\vbox to 0pt{\vss\parindent0cm \leftskip=\hsize 
\advance\leftskip3mm
\advance\hsize 4cm\strut\openup-4pt 
\rightskip 0cm plus 1cm minus 0.5cm ?  #1 ?\strut}}
         {\tt
            \escapechar = -1
            \nobreak\hskip0pt
            \expandafter\string\csname#1\endcsname
            \nobreak\hskip0pt
         }%
      }%
   \fi
   \csname\@citelabel{#1}\endcsname
   \@notfirstcitationtrue
}%
%
%
\def\@citelabel#1{b@#1}%
%
%
\def\@citedef#1#2{\expandafter\gdef\csname\@citelabel{#1}\endcsname{#2}}%
%
%
%
\def\@readbblfile{%
   \ifx\@itemnum\@undefined
      \@innernewcount\@itemnum
   \fi
   \begingroup
      \def\begin##1##2{%
         \setbox0 = \hbox{\biblabelcontents{##2}}%
         \biblabelwidth = \wd0
      }%
      \def\end##1{}
      %
      %
      \@itemnum = 0
      \def\bibitem{\@getoptionalarg\@bibitem}%
      \def\@bibitem{%
         \ifx\@optionalarg\empty
            \expandafter\@numberedbibitem
         \else
            \expandafter\@alphabibitem
         \fi
      }%
      \def\@alphabibitem##1{%
         \expandafter \xdef\csname\@citelabel{##1}\endcsname {\@optionalarg}%
         \ifx\biblabelprecontents\@undefined
            \let\biblabelprecontents = \relax
         \fi
         \ifx\biblabelpostcontents\@undefined
            \let\biblabelpostcontents = \hss
         \fi
         \@finishbibitem{##1}%
      }%
      \def\@numberedbibitem##1{%
         \advance\@itemnum by 1
         \expandafter \xdef\csname\@citelabel{##1}\endcsname{\number\@itemnum}%
         \ifx\biblabelprecontents\@undefined
            \let\biblabelprecontents = \hss
         \fi
         \ifx\biblabelpostcontents\@undefined
            \let\biblabelpostcontents = \relax
         \fi
         \@finishbibitem{##1}%
      }%
      \def\@finishbibitem##1{%
         \biblabelprint{\csname\@citelabel{##1}\endcsname}%
         \@writeaux{\string\@citedef{##1}{\csname\@citelabel{##1}\endcsname}}%
         \ignorespaces
      }%
      %
      %
      \let\em = \bblem
      \let\newblock = \bblnewblock
      \let\sc = \bblsc
      \frenchspacing
      \clubpenalty = 4000 \widowpenalty = 4000
      \tolerance = 10000 \hfuzz = .5pt
      \everypar = {\hangindent = \biblabelwidth
                      \advance\hangindent by \biblabelextraspace}%
      \bblrm
      \parskip = 1.5ex plus .5ex minus .5ex
      \biblabelextraspace = .5em
      \bblhook
      \input \bblfilebasename.bbl
   \endgroup
}%
%
%
\@innernewdimen\biblabelwidth
\@innernewdimen\biblabelextraspace
%
%
%
\def\biblabelprint#1{%
   \noindent
   \hbox to \biblabelwidth{%
      \biblabelprecontents
      \biblabelcontents{#1}%
      \biblabelpostcontents
   }%
   \kern\biblabelextraspace
}%
%
%
%
\def\biblabelcontents#1{{\bblrm [#1]}}%
%
%
\def\bblrm{\rm}%
%
%
\def\bblem{\it}%
%
%
\def\bblsc{\ifx\@scfont\@undefined
              \font\@scfont = cmcsc10
           \fi
           \@scfont
}%
%
%
\def\bblnewblock{\hskip .11em plus .33em minus .07em }%
%
%
\let\bblhook = \empty
%
%
%
\def\printcitestart{[}
\def\printcitefinish{]}
\def\printbetweencitations{, }
\def\printcitenote#1{, #1}
%
%
%
\let\citation = \@gobble
%
%
%
\@innernewcount\@numparams
%
%
\def\newcommand#1{%
   \def\@commandname{#1}%
   \@getoptionalarg\@continuenewcommand
}%
%
%
\def\@continuenewcommand{%
   \@numparams = \ifx\@optionalarg\empty 0\else\@optionalarg \fi \relax
   \@newcommand
}%
%
%
\def\@newcommand#1{%
   \def\@startdef{\expandafter\edef\@commandname}%
   \ifnum\@numparams=0
      \let\@paramdef = \empty
   \else
      \ifnum\@numparams>9
         \errmessage{\the\@numparams\space is too many parameters}%
      \else
         \ifnum\@numparams<0
            \errmessage{\the\@numparams\space is too few parameters}%
         \else
            \edef\@paramdef{%
               \ifcase\@numparams
                  \empty  No arguments.
               \or ####1%
               \or ####1####2%
               \or ####1####2####3%
               \or ####1####2####3####4%
               \or ####1####2####3####4####5%
               \or ####1####2####3####4####5####6%
               \or ####1####2####3####4####5####6####7%
               \or ####1####2####3####4####5####6####7####8%
               \or ####1####2####3####4####5####6####7####8####9%
               \fi
            }%
         \fi
      \fi
   \fi
   \expandafter\@startdef\@paramdef{#1}%
}%
%
%
%
%
\def\@readauxfile{%
   \if@auxfiledone \else 
      \global\@auxfiledonetrue
      \@testfileexistence{aux}%
      \if@fileexists
         \begingroup
            \endlinechar = -1
            \catcode`@ = 11
            \input \jobname.aux
         \endgroup
      \else
         \message{\@undefinedmessage}%
         \global\@citewarningfalse
      \fi
      \immediate\openout\@auxfile = \jobname.aux
   \fi
}%
%
%
\newif\if@auxfiledone
\ifx\noauxfile\@undefined \else \@auxfiledonetrue\fi
%
%
%
%
\@innernewwrite\@auxfile
\def\@writeaux#1{\ifx\noauxfile\@undefined \write\@auxfile{#1}\fi}%
%
%
%
\ifx\@undefinedmessage\@undefined
   \def\@undefinedmessage{No .aux file; I won't give you warnings about
                          undefined citations.}%
\fi
%
%
\@innernewif\if@citewarning
\ifx\noauxfile\@undefined \@citewarningtrue\fi
%
%
%
\catcode`@ = \@oldatcatcode


\def\widestnumber#1#2{}

\def\rm{\fam0 \tenrm}

\def\fakesubhead#1\endsubhead{\bigskip\noindent{\bf#1}\par}



%
%
%

%

\font\textrsfs=rsfs10
\font\scriptrsfs=rsfs7
\font\scriptscriptrsfs=rsfs5

\newfam\rsfsfam
\textfont\rsfsfam=\textrsfs
\scriptfont\rsfsfam=\scriptrsfs
\scriptscriptfont\rsfsfam=\scriptscriptrsfs

\edef\oldcatcodeofat{\the\catcode`\@}
\catcode`\@11

\def\Cal@@#1{\noaccents@ \fam \rsfsfam #1}

\catcode`\@\oldcatcodeofat


\expandafter\ifx \csname margininit\endcsname \relax\else\margininit\fi

\mn

\newpage

\head {Content} \endhead  \resetall 
\bn
\S0 $\quad$ Introduction
\mn
\S1 $\quad$ CON$({\frak a} > {\frak d})$
\mr
\item "{{}}"  [We prove the consistency mentioned in the title,
relying on the theory of CS iteration of nep forcing (from
\cite{Sh:630}, this is a concise version).]
\ermn
\S2 $\quad$ On CON$({\frak a} > {\frak d})$ revisited with FS, non-transitive
memory of non-well ordered \nl

$\quad$ length
\mr
\item "{{}}"  [Does not depend on \S1.
We define ``FSI template", a depth on their subsets on which
we shall do induction; we are interested just in the cases where the depth is
$< \infty$.  Now the iteration is defined and its properties are proved
simultaneously by induction on the depth.  After we have understood such
iterations sufficiently well, we proceed to prove the consistency in details]. 
\ermn
\S3 $\quad$ Eliminating the measurable
\mr
\item "{{}}"  [In \S2, for checking the criterion which appears there for
having ``${\frak a}$ large", we have used ultrapower by some 
$\kappa$-complete ultrafilter.
Here we construct templates of cardinality, e.g. $\aleph_3$ which
satisfy the criterion; by constructing them 
such that any sequence of $\omega$-tuples of
appropriate length has a (big) subsequence which is ``convergent".]
\ermn
\S4 $\quad$  On related cardinal invariants
\mr
\item "{{}}"  [We prove e.g. the consistency of ${\frak u} < {\frak a}$.
Here the forcing notions are not so definable.]
\endroster
\newpage

\head {\S0 Introduction} \endhead  \resetall \sectno=0
\bigskip

We deal with the theory of iteration of forcing notions for the
continuum and prove CON$({\frak a} > {\frak d})$ and related results.  We
present it in several perspectives; so \S2 + \S3 does not depend on \S1; and \S4
does not depend on \S1, \S2, \S3.  
In \S2 we introduce and investigate iterations which are of finite
support but with non transitive memory and linear, non well ordered
length and prove CON$({\frak a} > {\frak d})$ using a measurable.
In \S4 we answer also related questions
(${\frak u} < {\frak a},{\frak i} < {\frak a}$); in \S3, relying on \S2 
we eliminate the use of a measurable, and in \S1 we rely heavily on \cite{Sh:630}.

Very basically, the difference between ${\frak a}$ and ${\frak b},{\frak d}$
which we use is that ${\frak a}$ speaks on a set, whereas ${\frak b}$
is witnessed by a sequence and ${\frak d}$ by a quite directed family;
it essentially deals with cofinality;
so every unbounded subsequence is a witness as well, i.e.
the relevant relation is transitive; when ${\frak b} =
{\frak d}$ things are smooth, otherwise the situation is still similar.
This manifests itself by using ultrapowers for
some $\kappa$-complete ultrafilter (in model theoretic outlook), and by
using ``convergent sequence" (see \cite{Sh:300}, or the existence of
Av, the average, in \cite{Sh:c}) in \S2, \S3, respectively.
The meaning of ``model theoretic outlook", is that by experience set
theorists starting to hear an explanation of the forcing tend to think
of an elementary embedding
$\bold j:\bold V \rightarrow M$ and then the limit practically does not
make sense (though of course we can translate).  Note that ultrapowers by e.g.
an ultrafilter on $\kappa$, preserve any witness for a cofinality of a linear
order being $\ge \kappa^+$ (or the cofinality of a $\kappa^+$-directed partial
order), as the set of old elements is cofinal and a cofinal subset of a cofinal
subset is a cofinal subset.   On the other hand, the ultrapower always
``increase" a set of cardinality at least the completeness of the ultrafilter.
\bn
\centerline {$* \qquad * \qquad *$}
\bn
This is one of the oldest problems on cardinal invariants of the continuum
(see \cite{vD}).  It was mostly thought that consistently ${\frak a} >
{\frak d}$ and that the natural way to proceed is by CS iteration $\langle P_i,
{\underset\tilde {}\to Q_i}:i < \omega_2 \rangle$ of proper
${}^\omega \omega$-bounding forcing notions, starting with $\bold V \models$ GCH,
and $|P_i| = \aleph_1$ for $i < \omega_2$ and 
${\underset\tilde {}\to Q_i}$ ``deal" with one
MAD family ${\Cal A}_i \in V^{P_i},{\Cal A}_i \subseteq [\omega]^{\aleph_0}$,
adding an infinite subset of $\omega$ almost disjoint to every 
$A \in {\Cal A}_i$.  The needed iteration theorem holds
by \cite[Ch.V,\S4]{Sh:f}, saying that in $\bold V^{P_{\omega_2}},{\frak d} =
{\frak b} = \aleph_1$ and no cardinal is collapsed, \ub{but} 
the single step forcing is not known to exist.  This has been
explained in details in \cite{Sh:666}.
\bn
We do not go in this
way but in a totally different direction involving making the continuum
large, so we still do not know \nl
\ub{\stag{0.1} Problem}  Is $ZFC + 2^{\aleph_0} + \aleph_2 + {\frak a} >
{\frak d}$ consistent?

To clarify our idea, let $D$ be a normal ultrafilter on $\kappa$, a measurable
cardinal and consider a c.c.c. forcing notion $P$ and
\mr
\item "{$(a)$}"  a sequence ${\underset\tilde {}\to {\bar f}} = \langle
{\underset\tilde {}\to f_\alpha}:\alpha < \kappa^+ \rangle$ of $P$-names
such that \nl
$\Vdash_P ``\langle {\underset\tilde {}\to f_\alpha}:\alpha < \kappa^+
\rangle$ is $<^*$-increasing cofinal in ${}^\omega \omega$" \nl
(so ${\underset\tilde {}\to {\bar f}}$ exemplifies $\Vdash_P ``{\frak b} =
{\frak d} = \kappa^+"$)
\sn
\item "{$(b)$}"  a sequence $\langle {\underset\tilde {}\to A_\alpha}:\alpha
< \alpha^* \rangle$ of $P$-names such that
\nl
$\Vdash_P ``\{{\underset\tilde {}\to A_\alpha}:\alpha < \alpha^*\}$ is
MAD that is $\alpha \ne \beta \Rightarrow {\underset\tilde {}\to A_\alpha} \cap
{\underset\tilde {}\to A_\beta}$ is finite and
${\underset\tilde {}\to A_\alpha} \in [\omega]^{\aleph_0}"$.
\ermn
Now $P_1 = P^\kappa/D$ also is a c.c.c. forcing notion by 
{\L}o\'s theorem for
$L_{\kappa,\kappa}$; let $\bold j:P \rightarrow P_1$ be the canonical
embedding; moreover, under the canonical identification we have
$P \prec_{L_{\kappa,\kappa}} P_1$.  So also $\Vdash_{P_1} ``
{\underset\tilde {}\to f_\alpha} \in {}^\omega \omega"$, recalling that
${\underset\tilde {}\to f_\alpha}$ actually consists of $\omega$ maximal
antichains of $P$ (or think of $({\Cal H}(\chi),\in)^\kappa/D,\chi$ large
enough).  Similarly $\Vdash_{P_1} ``{\underset\tilde {}\to f_\alpha} <^*
{\underset\tilde {}\to f_\beta}$ if $\alpha < \beta < \kappa^+"$.

Now, if $\Vdash_{P_1} ``\underset\tilde {}\to g \in {}^\omega \omega"$,
then $\underset\tilde {}\to g = \langle {\underset\tilde {}\to g_\varepsilon}:
\varepsilon < \kappa \rangle/D,\Vdash_P 
``{\underset\tilde {}\to g_\varepsilon} \in {}^\omega \omega"$ so for some
$\alpha^* < \kappa^+$ we have $\Vdash_P 
``{\underset\tilde {}\to g_\varepsilon} <^* {\underset\tilde {}\to f_\alpha}$
for $\varepsilon < \kappa"$ hence by {\L}o\'s theorem
$\Vdash_{P_1} ``\underset\tilde {}\to g <^* f_\alpha"$ (so before the
identification this means $\Vdash_{P_1} ``\underset\tilde {}\to g <^*
\bold j(f_\alpha)"$), so $\langle
{\underset\tilde {}\to f_\alpha}:\alpha < \kappa^+ \rangle$ exemplifies also
$\Vdash_{P_1} ``{\frak b} = {\frak d} = \kappa^+"$.

On the other hand $\langle {\underset\tilde {}\to A_\alpha}:\alpha <
\alpha^* \rangle$ cannot exemplify that ${\frak a} \le \kappa^+$ in
$\bold V^{P_1}$ because $\alpha^* \ge \kappa^+$ (as $ZFC \models {\frak b} \le 
{\frak a}$) so $\langle {\underset\tilde {}\to A_\alpha}:\alpha < \kappa
\rangle/D$ exemplifies that $\Vdash_{P_1} 
``\{{\underset\tilde {}\to A_\alpha}:\alpha < \alpha^*\}$ is not MAD".

Our original idea here is to start with a FS iteration $\bar Q^0 = \langle
P^0_i,{\underset\tilde {}\to Q^0_i}:i < \kappa^+ \rangle$ of nep c.c.c.
forcing notions, ${\underset\tilde {}\to Q^0_i}$ adding a dominating real, (e.g.
dominating real = Hechler forcing), for $\kappa$ a measurable cardinal and let
$D$ be a $\kappa$-complete uniform ultrafilter on $\kappa$ and $\chi >>
\kappa$.  Then let $L_0 = \kappa^+,\bar Q^1 = \langle P^1_i,Q^1_i:i \in L_1
\rangle$ be $\bar Q^0$ as interpreted in $({\Cal H}(\chi),\in,<^*_\chi)
^\kappa/D$, it
looks like $\bar Q^0$ replacing $\kappa^+$ by $(\kappa^+)^\kappa/D$.  We look
at Lim$(\bar Q^0) = \dbcu_i P_i$ as a subforcing of Lim$(\bar Q^1)$
identifying ${\underset\tilde {}\to Q_i}$ with 
${\underset\tilde {}\to Q_{{\bold j}_0(i)}},\bold j_0$ the canonical
elementary embedding of $\kappa^+$ into $(\kappa^+)^\kappa/D$ (no Mostowski
collapse!).  We continue to define $\bar Q^n$ and then $\bar Q^\omega$ as the
following limit: for the original $i \in \kappa^+$, we use the definition,
otherwise we use direct limit (``founding fathers priviledge" you may say ).  
So $P^i = \text{ Lim}(\bar Q^i)$ is 
$\lessdot$-increasing, continuous when cf$(i) > \aleph_0$; so now we
have a kind of iteration with non transitive memory and not well founded.  We continue
$\kappa^{++}$ times.  Now in $\bold V^{\text{Lim}(\bar Q^{\kappa^{++}})}$, the 
original $\kappa^+$ generic reals exemplify ${\frak b} = {\frak d} =
\kappa^+$, so we know that ${\frak a} \ge \kappa^+$.  To finish assume
$p \Vdash ``\{{\underset\tilde {}\to A_\gamma}:\gamma < \kappa^+\}
\subseteq [\omega]^{\aleph_0}$ is a MAD family".  Each name
${\underset\tilde {}\to A_\gamma}$ is a ``countable object" and so depends
on countably many co-ordinates, so all of them are in Lim$(\bar Q^i)$ for 
some $i < \kappa^{++}$.  In the next stage, $\bar Q^{i+1},
\langle {\underset\tilde {}\to A_\gamma}:\gamma < \kappa \rangle/D$ is a
name of an infinite subset of $\omega$ almost disjoint to
${\underset\tilde {}\to A_\beta}$ for each $\beta < \kappa^+$, contradiction.

All this is a reasonable scheme.  This is done in \S1 but relay on
``nep forcing" from \cite{Sh:630}.  But a self contained another approach in \S2,\S3, 
where the meaning of the iteration is more on the surface (and also, in \S3,
help to eliminate the use of large cardinals).  In \S4 we deal with
the case of an additional cardinal invariant, ${\frak u}$.

Note that just using FS iteration on a non well-ordered linear order $L$
(instead of an ordinal) is impossible by a theorem of Hjorth.  On nonlinear
orders for iterations (history and background) see \cite{RoSh:670}.  
On iteration with
nontransitive memory see \cite{Sh:592}, \cite{Sh:619} and in particular
\cite[\S3]{Sh:619}.

I thank Heike Mildenberger and Juris Steprans for their comments.
\newpage

\head {\S1 On Con$({\frak a} > {\frak d})$} \endhead  \resetall \sectno=1
\bigskip

In this section, we look at it in the context of \cite{Sh:630} and we 
use a measurable.
\bigskip

\definition{\stag{da.1} Definition}  1) Given sets $A_\ell$ of ordinals
for $\ell < n$, we say ${\Cal T}$ is an $(A_0,\dotsc,A_{n-1})$-tree if
${\Cal T} = \dbcu_{k < \omega} {\Cal T}_k$ where ${\Cal T}_k \subseteq 
\{(\eta_0,\dotsc,\eta_\ell,\dotsc,\eta_{n-1}):
\eta_\ell \in {}^k(A_\ell)$ for $\ell < n\}$ and ${\Cal T}$
is ordered by $\bar \eta \le_{\Cal T} \bar \nu \Leftrightarrow \dsize 
\bigwedge_{\ell < n} \eta_\ell \trianglelefteq \nu_\ell$ and we let 
$\bar \eta \upharpoonleft k_1 =: \langle \eta_\ell \restriction 
k_1:\ell < n \rangle$ and demand $\bar \eta \in {\Cal T}_k 
\and k_1 < k \Rightarrow \bar \eta \upharpoonleft k_1 \in
{\Cal T}_{k_1}$.  We call ${\Cal T}$ locally countable if $k \in [1,\omega) 
\and \bar \eta \in {\Cal T}_k \Rightarrow |\{\bar \nu \in {\Cal T}_{k+1}:
\bar \eta \le_{\Cal T} \bar \nu\}| \le \aleph_0$.
Let lim$({\Cal T}) = \{ \langle \eta_\ell:\ell < n \rangle:\eta_\ell \in
{}^\omega(A_\ell)$ for $\ell < k$ and $m < \omega \Rightarrow 
\langle \eta_\ell \restriction m:\ell < n \rangle \in {\Cal T}\}$.  
Lastly for $n_1 \le n$ we let prj lim$_{n_1}({\Cal T}) = 
\{\langle \eta_\ell:\ell < n_1 \rangle:$ for some 
$\eta_{n_1},\dotsc,\eta_{n-1}$ we have $\langle \eta_\ell:\ell < n 
\rangle \in \text{ lim}({\Cal T})\}$; and if $n_1$ is omitted we mean 
$n_1 = n-1$.
\nl
2)
$$
\align
{\frak K} = \bigl\{ \bar{\Cal T}:&\text{for some sets } A,B 
\text{ of ordinals we have} \\
  &(i) \quad \bar{\Cal T} = ({\Cal T}_1,{\Cal T}_2), \\
  &(ii) \quad {\Cal T}_1 \text{ is a locally countable } (A,B)\text{-tree},\\
  &(iii) \quad {\Cal T}_2 \text{ is a locally countable } 
(A,A,B)\text{-tree, and} \\
  &(iv) \quad Q_{\bar{\Cal T}} =: (\text{prj lim}({\Cal T}_1), \text{ prj lim}
({\Cal T}_2)) \text{ is a c.c.c. forcing notion} \\
  &\qquad \text{absolutely under c.c.c. forcing notions (see below)} \bigr\}
\endalign
$$
\mn
2A) We say that $Q_{\bar{\Cal T}}$ is c.c.c. absolutely for c.c.c. forcing if:
for c.c.c. forcing notions $P \lessdot R \text{ we have } 
Q^{\bold V^P}_{\bar{\Cal T}} \lessdot Q^{\bold V^R}_{\bar{\Cal T}}$ so membership, order,
nonorder, compatibility, noncompatibility and being predense over $p$ 
are preserved (the $Q_{\bar{\Cal T}}$'s are snep, from \cite{Sh:630} with 
slight restriction).  Similarly we define ``$Q_{\bar{\Cal T}} \lessdot
Q_{\bar{\Cal T}}$ absolutely under c.c.c. forcing".
\nl
3) For a set or class A of ordinals, ${\frak K}^\kappa_A$ 
is the family of $\bar{\Cal T} \in {\frak K}$ which
are a pair of objects, the first an $(A,B)$-tree and the second an
$(A,A,B)$-trees for some $B$ such 
that $|{\Cal T}_1| \le \kappa,|{\Cal T}_2| \le \kappa$.  For a
cardinal $\kappa$ and a pairing function pr with inverses pr$_1$, pr$_2$, let 
${\frak K}^\kappa_{\text{pr}_1,\gamma} = 
{\frak K}^\kappa_{\{\alpha:\text{pr}_1(\alpha)=\gamma\}}$ and
${\frak K}^\kappa_{\text{pr}_1,<\gamma} = 
{\frak K}^\kappa_{\{\alpha:\text{pr}_i(\alpha) < \gamma\}}$.  Let $|\bar{\Cal T}| =
|{\Cal T}_1| + |{\Cal T}_2|$. \nl
4) Let $\bar{\Cal T},\bar{\Cal T}' \in {\frak K}$, we say $\bold f$ is 
an isomorphism from $\bar{\Cal T}$ onto $\bar{\Cal T}'$ \ub{when}:
$\bold f = (f_1,f_2)$ and for $m = 1,2$ we have: $f_m$ is a one-to-one
function from ${\Cal T}_m$ onto ${\Cal T}'_m$ preserving the level (in the
respective trees), preserving the relations 
$x = y \upharpoonleft k,x \ne y \upharpoonleft k$
and if $f_2((\eta_1,\eta_2,\eta_3)) = (\eta'_1,\eta'_2,\eta'_3),f_1((\nu_1,
\nu_2)) = (\nu'_1,\nu'_2)$ then $[\eta_1 = \nu_1 \Leftrightarrow \eta'_1 =
\nu'_1],[\eta_2 = \nu_1 \Leftrightarrow \eta'_2 = \nu'_1]$. \nl
In this case let $\hat{\bold f}$ be the isomorphism induced by $\bold f$ from
$Q_{\bar{\Cal T}}$ onto $Q_{\bar{\Cal T}'}$.
\enddefinition
\bigskip

\definition{\stag{da.2} Definition}  For $\bar{\Cal T}',\bar{\Cal T}'' 
\in {\frak K}$ let $\bar{\Cal T}' \le_{\frak K} \bar{\Cal T}''$ mean:
\mr
\item "{$(a)$}"  ${\Cal T}'_\ell \subseteq {\Cal T}''_\ell$ (as trees) 
for $\ell =1,2$
\sn
\item "{$(b)$}"  if $\ell \in \{1,2\}$ and $\bar \eta \in {\Cal T}''_\ell
\backslash {\Cal T}'_\ell$ and $\bar \eta \upharpoonleft k \in
{\Cal T}'_\ell$ \ub{then} $k \le 1$
\sn
\item "{$(c)$}"  $Q_{\bar{\Cal T}'} \lessdot Q_{\bar{\Cal T}''}$ 
(absolutely under c.c.c. forcing); note that by (a) + (b) we have: \nl
$x \in Q_{\bar{\Cal T}'} \Rightarrow x \in
Q_{\bar{\Cal T}''}$ and $Q_{\bar{\Cal T}'} \models x \le y \Rightarrow
Q_{\bar{\Cal T}''} \models x \le y$).
\endroster
\enddefinition
\bigskip

\remark{Remark}  The definition is tailored such that the union of an
increasing chain will give a forcing notion which is the union.
\endremark
\bigskip

\proclaim{\stag{da.3} Claim/Definition}  0) $\le_{\frak K}$ is a
partial order of ${\frak K}$. \nl
1) Assume $\langle \bar{\Cal T}[i]:
i < \delta \rangle$ is $\le_{\frak K}$-increasing and 
$\bar{\Cal T}$ is defined by $\bar{\Cal T} = \dbcu_i \bar{\Cal T}[i]$ 
that is ${\Cal T}_m = \dbcu_{i < \delta} {\Cal T}_m[i]$ for $m=1,2$ \ub{then}
\mr
\item "{$(a)$}"  $i < \delta \Rightarrow \bar{\Cal T}[i] 
\le_{\frak K} \bar{\Cal T}$
\sn
\item "{$(b)$}"  $Q_{\bar{\Cal T}} = \dbcu_{i < \delta} 
Q_{\bar{\Cal T}[i]}$.
\ermn
2) Assume $\bar{\Cal T}',\bar{\Cal T} \in {\frak K}$.  \ub{Then} 
there is $\bar{\Cal T}'' \in {\frak K}$ such that $\bar{\Cal T}' 
\le_{\frak K} \bar{\Cal T}''$ and $Q_{\bar{\Cal T}''}$
is isomorphic to $Q_{\bar{\Cal T}'} * 
{\underset\tilde {}\to Q_{\bar{\Cal T}}}$ and this is absolute by c.c.c.
forcing.  Moreover, there is such an isomorphism
extending the identity map from $Q_{\bar{\Cal T}'}$ into
$Q_{\bar{\Cal T}''}$.
\nl
3) There is $\bar{\Cal T} \in {\frak K}^{\aleph_0}_\omega$ such that 
$Q_{\bar{\Cal T}}$ is the trivial forcing. \nl
4) There is $\bar{\Cal T} \in {\frak K}^{\aleph_0}_\omega$ such that $Q_{\bar{\Cal T}}$
is the dominating real forcing.
\endproclaim
\bigskip

\demo{Proof}  See \cite{Sh:630}.
\enddemo
\bigskip

\proclaim{\stag{da.4} Claim} 1) Assume $\bar{\Cal T}[\gamma] \in {\frak K}
_{\text{pr}_1,\gamma}$ for $\gamma < \gamma(*)$.  \ub{Then} for each $\alpha
\le \gamma(*)$ there is $\bar{\Cal T} \langle \alpha \rangle
\in {\frak K}_{\text{pr}_1,\gamma < \gamma(*)}$ such that 
$Q_{\bar{\Cal T} \langle \alpha \rangle}$ is $P_\alpha$ where 
$\langle P_\gamma,{\underset\tilde {}\to Q_\beta}:
\gamma \le \gamma(*),\beta < \gamma(*) \rangle$ is an FS iteration and
${\underset\tilde {}\to Q_\beta} = (Q_{\bar{\Cal T}[\beta]})^{{\bold V}[P_\beta]}$ 
and $\bar{\Cal T} \langle \alpha \rangle \in {\frak K}_{\text{pr}_1,
< \alpha}$ and $\bar{\Cal T} \langle \alpha_1 \rangle \le_{\frak K}
\bar{\Cal T} \langle \alpha_2 \rangle$ for $\alpha_1 \le \alpha_2 \le
\gamma(*),\bar{\Cal T}[\gamma] \le_{\frak K} \bar{\Cal T} \langle \alpha
\rangle$ for $\gamma < \alpha \le \gamma(*)$.  We write $\bar{\Cal T}
\langle \alpha \rangle = \dsize \sum_{\gamma < \alpha} \bar{\Cal T}[\gamma]$.
\nl
2) In part (1), for each $\gamma < \gamma(*)$ there is $\bar{\Cal T}' \in
{\frak K}_{\text{pr}_1,\gamma}$ such that $\bar{\Cal T}',\bar{\Cal T}$ 
are isomorphic over $\bar{\Cal T}[\gamma]$ hence $Q_{\bar{\Cal T}'},
Q_{\bar{\Cal T}}$ are isomorphic over $Q_{\bar{\Cal T}[\gamma]}$. \nl
3) If in addition ${\Cal T}[\gamma] \le_{\frak K} {\Cal T}'[\gamma] \in 
{\frak K}_{\text{pr}_1,\gamma}$ for $\gamma < \gamma(*)$ 
and $\langle P_\gamma,{\underset\tilde {}\to Q'_\beta}:
\gamma \le \gamma(*)$, \nl
$\beta < \gamma(*)
\rangle$ is an FS iteration as above with $P'_{\gamma(*)} = Q_{\bar{\Cal T}'}$,
\ub{then} we find such $\bar{\Cal T}'$ with $\bar{\Cal T} \le_{\frak K} 
\bar{\Cal T}'$.
\endproclaim
\bigskip

\demo{Proof}  Straight.
\enddemo
\bigskip

\proclaim{\stag{da.6} Theorem}  Assume
\mr
\item "{$(a)$}"  $\kappa$ is a measurable cardinal
\sn
\item "{$(b)$}"  $\kappa < \mu = \text{ cf}(\mu) < \lambda =
\text{ cf}(\lambda) = \lambda^\kappa$ and
$(\forall \alpha < \mu)(|\alpha|^{\aleph_0} < \mu)$ for simplicity.
\ermn
\ub{Then} for some c.c.c. forcing notion $P$ of cardinality $\lambda$, in
$\bold V^P$ we have: $2^{\aleph_0} = \lambda,{\frak d} = {\frak b} = \mu$ and
${\frak a} = \lambda$.
\endproclaim
\bigskip

\demo{Proof}  We choose by induction on $\zeta \le \lambda$ the following
objects satisfying the following conditions:
\mr
\item "{$(a)$}"  a sequence $\langle \bar{\Cal T}[\gamma,\zeta]:\gamma < \mu
\rangle$
\sn
\item "{$(b)$}"  $\bar{\Cal T}[\gamma,\zeta] \in 
{\frak K}^\lambda_{\text{pr}_1,\gamma}$
\sn
\item "{$(c)$}"  $\xi < \zeta \Rightarrow \bar{\Cal T}[\gamma,\xi] 
\le_{\frak K} \bar{\Cal T}[\gamma,\zeta]$
\sn
\item "{$(d)$}"  if $\zeta$ limit then $\bar{\Cal T}[\gamma,\zeta] =
\dbcu_{\xi < \zeta} \bar{\Cal T}[\gamma,\xi]$
\sn
\item "{$(e)$}"  if $\gamma < \mu,\zeta = 1$ \ub{then}
$Q_{\bar{\Cal T}}[\gamma,\zeta]$ is the dominating real forcing =
Hechler forcing
\sn
\item "{$(f)$}"  if $\gamma < \mu,\zeta = \xi +1 >1$ and $\xi$ is even, 
\ub{then} $\bar{\Cal T}[\gamma,\zeta]$ is isomorphic to 
$\bar{\Cal T}\langle \gamma +1,\xi \rangle$
over $\bar{\Cal T}[\gamma,\xi]$ say by $\bold j_{\gamma,\xi}$ where
$\bar{\Cal T} \langle \gamma +1,\xi \rangle =: 
\dsize \sum_{\beta \le \gamma} \bar{\Cal T}[\beta,\xi]$ and let 
$\hat{\bold j}_{\gamma,\xi}$ be the isomorphism induced from
$Q_{\bar{\Cal T} \langle \gamma +1,\xi \rangle}$ onto
$Q_{\bar{\Cal T}}[\gamma,\zeta]$ over $Q_{\bar{\Cal T}[\gamma,\xi]}$
\sn
\item "{$(g)$}"  if $\gamma < \mu,\zeta = \xi +1,\xi$ odd, \ub{then}
$\bar{\Cal T}[\gamma,\zeta]$ is almost isomorphic to 
$(\bar{\Cal T}[\gamma,\xi])^\kappa/D$ over $\bar{\Cal T}_{[\gamma,\xi]}$
say $\bold j_{\gamma,\xi}$ is an 
almost isomorphism from $(\bar{\Cal T}[\gamma,\xi])^\kappa/D$
onto $\bar{\Cal T}[\gamma,\zeta]$ such that by $\bold j_{\gamma,\xi} \, 
\langle x:\varepsilon < \kappa \rangle/D$ is mapped onto $x$.
\ermn
There is no problem to carry the definition.  Let $P_\zeta = 
Q_{{\bar{\Cal T}} \langle \mu,\zeta \rangle}$ where $\bar{\Cal T}
\langle \mu,\zeta \rangle =: \dsize \sum_{\gamma < \mu} \bar{\Cal T}
[\gamma,\zeta]$ for $\zeta \le \lambda,
P = P_\lambda$ and $P_{\gamma,\zeta} = Q_{{\bar{\Cal T}} \langle \gamma,
\zeta \rangle}$.  Now
\mr
\item "{$\boxtimes_1$}"  $|P| \le \lambda$ \nl
[why?  as we prove by induction on $\zeta \le \lambda$ that: each
$\bar{\Cal T}[\gamma,\zeta]$ and
$\dsize \sum_{\gamma \le \mu} \bar{\Cal T}[\gamma,\lambda]$ has
cardinality $\le \lambda$.  Hence for $\gamma < \mu$ we have: the forcing notion
$Q_{\bar{\Cal T}[\gamma,\lambda]}$ in the universe 
$\bold V^{Q_{\bar{\Cal T} \langle \gamma,\lambda \rangle}}$ has
cardinality $\le \lambda^{\aleph_0} = \lambda$]
\sn
\item "{$\boxtimes_2$}"  in $\bold V^P$ we have ${\frak b} = {\frak d} = \lambda$
\nl
[why?  let ${\underset\tilde {}\to \eta_\gamma}$ be the 
$Q_{\bar{\Cal T}[\gamma,1]}$-name of the dominating real (see clause (e)).  
As $\bar{\Cal T}[\gamma,1] \le_{\frak K} \bar{\Cal T}[\gamma,\lambda]$, 
clearly ${\underset\tilde {}\to \eta_\gamma}$ is also a
$Q_{\bar{\Cal T}[\gamma,\lambda]}$-name of a dominating real,
so $\Vdash_P ``{\underset\tilde {}\to \eta_\gamma}$ dominate
$({}^\omega \omega)^{\bold V[P_{\gamma,\lambda}]}"$.   But $\langle
P_{\gamma,\lambda}:\gamma < \mu \rangle$ is $\lessdot$-increasing with
union $P$ and cf$(\mu) = \mu > \aleph_0$ so
$\Vdash_P ``\langle
{\underset\tilde {}\to \eta_\gamma}:\gamma < \mu \rangle$ is
$<^*$-increasing and dominating", so the conclusion follows.]
\ermn
We shall prove below that ${\frak a} \ge \lambda$, together this finishes the
proof (note that it implies $2^{\aleph_0} \ge \lambda$ hence as
$\lambda = \lambda^{\aleph_0}$ by
$\boxtimes_1$ we get $2^{\aleph_0} = \lambda$)
\mr
\item "{$\boxtimes_3$}"  $\Vdash_P ``{\frak a} \ge \lambda"$.
\ermn
So assume $p \Vdash ``{\underset\tilde {}\to {\Cal A}} = 
\{{\underset\tilde {}\to A_i}:i < \theta\}$ is a MAD family,
i.e. ($\theta \ge \aleph_0$ and)
\mr
\widestnumber\item{$(iii)$}
\item "{$(i)$}"  ${\underset\tilde {}\to A_i} \in [\omega]^{\aleph_0}$,
\sn
\item "{$(ii)$}"  $i \ne j \Rightarrow |{\underset\tilde {}\to A_i} \cap
{\underset\tilde {}\to A_j}| < \aleph_0$ and
\sn
\item "{$(iii)$}"   ${\underset\tilde {}\to {\Cal A}}$ is maximal under
$(i) + (ii)$".
\ermn
Without loss of generality $\Vdash_P ``{\underset\tilde {}\to A_i} \in
[\omega]^{\aleph_0}"$. \nl
As always ${\frak a} \ge {\frak b}$, by $\boxtimes_2$ we know that 
$\theta \ge \mu$, and toward contradiction assume $\theta < \lambda$.  
For each $i < \theta$ and $m < \omega$ there is a maximal antichain 
$\langle p_{i,m,n}:n < \omega \rangle$ of $P$ and a sequence $\langle
\bold t_{i,m,n}:n < \omega \rangle$ of truth values such that $p_{i,m,n} 
\Vdash_P ``n \in {\underset\tilde {}\to A_i}$
iff $\bold t_{i,m,n}$ is truth".  We can find a countable $w_i \subseteq \mu$
such that: $\bigl[ \dbcu_{m,n < \omega} \text{ Dom}(p_{i,m,n}) \subseteq 
w_i \bigr],p_{i,m,n} \in Q_{\cup\{\bar{\Cal T}[\gamma,\lambda]:\gamma
\in w_i\}}$, moreover, 
$\gamma \in \text{ Dom}(p_{i,m,n}) \Rightarrow p_{i,m,n}(\gamma)$ is
a \nl
$Q_{\sum\{\bar{\Cal T}[\beta,\lambda]:\beta \in \gamma \cap w_i\}}$-name.  
Note that $Q_{\sum\{\bar{\Cal T}[\beta,\lambda]:\beta \in \gamma \cap w_i,
i < \theta\}} \lessdot Q_{\sum\{\bar{\Cal T}_\beta:\beta < \gamma\}}$, see
\cite{Sh:630}.
\mn
Clearly for some even $\zeta < \lambda$, we have $\{p_{i,m,n}:i < \theta,m <
\omega$ and $n < \omega \} \subseteq Q_{\sum\{\bar{\Cal T}[\beta,\zeta]:
\beta < \mu\}}$.  Now for some stationary
$S \subseteq \{\delta < \mu:\text{cf}(\delta) = \kappa\}$ and $w^*$ we have:
$\delta \in S \Rightarrow w_\delta \cap \delta = w^*$ and $\alpha <
\delta \in S \Rightarrow w_\alpha \subseteq \delta$.  Let $\langle
\delta_\varepsilon:\varepsilon < \kappa \rangle$ be an increasing sequence of
members of $S$, and $\delta^* = \dbcu_{\varepsilon < \kappa} 
\delta_\varepsilon$.  The definition of $\langle \bar{\Cal T}
[\gamma,\zeta +1]:\gamma < \mu \rangle,\langle \bar{\Cal T}
[\gamma,\zeta + 2]:\gamma < \mu \rangle$ was made to get a name of 
an infinite $\underset\tilde {}\to A \subseteq \omega$
almost disjoint to every ${\underset\tilde {}\to A_\beta}$ for $\beta <
\theta$ (in fact $(\dsize \sum_{\gamma < \mu} Q_{\bar{\Cal T}[\gamma,\zeta]})
^\kappa /D$ can be $\lessdot$-embedded into $\dsize \sum_{\gamma < \mu}
Q_{\bar{\Cal T}[\gamma,\zeta + 2]})$.  \hfill$\square_{\scite{da.6}}$\margincite{da.6}
\enddemo
\bigskip

\remark{Remark}  In later proofs in \S2 we give more details.
\endremark
\newpage

\head {\S2 \\
On Con$({\frak a} > {\frak d})$ revisited with FS, with non transitive
memory, non-well ordered length} \endhead  \resetall 
\bigskip

We first define the FSI-templates, telling us how do we iterate along a
linear order $L$; we think of having for each $t \in L$, a forcing notion
$Q_t$, say adding a generic ${\underset\tilde {}\to \nu_t}$, and $Q_t$ will
really be $\cup\{Q^{\bold V[\langle {\underset\tilde {}\to \nu_s}:s \in A
\rangle]}:A \in I_t\}$ where $I_t$ an ideal of subsets of $\{s:s <_L t\}$;
so $Q_t$ in the nice case is a definition.  In our application this definition
is constant, but we treat a more general case, so 
${\underset\tilde {}\to Q_t}$ may be defined using parameters from
$\bold V[\langle {\underset\tilde {}\to \nu_s}:s \in K_t \rangle],K_t$ a
subset of $\{s:s <_L t\}$ so the reader may consider only the case
$t \in L \Rightarrow K_t = \emptyset$.  In part (3) instead distinguishing
``$\zeta$ odd, $\zeta$ even" we can consider the two cases for each
$\zeta$.  
The depth of $L$
is the ordinal on which our induction rests (as otp$(L)$ is inadequate).
\bigskip

\definition{\stag{ad.1} Definition}  1)  An FSI-template (= finite support
iteration template) ${\frak t}$ is a sequence 
$\langle I_t:t \in L \rangle = \langle
I^{\frak t}_t:t \in L^{\frak t} \rangle = \langle I_t[{\frak t}]:t \in
L[{\frak t}] \rangle$ such that
\mr
\item "{$(a)$}"  $L$ is a linear order (but we may write $x \in {\frak t}$
instead of $x \in L$ and $x <_{\frak t} y$ instead of $x <_L y$)
\sn
\item "{$(b)$}"  $I_t$ is an ideal of subsets of $\{s:L \models s < t\}$.
\ermn
We say ${\frak t}$ is locally countable if $t \in L^{\frak t} \and (\forall
B \in [A]^{\aleph_0})(B \in I_t) \Rightarrow A \in I_t$ and we say
${\frak t}^1,{\frak t}^2$ are equivalent if $L^{{\frak t}^1} =
L^{{\frak t}^2}$ and $t \in L^{{\frak t}^1} \and |A| \le \aleph_0
\Rightarrow (A \in I^{{\frak t}^1}_t \equiv A \in I^{{\frak t}^2}_t)$. \nl
2) Let ${\frak t}$ be an FSI-template.
\mr
\item "{$(c)$}"  We say $\bar K = \langle K_t:t \in L^{\frak t} \rangle$ is
a ${\frak t}$-memory choice if
{\roster
\itemitem{ $(i)$ }  $K_t \in I^{\frak t}_t$ is countable 
\sn
\itemitem{ $(ii)$ }  $s \in K_t \Rightarrow K_s \subseteq K_t$.
\endroster}
\item "{$(d)$}"  We say $L \subseteq L^{\frak t}$ is $\bar K$-closed if
$t \in L \Rightarrow K_t \subseteq L$
\sn
\item "{$(e)$}"  for $\bar K$ a ${\frak t}$-memory choice and $L \subseteq
L^{\frak t}$ which is $\bar K$-closed we say $\bar K' = \bar K \restriction L$ if
Dom$(\bar K') = L$ and
$K'_t$ is $K_t$ for $t \in L$, (it is a 
$({\frak t} \restriction L)$-memory choice, see part (5)).
\ermn
3) For an FSI-template ${\frak t}$ and ${\frak t}$-memory choice $\bar K$
and $\bar K$-closed $L \subseteq L^{\frak t}$ we define Dp$_{\frak t}(L,
\bar K)$, the ${\frak t}$-depth (or $({\frak t},\bar K)$-depth)
of $L$ by defining by induction on the 
ordinal $\zeta$ when Dp$_{\frak t}(L,\bar K) \le \zeta$.
\mn
\ub{For $\zeta=0$}:  Dp$_{\frak t}(L,\bar K) \le \zeta$ when $L = \emptyset$.
\mn
\ub{For $\zeta$ odd}:  
Dp$_{\frak t}(L,\bar K) \le \zeta$ iff:
\mr
\item "{$(a)$}"  there is $L^*$ such that: $L^* \subseteq L,(\forall t \in L)
(\forall A \in I^{\frak t}_t)(A \cap L^* = \emptyset)$ hence $L
\backslash L^*$ is $\bar K$-closed and
Dp$_{\frak t}(L \backslash L^*,\bar K) < \zeta$ and for every $t \in L^*$ 
we have
{\roster
\itemitem{ $\boxtimes_{t,L}$ }  in $\{A \in I^{\frak t}_t:A \subseteq
L\}$ there is a maximal element \footnote{we can use less, it seems not needed
at the moment.  We can go deeper to names of depth $\le \varepsilon$
inductively on $\varepsilon < \omega_1$, as in \cite[\S3]{Sh:619}, or in a more
particular way to make the point this is used here true, and/or make
$I^{\frak t}_t$ only closed under unions (but not subsets), etc. \nl
Note that e.g. Lim$_{\frak t}(\bar Q)$ is well defined when
$L^{\frak t}$ is well ordered.} and it is $\bar K$-closed,
\endroster}
\ermn
\ub{For $\zeta > 0$ even}:  Dp$_{\frak t}(L,\bar K) \le \zeta$ iff:
\mr
\item "{$(b)$}"  there is a directed partial order $M$ and a sequence $\langle
L_a:a \in M \rangle$ with union $L$ such that $M \models a \le b 
\Rightarrow L_a \subseteq L_b$, each $L_b$ is $\bar K$-closed, 
$(\forall b \in M)(\zeta > \text{ Dp}_{\frak t}(L_b,\bar K)$) and $t \in L \and 
A \in I_t \and A \subseteq L \Rightarrow (\exists a \in M) \, A \subseteq 
L_a$.
\ermn
So Dp$_{\frak t}(L,\bar K) = \zeta$ iff Dp$_{\frak t}(L,\bar K) \ge \zeta \and
(\forall \xi < \zeta) \text{ Dp}_{\frak t}(L,\bar K) \nleq \xi$ and 
Dp$_{\frak t}(L,\bar K) = \infty$ iff ($\forall$ ordinal $\zeta$) 
[Dp$_{\frak t}(L,\bar K) \nleq \zeta]$.
\nl
4) We say $\bar K$ is a smooth ${\frak t}$-memory choice if
Dp$_{\frak t}(L^{\frak t},\bar K) < \infty$ and $\bar K$ a ${\frak t}$-memory
choice. \nl
5) If $\bar K$ is omitted we mean $K_t = \emptyset$ for $t \in L^{\frak t}$.
We say ${\frak t}$ is smooth if the trivial $\bar K$ is a smooth
${\frak t}$-memory choice.  For $L \subseteq L^{\frak t}$ let ${\frak t}
\restriction L = \langle I_t \cap {\Cal P}(L):t \in L \rangle$. \nl
6) Let $L_1 \le_{\frak t} L_2$ mean $L_1 \subseteq L_2 \subseteq
L^{\frak t}$ and $t \in L_1
\and A \in I^{\frak t}_t \Rightarrow A \cap L_2 \subseteq L_1$.
\enddefinition
\bigskip

\definition{\stag{ad.2} Definition}  Let ${\frak t} = \langle I_t:t \in
L^{\frak t} \rangle$ be a FS iteration template and $\bar K$ a
${\frak t}$-memory choice. \nl
1) We say $\bar L$ is a $({\frak t},\bar K)$-representation of $L$ if:
\mr
\item "{$(a)$}"  $L \subseteq L^{\frak t}$ is $\bar K$-closed
\sn
\item "{$(b)$}"  $\bar L = \langle L_a:a \in M \rangle$
\sn
\item "{$(c)$}"  $M$ is a directed partial order
\sn
\item "{$(d)$}"  $\bar L$ is increasing, that is $a <_M b \Rightarrow L_a \subseteq L_b$
\sn
\item "{$(e)$}"  $L = \dbcu_{a \in M} L_a$
\sn
\item "{$(f)$}"  each $L_a$ is $\bar K$-closed
\sn
\item "{$(g)$}"  if $t \in L,A \in I^{\frak t}_t,A \subseteq L$ then
$(\exists a\in M)(A \subseteq L_a)$.
\ermn
2) We say $(L^*,\bar A)$ is a $({\frak t},\bar K)-^*$representation of $L$
if
\mr
\item "{$(a)$}"  $L \subseteq L^{\frak t}$ is $\bar K$-closed
\sn
\item "{$(b)$}"  $L^* \subseteq L,\bar A = \langle A_t:t \in L^* \rangle$
\sn
\item "{$(c)$}"  if $t \in L$ and $A \in I^{\frak t}_t$ then $A \cap L^* =
\emptyset$ (so $(L \backslash L^*) \le_{\frak t} L$)
\sn
\item "{$(d)$}"  $A_t \in I^{\frak t}_t,A_t \subseteq L$ and $A_t$ is
maximal under those requirements
\sn
\item "{$(e)$}"  $L \backslash L^*$ is $\bar K$-closed (actually
follows from clause (d))
\sn
\item "{$(f)$}"  $A_t$ is $\bar K$-closed
\endroster
\enddefinition
\bigskip

\proclaim{\stag{ad.3} Claim}  Let ${\frak t}$ be an FSI-template and
$\bar K$ a ${\frak t}$-memory choice. \nl
0) The family of $\bar K$-closed sets is closed under (arbitrary) unions and
intersections.  Also if $L \subseteq L^{\frak t}$ then 
$L \cup \bigcup\{K_t:t \in L\}$ is $\bar K$-closed. \nl
1) If $L_2 \subseteq L^{\frak t}$ is $\bar K$-closed and $L_1$ is an initial
segment of $L_2$, \ub{then} $L_1$ is $\bar K$-closed. \nl
2) If $L_1 \subseteq L_2 \subseteq L^{\frak t}$ are $\bar K$-closed \ub{then}
\mr
\item "{$(\alpha)$}"  Dp$_{\frak t}(L_1,\bar K) \le 
\text{ Dp}_{\frak t}(L_2,\bar K)$, moreover
\sn
\item "{$(\beta)$}"  $(\exists t \in L_2)[L_1 \in I^{\frak t}_t]$ implies
Dp$_{\frak t}(L_1,\bar K) < \text{ Dp}_{\frak t}(L_2,\bar K)$.
\ermn
3)  If $L_1 \subseteq L_2 \subseteq L^{\frak t}$ are $\bar K$-closed
\ub{then} ${\frak t} \restriction L_2$ is an FSI-template, $L_1$ is
$({\frak t} \restriction L_2)$-closed and
Dp$_{{\frak t} \restriction L_2}(L_1,\bar K \restriction L_2) 
= \text{ Dp}_{\frak t}(L_1,\bar K)$.  ${{}}$  \hfill$\square_{\scite{ad.3}}$\margincite{ad.3}
\endproclaim
\bigskip

\demo{Proof}  0), 1) Trivial - read the definitions. \nl
2) We prove by induction on $\zeta$ that
\mr
\item "{$(*)_\zeta(\alpha)$}"  if Dp$_{\frak t}(L_2,\bar K) = \zeta$ 
(and $L_1,L_2$ are $\bar K$-closed) \ub{then}
Dp$_{\frak t}(L_1,\bar K) \le \zeta$
\sn
\item "{$(\beta)$}"  if in addition $(\exists t \in L_2)(L_1 \in I^{\frak t}
_t)$ \ub{then} Dp$_{\frak t}(L_1,\bar K) < \zeta$.
\ermn
So assume Dp$_{\frak t}(L_2,\bar K) = \zeta$, so 
Dp$_{\frak t}(L_2,\bar K) \ngeq \zeta +1$ hence one of the following cases occurs.
\enddemo
\bn
\ub{First Case}:  $\zeta = 0$.  

Trivial; note that clause $(\beta)$ is empty.
\mn
\ub{Second Case}:  $\zeta$ is odd, $L_2$ has a $({\frak t},\bar K)$-$^*$representation
$(L^*,\bar A)$ such that Dp$_{\frak t}(L_2 \backslash L^*,\bar K) <
\zeta$; see Definition \scite{ad.2}(2).

Let $L^-_2 =: L_2 \backslash L^*$; if $L_1 \subseteq L^-_2$ then by the
induction hypothesis Dp$_{\frak t}(L_1,\bar K) 
\le \text{ Dp}_{\frak t}(L^-_2,\bar K) <
\zeta$, so assume $L_1 \nsubseteq L^-_2$ and so only clause $(\alpha)$ is
relevant.  Now letting
$L^-_1 = L_1 \backslash L^*$ we have $[L^-_1,L^-_2$ are $\bar K$-closed]
$\and L^-_1 \subseteq L^-_2 \and
\text{ Dp}_{\frak t}(L^-_2,\bar K) < \zeta$ 
hence Dp$_{\frak t}(L^-_1,\bar K) < \zeta$ by the induction hypothesis.
Let $L^*_1 = L_1 \cap L^*$, so $L^*_1 \subseteq L_1,L_1$ is $\bar K$-closed,
$L_1 \backslash L^*_1 = (L_2 \backslash L^*_2) \cap L_1$ is $\bar K$-closed,
Dp$_{\frak t}(L_1 \backslash L^*_1,\bar K) = \text{Dp}_{\frak t}(L^-_1,
\bar K) < \zeta$.  Also easily: $t \in L^*_1$ implies $A_t \cap L^-_1$ is
$\bar K$-closed and maximal in $\{A \in I^{\frak t}_t:A \subseteq L_1\}$ so
$(L^*_1,\langle A_t \cap L_1:t \in L^*_1 \rangle)$ is a $({\frak t},
\bar K)$-$^*$representation of $L_1$.  So clearly
Dp$_{\frak t}(L_1,\bar K) \le \text{ Dp}_{\frak t}(L^-_1,\bar K)+1 \le
\zeta$ if Dp$_{\frak t}(L^-_1,\bar K)+1$ is odd, and Dp$_{\frak
t}(L_1,\bar K) \le (\text{Dp}_{\frak t}(L^-_1,\bar K) + 1)+1 \le
\zeta$ if Dp$_{\frak t}(L^-_1,\bar K)+1$ is even hence $< \zeta$.
\bn
\ub{Third Case}:  $\zeta$ is even and
$\langle L^a:a \in M \rangle$ is a
$({\frak t},\bar K)$-representation of $L_2$ such that $a \in M
\Rightarrow$ Dp$_{\frak t}(L_a,\bar K) < \zeta$.
\nl
Let $L^a_2 =: L^a$ and $L^a_1 =: L^a \cap L_1$, so 
$\langle L^a_1:a \in M \rangle$ is increasing, $\dbcu_{a \in M} L^a_1 = L_1$ 
and each $L^a_1$ is $\bar K$-closed (as $L^a_2,L_1$ are $\bar
K$-closed, see part (0))
and $t \in L_1 \and A \in I^{\frak t}_t \and A \subseteq L_1 \Rightarrow
(\exists a \in M)(A \subseteq L^a_2 \cap L_1 = L^a_1)$.  
Also by the induction hypothesis, $b \in M \Rightarrow 
\text{ Dp}_{\frak t}(L^b_2,\bar K) < \zeta$.  By 
the last two sentences (and Definition \scite{ad.1}) we get 
Dp$_{\frak t}(L_1,\bar K) \le \zeta$, as required in clause
$(\alpha)$.  For clause $(\beta)$ we know that there is $t \in L_2$
such that $L_1 \in I^{\frak t}_t$, hence by clause (f) of
Definition \scite{ad.2}(1)) for some $b \in M$ we have 
$L_1 \subseteq L^b$ and we can use the induction hypothesis on $\zeta$
for $L_1,L^b$. \nl
3) Easy.    \hfill$\square_{\scite{ad.3}}$\margincite{ad.3} 
\bigskip

\proclaim{\stag{ad.4} Claim}  1) If for $\ell = 1,2$ we have 
$\bar L^\ell$ is a $({\frak t},\bar K)$-representation of $L$ and
$\bar L^\ell = \langle L^\ell_a:
a \in M_\ell \rangle$ and $M = M_1 \times M_2$ \ub{then} $\bar L = \langle
L_a \cap L_b:(a,b) \in M \rangle$ is a $({\frak t},\bar K)$-representation
of $L$. \nl
2) If $(L^*_\ell,\bar A^\ell)$ is a $({\frak t},\bar K)-^*$representation of 
$L$ for $\ell = 1,2$ and we let $L^- = L \backslash L^*_1 \backslash
L^*_2$ and $\bar A = \langle A_t:t \in L^*_1 \cup L^*_2 \rangle$ and
$A_t$ is $A^\ell_t$ if $t \in L^*_\ell$ (no contradiction!)
\ub{then}
\mr
\item "{$(a)$}"  $\bar A^1 \restriction (L^*_1 \cap L^*_2) = \bar A^2
\restriction (L^*_1 \cap L^*_2)$
\sn
\item "{$(b)$}"  $(L^*_1 \cup L^*_2,\bar A)$ is a 
$(t,\bar K)$-$^*$representation of $L$.
\endroster
\endproclaim
\bigskip

\demo{Proof}  1) Straight. \nl
2) Easy, too.  \hfill$\square_{\scite{ad.4}}$\margincite{ad.4}
\enddemo
\bn
\ub{\stag{ad.4a} Discussion}:  1) 
Our next aim is to define iteration for any $\bar K$-smooth FSI-template
${\frak t}$; for this we define and prove the relevant things; of course,
by induction on the depth.
In the following Definition \scite{ad.5}, in clause (A)(a), we
avoid relying on \cite{Sh:630}; moreover the reader may consider only the case $K_t =
\emptyset$, omit ${\underset\tilde {}\to \eta_t}$ and have 
${\underset\tilde {}\to Q_{t,\bar \varphi'_t}}$ be the dominating real
forcing = Hechler forcing.  \nl
2)  We may more generally than here allow
${\underset\tilde {}\to \eta_t}$ to be e.g. a sequence of ordinals,
and member of ${\underset\tilde {}\to Q_{t,\varphi,{\underset\tilde
{}\to \eta_t}}}$ be $\subseteq {\Cal H}_{< \aleph_1}(\text{Ord})$, and
even $K_t$ large but increasing $L$, we need more ``information" from
${\underset\tilde {}\to \eta_t} \restriction \text{ Lim}_{\frak t}(\bar Q \restriction
L)$.  We may change to: ${\underset\tilde {}\to Q_t}$ is a definition of nep c.c.c. forcing
(\cite{Sh:630}) or just ``Souslin c.c.c. forcing (= snep)" or just absolute 
enough c.c.c. forcing notion.
All those cases do not make real problems (but when the parameter
${\underset\tilde {}\to \eta_t}$ have length $\ge \kappa$ it change in
the ultrapower! i.e. $\bold j({\underset\tilde {}\to \eta_t})$ has
length $>$ length of ${\underset\tilde {}\to \eta_t}$). \nl
3) If we restrict ourselves to $\sigma$-centered forcing notions
(which is quite reasonable), we may in Definition \scite{ad.1}(3)(a)
omit $\boxtimes_{t,L}$ if in Definition \scite{ad.5} below in (A)(b)
second case we add that $t \in L^* \Rightarrow p \restriction (L
\backslash L^*)$ forces a value to ${\underset\tilde {}\to
f_t}(p(t))$ where ${\underset\tilde {}\to f_t}:{\underset\tilde {}\to
Q_t} \rightarrow w$ witnessed $\sigma$-centerness and is absolute
enough (or just assume $Q_t \subseteq \omega \times Q'_t,f_t(p(t)$ is
the first coordinate).  More carefully we can do this with $\sigma$-linked
instead $\sigma$-centered.
\bigskip

\definition{\stag{ad.5} Definition/Claim}  Let ${\frak t}$ be an FSI-template
and $\bar K = \langle K_t:t \in L^{\frak t} \rangle$ be a smooth
${\frak t}$-memory choice.

By induction on the ordinal $\zeta$ we shall define and prove
\mr
\item "{$(A)$}"  [Def] $\quad$ for $L \subseteq L^{\frak t}$ which is
$\bar K$-closed of $({\frak t},\bar K)$-depth $\le \zeta$ we define 
{\roster
\itemitem{ $(a)$ }  when $\bar Q = \langle
{\underset\tilde {}\to Q_{t,{\bar \varphi}_t,
{\underset\tilde {}\to \eta_t}}}:t \in L \rangle$
is a $({\frak t},\bar K)$-iteration of def-c.c.c. forcing notions, but we
can let ${\underset\tilde {}\to \eta_t}$ code $\bar \varphi_t$ so
usually omit $\bar \varphi_t$
\sn
\itemitem{ $(b)$ }  Lim$_{\frak t}(\bar Q)$ for $\bar Q$ as in (A)(a)  
\endroster}
\item "{$(B)$}"  [Claim] $\quad$ for $L_1 \subseteq L_2 \subseteq
L^{\frak t}$ which are $\bar K$-closed of $({\frak t},\bar K)$-depth 
$\le \zeta$ 
\nl

$\qquad \qquad$ and $({\frak t},\bar K)$-iteration of def-c.c.c. forcing
notions \nl

$\qquad \qquad \bar Q = \langle 
{\underset\tilde {}\to Q_{t,{\underset\tilde {}\to {\bar \varphi}_t}}},
{\underset\tilde {}\to \eta_t}:t \in L_2 \rangle$ and letting $\bar
Q^1 = \bar Q^2 \restriction L_1$ we prove:
{\roster
\itemitem{ $(a)$ }  $\bar Q \restriction L_1$ is a
$({\frak t},\bar K \restriction L_1)$-iteration of def-c.c.c. 
forcing notions
\sn
\itemitem{ $(b)$ }  Lim$_{\frak t}(\bar Q \restriction L_1) \subseteq
\text{ Lim}_{\frak t}(\bar Q)$ as quasi orders
\sn
\itemitem{ $(c)$ }  if $L_1 \le_{\frak t} L_2$ (see Definition \scite{ad.1}(6)) and
$p \in \text{ Lim}_{\frak t}(\bar Q^2)$, \ub{then} 
$p \restriction L_1 \in \text{ Lim}_{\frak t}(\bar Q^1)$ and
Lim$_{\frak t}(\bar Q^2) \models ``p \restriction L_1 \le p"$
\sn
\itemitem{ $(d)$ }  if $L_1 \le_{\frak t} L_2$ and
$p \in \text{ Lim}_{\frak t}(\bar Q)$ and
$\text{Lim}_{\frak t}(\bar Q^1 \restriction L_1) \models ``(p \restriction L_1) 
\le q"$ \ub{then} $q \cup (p \restriction (L_2 \backslash
L_1))$ is a lub of $\{p,q\}$ in $\text{Lim}_{\frak t}(\bar Q^2)$; hence
Lim$_{\frak t}(\bar Q \restriction L_1) \lessdot \text{ Lim}_{\frak t}
(\bar Q)$
\sn
\itemitem{ $(e)$ }  Lim$_{\frak t}(\bar Q \restriction L_1) \lessdot 
\text{ Lim}_{\frak t}(\bar Q)$, that \footnote{here we do not assume
$L_1 \le_{\frak t} L_2$,} is
\sn
\itemitem{ ${{}}$ } $\qquad (i) \quad p \in \text{ Lim}_{\frak t}(\bar
Q \restriction L_1)
\Rightarrow p \in \text{ Lim}_{\frak t}(\bar Q)$
\sn
\itemitem{ ${{}}$ }  $\qquad (ii) \quad \text{ Lim}_{\frak t}(\bar Q
\restriction L_1)
\models p \le q \Rightarrow \text{ Lim}_{\frak t}(\bar Q) \models p \le q$
\sn
\itemitem{ ${{}}$ }  $\qquad (iii) \quad$ if ${\Cal I} \subseteq
\text{ Lim}_{\frak t}(\bar Q \restriction L_1)$ 
is predense in $\text{Lim}_{\frak t}(\bar Q \restriction L_1)$,
\ub{then} \nl

$\qquad \qquad \qquad$ ${\Cal I}$ is predense in 
Lim$(\bar Q)$ (hence if $p,q \in \text{ Lim}_{\frak t}(\bar Q)$ are \nl

$\qquad \qquad \qquad$ incompatible in 
Lim$_{\frak t}(\bar Q \restriction L_1)$ then they are incompatible \nl

$\qquad \qquad \qquad$ in Lim$(\bar Q)$)
\sn
\itemitem{ $(f)$ }  if $L_0 \subseteq L_2$ is $\bar K$-closed,
$L = L_0 \cap L_1$ and $p \in
\text{ Lim}_{\frak t}(\bar Q \restriction L_0)$ and $q \in \text{ Lim}
_{\frak t}(\bar Q \restriction L)$ satisfies \nl
$(\forall r \in \text{ Lim}_{\frak t}(\bar Q \restriction L))
[q \le r \rightarrow p,r$ are compatible in Lim$_{\frak t}(\bar Q
\restriction L_0)]$ \ub{then} \nl
$(\forall r \in \text{ Lim}_{\frak t}(\bar Q \restriction L_1))
[q \le r \rightarrow p,r$ are compatible in Lim$_{\frak t}(\bar Q
\restriction L_2)]$ \nl
\block 
\eightpoint [explanation: this means that if $q$ forces for
$\Vdash_{\text{Lim}_{\frak t}(\bar Q \restriction L_0)}$ that
$p \in \text{ Lim}_{\frak t}(\bar Q \restriction
L_0)/\text{Lim}_{\frak t}(\bar Q \restriction L)$ \ub{then} $q$ forces for 
$\Vdash_{\text{Lim}_{\frak t}(\bar Q \restriction L_1)}$ that 
$p \in \text{ Lim}_{\frak t}(\bar Q)/\text{Lim}_{\frak t}(\bar Q \restriction L_1)$.]
\endblock
\tenpoint
\sn
\itemitem{ $(g)$ }  if $\langle L_a:a \in M_1 \rangle$ is a
$({\frak t},\bar K)-$representation of $L_1$ \ub{then} 
Lim$_{\frak t}(\bar Q \restriction L_1) =
\dbcu_{a \in M_1} \text{ Lim}_{\frak t}(\bar Q \restriction L^1_a)$
\sn
\itemitem{ $(h)$ }   if 
$(L^*,\bar A)$ is a $({\frak t},\bar K)$-$^*$representation
of $L_1$, \ub{then} Lim$_{\frak t}(\bar Q \restriction L_1)$ is as defined in
(A)(b) of our definition below, second case, from $(L^*,\bar A)$
\sn
\itemitem{ $(i)$ }  $(\alpha) \quad $if $p_1,p_2 \in \text{ Lim}_{\frak t}(\bar Q)$ and
$t \in \text{ Dom}(p_1) \cap \text{ Dom}(p_2) \Rightarrow p_1(t) =$ \nl

\hskip32pt  $p_2(t)$, \ub{then} $q=p_1 \cup p_2$ (i.e. $p_1 \cup (p_2 \backslash
(\text{ Dom}(p_1)))$ belongs to \nl

\hskip32pt  Lim$_{\frak t}(\bar Q)$ and is a l.u.b. of $p_1,p_2$
\sn
\itemitem{ ${{}}$ }  $(\beta) \quad p \in \text{ Lim}_{\frak t}(\bar Q)$
\ub{iff} $p$ is a function with domain a finite subset of \nl

\hskip32pt $L_2$ such
that for every $t \in \text{ Dom}(p)$ for some $A \in I^{\frak t}_t,A$
is $\bar K$-closed \nl

\hskip32pt  and $K_t \subseteq A$ and
$\Vdash_{\text{Lim}_{\frak t}(\bar Q \restriction A)} ``p(t) \in
Q_{t,{\underset\tilde {}\to \eta_t}}"$
\sn
\itemitem{ ${{}}$ } $(\gamma) \quad \text{Lim}_t(\bar Q) \models p \le q$
\ub{iff} $p,q \in \text{ Lim}_t(\bar Q)$ and for every $t \in \text{
Dom}(p)$ \nl

\hskip32pt   we have $t \in \text{ Dom}(q)$ and for some $\bar K$-closed
$A \in I^{\frak t}_t$ we have\nl

\hskip32pt   $q \restriction A \in \text{ Lim}_{\frak
t}(\bar Q \restriction A)$ and $q \Vdash_{\text{Lim}_{\frak t}(\bar Q
\restriction A)} ``p(t) \le q(t)$ \nl

\hskip32pt   in $Q_{t,{\underset\tilde {}\to
\eta_t}}$ (as interpreted in $\bold V^{\text{Lim}_{\frak t}(\bar Q
\restriction A)}$ of course)"
\sn
\itemitem{ $(j)$ }  Lim$_{\frak t}(\bar Q)$ is a c.c.c. forcing notion
and Lim$_{\frak t}(\bar Q) = \cup\{\text{Lim}_{\frak t}(\bar Q
\restriction L):L \in [L_2]^{\le \aleph_0}\}$
\sn
\itemitem{ $(k)$ }  Lim$_{\frak t}(\bar Q)$ has cardinality $\le
|L_2|^{\aleph_0}$ (here we use the assumption that ${\underset\tilde
{}\to \eta_t}$ and members of 
${\underset\tilde {}\to Q_{t,{\underset\tilde {}\to \eta_t}}}$ 
are reals; see definition in
(A)(a)(i)+(i)) below).
\endroster}
\endroster
\enddefinition
\mn
Let us carry the induction. 
\sn
\ub{Part (A)}: [Definition]

So assume Dp$_{\frak t}(L,\bar K) \le \zeta$.  If Dp$_{\frak t}(L) < \zeta$ we have
already defined being $({\frak t},\bar K)$-iteration and Lim$_{\frak t}(\bar Q
\restriction L)$, so assume Dp$_{\frak t}(L) = \zeta$.
\mn
\ub{Clause(A)(a)}: 
\mr
\widestnumber\item{$(iii)$}
\item "{$(i)$}"  ${\underset\tilde {}\to \eta_t}$ is a Lim$_{\frak t}(\bar Q
\restriction K_t)$-name of a real (i.e. from ${}^\omega 2$, used as a
parameter)  \nl
(legal as $K_t \subseteq L \and K_t \in I_t \and t \in L$ hence by
\scite{ad.3}(2), clause $(\beta)$ we have
Dp$_{\frak t}(K_t,\bar K) < \text{ Dp}_{\frak t}(K_t \cup \{t\},\bar K) \le
\text{ Dp}_{\frak t}(L,\bar K) \le \zeta$ so
Lim$_{\frak t}(\bar Q \restriction L_t)$ is a well defined forcing notion by the
induction hypothesis and \scite{ad.3}(2), clause $(\beta)$)
\sn
\item "{$(ii)$}"  $\bar \varphi_t$ is a pair of formulas with the
parameters ${\underset\tilde {}\to \eta_t}$ defining in
$\bold V^{\text{Lim}_{\frak t}(\bar Q \restriction K_t)}$ a forcing notion
denoted by $Q_{t,\bar \varphi_t,{\underset\tilde {}\to \eta_t}}$ whose
set of elements is $\subseteq {\Cal H}(\aleph_1)$
\sn
\item "{$(iii)$}"  in $\bold V^{\text{Lim}_{\frak t}(\bar Q
\restriction K_t)}$, if $P' \lessdot P''$ are c.c.c. forcing notions \ub{then}
$Q_{t,\bar \varphi_t,\eta_t}$ as interpreted
in $(\bold V^{\text{Lim}_{\frak t}(\bar Q \restriction K_t)})^{P'}$ 
is a c.c.c. forcing notion there, and it is a $\lessdot$-subforcing
of $(P''/P') * {\underset\tilde {}\to Q_{t,\bar \varphi,\eta_t}}$ where
$Q_{t,\bar \varphi,\eta_t}$ mean as interpreted
in $(\bold V^{\text{Lim}_{\frak t}(\bar Q \restriction K_t)})^{P''}$ (i.e. $``p
\le q"$,``$p,q$ incompatible", ``$\langle p_n:n < \omega \rangle$ is
predense" (so the sequence is from the smaller universe) are preserved).
\endroster
\mn
\ub{Clause (A)(b)}:
\sn
\ub{First Case}:  $\zeta =0$.

Trivial
\mn
\ub{Second Case}:  $\zeta > 0$ odd.

So let $(L^*,\bar A)$ be a $({\frak t},\bar K)$-$^*$representation of $L$.
\nl
Define

$$
\align
p \in \text{ Lim}_{\frak t}(\bar Q) \text{ \ub{iff} } p:&p 
\text{ is a finite function, Dom}
(p) \subseteq L,p \restriction (L \backslash L^*) \in 
\text{ Lim}_{\frak t}(\bar Q \restriction (L \backslash L^*)) \\
  &\text{and if } t \in L^* \cap \text{ Dom}(p), \text{ then } p(t) 
\text{ is a Lim}_{\frak t}(\bar Q \restriction A_t) \text{-name} \\
  &\text{of a member of } Q_{t,{\bar \varphi}_t,
{\underset\tilde {}\to \eta_t}} 
\endalign
$$
\mn
and the order is \nl
Lim$_{\frak t}(\bar Q) \models p \le q$ \ub{iff}
\mr
\item "{$(i)$}"  Lim$_{\frak t}(\bar Q \restriction (L \backslash L^*)) \models
``(p \restriction (L \backslash L^*) \le (q \restriction (L \backslash L^*))"$ 
and
\sn
\item "{$(ii)$}"  if $t \in L^* \cup \text{ Dom}(p)$ then $q \restriction A_t
\Vdash_{\text{Lim}_{\frak t}(\bar Q \restriction A_t)} ``p(t) \le q(t)"$.
\ermn
Clearly Lim$_{\frak t}(\bar Q)$ is a quasi order.
But we should prove that Lim$_{\frak t}(\bar Q)$ is well defined,
which means that the definition does not depend on the representation.  So we prove
\mr
\item "{$\boxtimes_1$}"  if Dp$_{\frak t}(L,\bar K) = \zeta$ and for
$\ell=1,2$ we have $(L^*_\ell,\bar A^\ell)$ is a 
$({\frak t},\bar K)$-$^*$representation of $L$ with 
Dp$_{\frak t}(L \backslash L^*_\ell,\bar K) < 
\zeta$ and $Q^\ell$ is Lim$_{\frak t}(\bar Q \restriction L)$
as defined by $(L^*_\ell,\bar A^\ell)$ above, \ub{then} $Q^1 = Q^2$.
\ermn
This is immediate by Claim \scite{ad.4}(2) and 
the induction hypothesis clause (B)(h).
\mn
\ub{Third Case}:  $\zeta$ even $> 0$.

So there are a directed partial order $M$ and $\bar L = \langle L_a:a \in M \rangle$ a 
$({\frak t},\bar K)$-representation of $L$ such that
$a \in M \Rightarrow \text{ Dp}_{\frak t}(L_a,\bar K) < \zeta$.  By
the induction hypothesis, $a \le_M b \Rightarrow L_a \subseteq L_b \and
\text{ Lim}_{\frak t}(\bar Q \restriction L_a) \subseteq \text{
Lim}_{\frak t}(\bar Q \restriction L_b)$.

We let Lim$_{\frak t}(\bar Q \restriction L) = \dbcu_{a \in M}
\text{ Lim}_{\frak t}(\bar Q \restriction L_a)$, so we have to prove
\mr
\item "{$\boxtimes_2$}"  the choice is of $\bar L$ is immaterial.
\ermn
So we just assume that for $\ell =1,2$ we have: $M_\ell$ is a directed
partial order,
$\bar L^\ell = \langle L^\ell_a:a \in M_\ell \rangle,L^\ell_a \subseteq L,
M_\ell \models a \le b \Rightarrow L^\ell_a \subseteq L^\ell_b$ and
$(\forall t \in L)(\forall A \in I_t)[A \subseteq L \rightarrow (\exists a
\in M_\ell)(A \subseteq L^\ell_a)$ and Dp$_{\frak t}(L^\ell_a,\bar K) < 
\zeta$. \nl
We should prove that $\dbcu_{a \in M_1} \text{ Lim}_{\frak t}(\bar Q
\restriction L^1_a),\dbcu_{a \in M_2} \text{ Lim}_{\frak t}(\bar Q
\restriction L^2_a)$ are equal, as quasi orders of course.

Now $M =: M_1 \times M_2$ with $(a_1,a_2) \le (b_1,b_2) \Leftrightarrow a_1
\le_{M_1} b_1 \and a_2 \le_{M_2} b_2$, is a directed partial order.  We let
$L_{(a_1,a_2)} = L^1_{a_1} \cap L^2_{a_2}$, so clearly $L_{(a_1,a_2)}
\subseteq L^{\frak t}$, Dp$_{\frak t}(L_{(a_1,a_2)},\bar K) 
< \zeta$ and $(a_1,a_2) \le_M
(b_1,b_2) \Rightarrow L_{(a_1,a_2)} \subseteq L_{(b_1,b_2)}$ and
$\langle L_{(a_1,a_2)}:(a_1,a_2) \in M \rangle$ is a $({\frak t},
\bar K)$-representation of $L$ by Claim \scite{ad.4}(1). So by
transitivity of equality, it is enough to prove for $\ell =1,2$ that
$\dbcu_{a \in M_\ell} \text{ Lim}_{\frak t}(\bar Q \restriction L^\ell_a),
\dbcu_{(a,b) \in M} \text{ Lim}_{\frak t}(\bar Q \restriction
L_{(a,b)})$ are equal as quasi orders.  
By the symmetry in the situation without loss of generality $\ell=1$.
\sn
Now for every $a \in M_1,\bar L= \langle L_{(a,b)}:b \in M_2 \rangle$ 
satisfies: $L^1_a \subseteq L$, Dp$(L^1_a) < \zeta,L_{(a,b)} \subseteq
L^1_a,L^1_a = \dbcu_{b \in M_2} L_{(a,b)},b_1 \le_{M_2} b_2 \Rightarrow
L_{(a,b_1)} \subseteq L_{(a,b_2)}$.  Also we know that $(\forall t \in
L)(\forall A \in I^{\frak t}_t)(\exists b \in M_2)(A \subseteq L
\rightarrow A \subseteq L_b)$ hence
$(\forall t \in L^1_a)(\forall A \in I^{\frak t}_t)
(A \subseteq L^1_a \rightarrow (\exists b \in M_2)(A \subseteq L_{(a,b)}))$.
Hence by the induction hypothesis for clause (B)(g) we have
Lim$_{\frak t}(\bar Q \restriction L^1_a),\dbcu_{b \in L_2}
\text{ Lim}_{\frak t}(\bar Q \restriction L_{(a,b)})$ are equal as quasi
orders.  As this holds for every $a \in M_1$ and $M_1$ is directed
we get $\dbcu_{a \in M_1}$ Lim$_{\frak t}
(\bar Q \restriction L^1_a),\dbcu_{a \in M_1} \, \dbcu_{b \in M_2}
\text{ Lim}_{\frak t}(\bar Q \restriction L_{(a,b)})$ are equal as quasi orders.  But 
the second is equal to $\dbcu_{(a,b) \in M} \text{ Lim}_{\frak t}(\bar Q
\restriction L_{(a,b)})$ so we are done.
\bn
\ub{Part (B)}:
\sn
\ub{\bf First Case}:  $\zeta=0$.

Trivial.
\mn
\ub{\bf Second Case}: $\zeta > 0$ is odd.

So by Definition \scite{ad.1}(3) there are $L^*,\langle A_t:t \in L^*
\rangle$ as in the appropriate case there.  Let $\langle t^*_i:i <
i(*) \rangle$ list $L^*$ with no repetitions.  So
$\{A \in I^{\frak t}_{t^*_i}:A \subseteq L_2\}$ has the maximal member
$A_{t^*_i}$ and $t \in L^* \Rightarrow t \notin A_{t^*_i}$ and 
Dp$_{\frak t}(L^-_2,\bar K) \le \xi =: \zeta -1$ 
where $L^-_2 =: L_2 \backslash L^* = L_2 \backslash \{t^*_i:i < i^*\}$, so
$i < i(*) \Rightarrow A_{t^*_i} \subseteq L^-_2$, that is
$(L^*_2,\bar A^1)$ is a $({\frak t},\bar K)$-$^*$representation of $L_2$ where
$L^*_2 = L_2 \backslash L^-_2,\bar A^2 = \langle A_{t^*_i}:i < i(*) \rangle$.
So we have already defined Lim$_{\frak t}(\bar Q)$.  We shall use freely the
uniqueness in the second case in the definition (A)(b).
Let $L^*_1 = L^*_2 \cap L_1,L^-_1 = L_1 \cap L^-_2$ and $\bar A^1 = \langle
A^1_t:t \in L^*_1 \rangle$, with $A^1_2 = A^2_t \cap L_1$ and $A^2_t = A_t$.
\mn
\ub{Clause (B)(a)}: \nl

Easy.
\mn
\ub{Clause (B)(b)}: \nl

If $L_1 = L^-_2$ this follows by the definition of Lim$_{\frak t}(\bar
Q \restriction L_2)$.

If $L_1 \subseteq L^-_2$ this is very easy by the induction hypothesis
and the previous sentence.  
Otherwise, clearly $(L^*_1,\bar A^1)$ is a $({\frak t},\bar K)$-$^*$representation of
$L_2$ so by clause $(B)(h)$ when Dp$_{\frak t}(L_1,\bar K) < \zeta$, and by
uniqueness proved in part (A) otherwise, 
we have: Lim$_{\frak t}(\bar Q \restriction L_1)$ is
defined as in (A)(b) second case for $(L^* \cap L_1,\bar A \restriction L_1)$.
By the induction hypthesis, Lim$_{\frak t}(\bar Q \restriction L^-_1) 
\lessdot \text{ Lim}_{\frak t}(\bar Q \restriction L^-_2)$ 
hence for each $i < i(*)$,
\mr
\item "{$(*)$}"  the forcing notion 
${\underset\tilde {}\to Q_{t^*_i,\bar \varphi^{t^*_i}}},
{\underset\tilde {}\to \eta_{t^*_i}}$
as interpreted in $\bold V^{\text{Lim}_{\frak t}(\bar Q \restriction
L^-_1)}$ is a sub-quasi order of the same forcing notion interpreted in
$\bold V^{\text{Lim}_{\frak t}(\bar Q \restriction L^-_2)}$.
\ermn
Looking at the definitions of Lim$_{\frak t}(\bar Q \restriction L_2)$, Lim$_{\frak
t}(\bar Q \restriction L_1)$ using $(L^*_2,\bar A^2)$ and $(L^*_1,\bar
A^1)$, O.K. by the uniqueness we can easily finish. 
\mn
\ub{Clause (B)(c),(d)}: \nl

Straight.
\mn
\ub{Clause (e)}:  \nl

If $L_1 \subseteq L^-_2$ then Lim$_{\frak t}(\bar Q \restriction L_1)
\lessdot \text{ Lim}_{\frak t}(\bar Q \restriction L^-_2)$ by the induction
hypothesis and Lim$_{\frak t}(\bar Q \restriction L^-_2) \lessdot
\text{ Lim}_{\frak t}(\bar Q \restriction L_2)$ by the definition in part (A)
so we are done.

So assume that $L_1 \nsubseteq L^-_2$, so $(L^*_1,\bar A^2)$ is a
$({\frak t},\bar K)$-$^*$representation of $L_1$ so the definition in clause
(A)(b) second case apply.  Consider Lim$_{\frak t}(\bar Q \restriction
L^-_1) \lessdot \text{ Lim}_{\frak t}(\bar Q \restriction L^-_2)$ which hold
by the induction hypothesis and the definitions of Lim$_{\frak t}(\bar Q
\restriction L^-_\ell)$ according to the $({\frak t},\bar
K)-^*$representation 
$(L^*_\ell,\bar A^\ell)$.

Now in $(*)$ above we can add
\mr
\item "{$(*)^+$}"  if $j(1) \le j(2) \le i,L_{1,j} = L^-_1 \cup
\{t^*_\varepsilon:\varepsilon < j,t^*_\varepsilon \in L_1\},L_{2,j} =
L^-_2 \cup\{t^*_\varepsilon:\varepsilon < j\},
{\Cal I} \in \bold V^{\text{Lim}_{\frak t}(\bar Q
\restriction L^-_1)}$ is a predense subset of ${\underset\tilde {}\to
Q_{t^*_i,\bar \varphi^{t^*_i},{\underset\tilde {}\to \eta_{t^*_i}}}}$
as interpreted in $\bold V^{\text{Lim}_{\frak t}(\bar Q \restriction
L_{1,j(1)})}$, \ub{then} ${\Cal I}$ is also a predense subset of $Q_{t^*_i,\bar
\varphi_{t^*_i},{\underset\tilde {}\to \eta_{t^*_i}}}$ as interpreted
in $\bold V^{\text{Lim}_{\frak t}(\bar Q \restriction L_{2,j(2)})}$ is.
\ermn
So the conclusion is immediate.
\mn
\ub{Clause (B)(f)}: \nl

Let $L_0,L = L_1 \cap L,q \in \text{ Lim}_{\frak t}(\bar Q
\restriction L), p \in \text{ Lim}_{\frak t}(\bar Q \restriction
L_0)$ be as there; clearly $L_0 \cup L_1$ is $\bar K$-closed.

If Dp$_{\frak t}(L_0 \cup L_1,\bar K) < \zeta$ we can use the induction
hypothesis and clause (B)(e) which we have already proved;
so assume that this fails, so Dp$_{\frak t}(L \cup L_1,\bar K) =
\zeta$ and so let $(L^*,\bar A)$ witness this.  Now using $\boxtimes_1$ and
the induction hypothesis for clause (B)(h) we can prove it by
induction on $i(*)$ thus reducing it to the case 
$i(*) +1,t^*_i \in L$ which is easy using $(*)^+$ from
above so we are done.
\mn
\ub{Clause (B)(g)}: \nl

Again using $\boxtimes_1$ and the induction hypothesis for clause (B)(h).
\mn
\ub{Clause (B)(h)}: \nl

Straight.
\bn
\ub{Clause (B)(i)}: \nl

Easy.
\mn
\ub{Clause (B)(j)}: \nl

Let $p_\alpha \in \text{ Lim}_{\frak t}(\bar Q)$ for $\alpha <
\omega_1$; set $w_\alpha =: \{i:t^*_i \in \text{ Dom}(p_\alpha)\}$, so
without loss of generality $\langle w_\alpha:\alpha < \omega_1 \rangle$ form a
$\Delta$-system with heart $w$; let $p'_\alpha = p_\alpha
\restriction (L^-_2 \cup w)$, and easily it suffices to prove that for
some $\alpha \ne \beta,p'_\alpha,p'_\beta$ are compatible in
Lim$_{\frak t}(\bar Q \restriction (L^-_2 \cup w))$  (if $q$ is a
common upper bound of $p'_\alpha,p'_\beta$ in Lim$_{\frak t}(\bar Q
\restriction (L^-_2 \cup w))$, then $q^+ = q \cup (p'_\alpha
\restriction (L_2 \backslash L^-_2 \backslash w)) \cup (p'_\beta
\restriction (L_2 \backslash L^-_2 \backslash w))$ is as required by
clauses $(\ell),(m)$ which is said below easily holds).  
We can do this by induction on $|w|$
and (using the uniqueness proved in (A)(b) above) we can reduce this
to the case $w$ is a singleton, say $\{t^*_0\}$.  So $p^-_\alpha =
p'_\alpha \restriction L^-_2 \in \text{ Lim}_{\frak t}(\bar Q
\restriction L^-_2)$ for $\alpha < \omega_1$ hence for some $G_2 \subseteq
\text{ Lim}_{\frak t}(\bar Q \restriction L^-_2)$ generic over $\bold
V$, the set $u = \{\alpha < \omega_1:p^-_\alpha \in G_2\}$ is
uncountable; now as Lim$_{\frak t}
(\bar Q \restriction A_{t^*_0}) \lessdot \text{ Lim}_{\frak t}(\bar
Q \restriction L^-_2)$, clearly $G^* = G_2 \cap \text{ Lim}_{\frak
t}(\bar Q \restriction A_{t^*_0})$ is generic over $\bold V$ and
$\alpha \in u \Rightarrow p^-_\alpha(t^*_0)[G^*] \in 
{\underset\tilde {}\to Q_{t^*_0,\bar \varphi_{t^*_0},{\underset\tilde {}\to
\eta^*_0}}}[\bold V[G^*]] \subseteq 
{\underset\tilde {}\to Q_{t^*_0,\bar \varphi_{t^*_0},{\underset\tilde {}\to
\eta_{t^*_0}}}}[Q][\bold V[G_2]]$.

Hence by (A)(a)(iii) below for some $\alpha \ne \beta$ from
$u,p^-_\alpha(t^*_0)[G^*],p^-_\alpha(t^*_0)[G^*]$ are compatible in
${\underset\tilde {}\to Q_{t^*_0,\bar \varphi_{t^*_0},{\underset\tilde {}\to
\eta_{t^*_0}}[G^*]}}[\bold V[G_2]]$, hence in
${\underset\tilde {}\to Q_{t^*_0,\bar \varphi_{t^*_0},{\underset\tilde {}\to
\eta_{t^*_0}}[G^*]}}[\bold V[G^*]]$, and we can easily finish.
\mn
\ub{Clause (k),(l),(m)}:

Easy. 
\bn
\ub{\bf Third Case}:  $\zeta >0$ even.

So let $\langle L^2_a:a \in M \rangle$ be a $({\frak t},
\bar K)$-representation of $L_2$ with $a \in M \Rightarrow 
\text{ Dp}_{\frak t}(L_a,\bar K) < \zeta$.
\mn
\ub{Clause (B)(a)}: \nl

Trivial.
\mn
\ub{Clause (B)(b)}: \nl

Clearly Dp$_{\frak t}(L_2,\bar K) \le \zeta$ by Claim
\scite{ad.3}(2)$(\alpha)$ hence Lim$_{\frak t}(\bar Q \restriction
L_1)$ is well defined by (A)(b) above Lim$_{\frak t}(\bar Q) = 
\text{ Lim}_{\frak t}(\bar Q \restriction L_2) = \dbcu_{a \in M_2}$
Lim$_{\frak t}(\bar Q \restriction L^2_a)$ as quasi orders. \nl
Clearly $\langle L^1_a = L_1 \cap L^2_a:a \in M  
\rangle$ is a $({\frak t},\bar K)$-representation of $L_1$ hence 
by the induction hypothesis (if Dp$_{\frak t}(L_1,\bar K) < \zeta$) or
by the uniqueness proved in (A)(b) (if Dp$_{\frak t}(L_1,\bar K) = \zeta$) we know that
Lim$_{\frak t}(\bar Q \restriction L_1) = 
\dbcu_{a \in M}$ Lim$_{\frak t}(\bar Q \restriction
L^1_a)$ as quasi orders and by the induction hypothesis for (B)(b) we know
Lim$_{\frak t}(\bar Q \restriction L^1_a) \subseteq$ Lim$_{\frak t}(\bar Q
\restriction L^2_a)$ as quasi orders (for $a \in M$), and we can easily finish.
\mn
\ub{Clause (B)(c),(d)}: \nl

Use the proof of clause (B)(b) noting that $L^1_a \le_{\frak t} L^2_a$
and so we can use the induction hypothesis (i.e. if $p \in \text{
Lim}_{\frak t}(\bar Q \restriction L_2)$, as $M$ is directed there is
$a \in M$ such that Dom$(p) \subseteq L^2_a$, now $a \le_M b
\Rightarrow p \restriction L^1_b = p \restriction L^1_a$ and we can
finish easily).
\mn
\ub{Clause (B)(e)}: \nl

The statements (i) + (ii) holds by clause (b). \nl
The statement (iii) holds: let ${\Cal I}$ be a predense subset of
Lim$_{\frak t}(\bar Q \restriction L_1)$, let $p \in \text{
Lim}_{\frak t}(\bar Q)$, so for some $a \in M$ we have $p \in \text{
Lim}_{\frak t}(\bar Q \restriction L^2_a)$.  By the induction
hypothesis applying clause (B)(e) to $p,L^1_a,L^2_a$ there is $q \in
\text{ Lim}_{\frak t}(\bar Q \restriction L^1_a)$ as there.  Now by the
assumption on ``${\Cal I} \subseteq \text{ Lim}_{\frak t}(\bar Q
\restriction L_1)$ is dense", as $q \in \text{ Lim}_{\frak t}(\bar Q
\restriction L_1)$ (by clause (B)(b)) we can find $q_0 \in {\Cal I}$
and $q_1$ such that Lim$_{\frak t}(\bar Q \restriction L_1) \models
q_0 \le q_1 \and q \le q_1$, so for some $b \in M$ we have
$q,q_0,q_1 \in L^1_b$ and $a \le_M b$ (as $M$ is directed).  Now we
consider $p,q,L^1_a,L^2_a,L^1_b,L^2_b$ and apply by clause (B)(f).
\mn
\ub{Clause (B)(f)}: \nl

Easy to check using clause (f) for the $L^2_a$'s, which holds by the
induction hypothesis.
\mn
\ub{Clause (B)(g)}: \nl

Let $M_2 =: M$.
For each $a_1 \in M_1$, clearly Dp$_{\frak t}(L_a,\bar K) \le \zeta$ as
$L_{a_1} \subseteq L_1$ and $\langle L_{a_1} \cap L^2_{a_2}:a_2 \in M_2
\rangle$ is a $({\frak t},\bar K)$-representation of $L_a$ hence by (A)(b)
we know Lim$_{\frak t}(\bar Q \restriction L_{a_1}) = \dbcu_{a_2 \in M_2}
\text{ Lim}_{\frak t}(\bar Q \restriction (L_{a_1} \cap L^2_{a_2}))$.
The rest should be clear.
\mn
\ub{Clause (B)(h)}: \nl

Easy.
\mn
\ub{Clause (B)(i)}:

Easy.
\mn
\ub{Clause (B)(j)}:

So let $p_\alpha \in \text{ Lim}_{\frak t}(\bar Q)$ for $\alpha <
\omega_1$; let $w_\alpha = \text{ Dom}(p_\alpha)$ and without loss of generality
$\langle w_\alpha:\alpha < \omega_1 \rangle$ is a $\Delta$-system with heart
$w$.  So for some $a \in M$ we have $w \subseteq L^2_a$.  
For each $\alpha$, for some $a_\alpha \in M$ we have $a \le_M
a_\alpha$ and $p_\alpha \in \text{ Lim}_{\frak t}(\bar Q \restriction
L^2_{a_\alpha})$, so by clause (e) for $L^2_{a_\alpha},L^2_a$ (which
holds by the induction hypothesis, there is $p^+_\alpha \in \text{
Lim}_{\frak t}(\bar Q \restriction L^2_a)$ such that $p^*_\alpha
\Vdash_{\text{Lim}_{\frak t}(\bar Q \restriction L_a)} ``p \in \text{
Lim}_{\frak t}(\bar Q \restriction L^2_{a_\alpha})/\text{Lim}_{\frak
t}(\bar Q \restriction L^2_a)"$ and by the induction hypothesis for some $\alpha
< \omega_1$ there is $q \in \text{ Lim}_{\frak t}(\bar Q \restriction
L^2_a)$ which is (there) above $p^*_\alpha$ and above $p^*_\beta$.

Let $q^+ = q \cup (p_\alpha \restriction (w_\alpha \backslash w)) \cup
(p_\beta \restriction (w_\beta \backslash w))$ and let $p^+_\alpha = q
\cup (p_\alpha \restriction (w_\alpha \backslash w))$ and $p^+_\beta
= q \cup (p_\beta \restriction (w_\beta \backslash w))$.  Clearly $L_a
\le_{\frak t} L_2$ hence by clause $(i)(\beta) + (\gamma)$ for $\bar Q
\restriction (L^2_{a_\alpha} \cup L^2_{a_\beta})$ 
we have $p^+_\alpha \in \text{ Lim}_{\frak t}(\bar Q),
q \le p^+_\alpha,p_\alpha \le p^+_\alpha$ and
similarly $p^+_\beta \in \text{ Lim}_{\frak t}(\bar Q),q \le
p^+_\beta,p_\beta \le p^+_\beta$ clause $(B)(i)(\alpha)$ our $q =
p^+_\alpha \cup p^+_\beta$ is as required.
\mn
\ub{Clause (k)}:

Easy.     \hfill$\square_{\scite{ad.5}}$\margincite{ad.5}
\bigskip

\proclaim{\stag{ad.6} Claim}  1) Assume
\mr
\item "{$(a)$}"  ${\frak t}$ is an FSI-template,
Dp$_{\frak t}(L,\bar K) < \infty$ i.e. $\bar K$ is a smooth
${\frak t}$-memory choice
\sn
\item "{$(b)$}"  $\bar Q= \langle {\underset\tilde {}\to Q_{t,
{\underset\tilde {}\to \eta_t}}}:t \in L \rangle$ is a $({\frak t},
\bar K$)-iteration of def-c.c.c. forcing notions
\sn
\item "{$(c)$}"  $L_1,L_2 \subseteq L$ and $L_1 < L_2$ (that is $(\forall
t_1 \in L_1)(\forall t_2 \in L_2)(L^{\frak t} \models t_1 < t_2))$ and
$t \in L_2 \Rightarrow L_1 \in I^{\frak t}_t$ or at least $t \in L_2 \and
L' \subseteq L_1 \and |L'| \le \aleph_0 \Rightarrow L' \in I^{\frak t}_t$
and $L = L_1 \cup L_2$.
\ermn
\ub{Then}
\mr
\item "{$(\alpha)$}"  Lim$_{\frak t}(\bar Q)$ is actually a definition of
a forcing (in fact c.c.c. one) so meaningful in bigger universes, moreover
for extensions $\bold V_1 \subseteq \bold V_2$ of
$\bold V = \bold V_0$ (with the same ordinals of course), 
we \footnote{of course possibly $L_1 = \emptyset$} 
get $[\text{Lim}_{\frak t}(\bar Q)]^{{\bold
V}_1} \lessdot [\text{Lim}_{\frak t}(\bar Q)]^{{\bold V}_2}$
\sn
\item "{$(\beta)$}"  Lim$_{\frak t}(\bar Q)$ is in fact 
$Q_1 * {\underset\tilde {}\to Q_2}$ where $Q_1 = \text{ Lim}_{\frak t}
(\bar Q \restriction L_1)$ and $Q_2 = [\text{Lim}_{\frak t}(\bar Q 
\restriction L_2)]^{{\bold V}[{\underset\tilde {}\to G_{Q_1}}]}$
(composition).
\ermn
2) Assume clauses (a), (b) of part (1) and
\mr
\item "{$(c)_2$}"  $L$ has a last element $t^*$ and let $L^- = L
\backslash \{t^*\}$.
\ermn
\ub{Then} for any $G^- \subseteq \text{ Lim}_{\frak t}(\bar Q \restriction 
L^-)$ generic over $\bold V$, letting $\eta_{t^*} = 
{\underset\tilde {}\to \eta_{t^*}}[G^-]$ in $\bold V[G^-]$ we have: the forcing notion
Lim$_{\frak t}(\bar Q)/G^-$ is equivalent to 
$\cup \{Q^{\bold V[G^-_A]}_{t^*,\eta_{t^*}}:A \in I^{\frak t}_{t^*}$ 
is $\bar K$-closed$\}$ where $G^-_A =: G^- \cap
\text{ Lim}_{\frak t}(\bar Q \restriction A)$ and $\eta_{t^*_1} =
{\underset\tilde {}\to \eta_{t^*}}[G^-]$. \nl
3) Assume clauses (a), (b) of part (1) and
\mr
\item "{$(c)_3$}"  $\langle L_i:i < \delta \rangle$ is an increasing
continuous sequence of initial segments of $L$ with union $L$ and $\delta$
is a limit ordinal.
\ermn
\ub{Then} Lim$_{\frak t}(\bar Q)$ is
$\dbcu_{i < \delta} \text{ Lim}_{\frak t}(\bar Q \restriction L_i)$,
moreover $\langle \text{ Lim}_{\frak t}(\bar Q \restriction L_i):i < \delta
\rangle$ is $\lessdot$-increasing continuous.  \nl
4) Assume ${\frak t}^1,{\frak t}^2$ are FSI-templates, $L^{{\frak
t}^1} = L^{{\frak t}^2}$ call it $L$ and for every $t \in L,I^{{\frak
t}^1}_t \cap [L]^{\le \aleph_0} = I^{{\frak t}^2}_t \cap [L]^{\le
\aleph_0}$ and $\bar K$ is smooth ${\frak t}^\ell$-memory choice and
$\bar Q = \langle Q_{t,\bar \varphi_t,{\underset\tilde {}\to
\eta_t}}:t \in L \rangle$ is a $({\frak t}^\ell,\bar K)$-iteration of
def-c.c.c. forcing notions for $\ell = 1,2$.  Then Lim$_{{\frak
t}^1}(\bar Q) = \text{Lim}_{{\frak t}^2}(\bar Q)$.
\endproclaim
\bigskip

\demo{Proof}  Straight (or read \cite{Sh:630}).  \hfill$\square_{\scite{ad.6}}$\margincite{ad.6}
\enddemo
\bn
We now give sufficient conditions for: ``if we force by Lim$_{\frak t}
(\bar Q)$ from \scite{ad.5}, then some cardinal invariants are small or
equal/bigger than some $\mu$.  The necessity of such a claim in our framework
is obvious; we deal with two-place relations only as this is the case in the
popular cardinal invariants, in particular those we deal with.
\proclaim{\stag{ad.7} Claim}  Assume ${\frak t},\bar K$ and $\bar Q = 
\langle {\underset\tilde {}\to Q_{t,{\underset\tilde {}\to \eta_t}}}:t \in 
L^{\frak t} \rangle$ are as in \scite{ad.5} and $P = \text{Lim}_{\frak t}
(\bar Q)$. \nl
1) Assume
\mr
\item "{$(a)$}"  $R$ is a Borel \footnote{here and below just enough
absoluteness is enough, of course}  two-place relation on ${}^\omega \omega$
(we shall use $<^*$)
\sn
\item "{$(b)$}"  $L^* \subseteq L^{\frak t}$
\sn
\item "{$(c)$}"  for every countable $\bar K$-closed $A \subseteq L^{\frak t}$ for some
$t \in L^*$ we have $A \in I^{\frak t}_t$
\sn
\item "{$(d)$}" for $t \in L^*$ and $\bar K$-closed $A \in K^{\frak
t}_t$ which include $K_t$, in $\bold V^{\text{Lim}_{\frak t}(\bar Q
\restriction A)}$  we have $\Vdash_{{\underset\tilde {}\to Q_{t,
{\underset\tilde {}\to \eta_t}}}} ``{\underset\tilde {}\to \nu_t} \in
{}^\omega \omega$ is an $R$-cover of the old reals, that is
$\rho \in ({}^\omega \omega)^{\bold V} \Rightarrow \rho R  
{\underset\tilde {}\to \nu_t}$" where
${\underset\tilde {}\to \nu_t}$ is a name in the forcing
${\underset\tilde {}\to Q_{t,{\underset\tilde {}\to \eta_t}}}$ i.e. in
$(Q_{t,{\underset\tilde {}\to \eta_t}[\underset\tilde {}\to
G}])^{\bold V[\underset\tilde {}\to G]},\underset\tilde {}\to G$ the
generic subset of Lim$_{\frak t}(\bar Q \restriction A)$; not
depending on $A$.  (Usually
${\underset\tilde {}\to \nu_t}$ is the generic
real of 
${\underset\tilde {}\to Q_{t,{\underset\tilde {}\to \eta_t}}}$, and
hence
${\underset\tilde {}\to Q_{t,{\underset\tilde {}\to \eta_t}}}$ is
interpreted in the universe $\bold V^{\text{Lim}_{\frak t}(\bar Q
\restriction A)}$, so ${\underset\tilde {}\to \eta_t}$ is determined
by the generic; normally this we assume absolutely).
\ermn
\ub{Then} $\Vdash_P ``(\forall \rho \in {}^\omega \omega)(\exists t \in L^*)
(\rho R {\underset\tilde {}\to \nu_t})$, i.e. 
$\{{\underset\tilde {}\to \nu_t}:t \in L^*\}$ is an $R$-cover, which, if
$R = <^*$ means ${\frak d} \le |L^*|$". \nl
1A) If we weaken assumption (d) to ``for some ${\underset\tilde {}\to
\nu_t}$ a Lim$_{\frak t}(\bar Q \restriction K_t)$-name" we get
$\Vdash_P ``(\forall \rho \in {}^\omega \omega)(\exists t)(\exists \nu
\in \bold V^{\text{Lim}_{\frak t}(\bar Q \restriction K_t)})[\rho R
\nu]$. \nl
2) Assume
\mr
\item "{$(a)$}"  $R$ is a Borel two-place relation on ${}^\omega \omega$
(we shall use $<^*$)
\sn
\item "{$(b)$}"  $\mu$ is a cardinality
\sn
\item "{$(c)$}"  if $L^* \subseteq L^{\frak t},|L^*| < \mu$ \ub{then} for some
$t \in L^{\frak t}$ and $\bar K$-closed $L^{**} \supseteq L^*$ 
we have $L^{**} \in I^{\frak t}_t$ and in $\bold V^{\text{Lim}_{\frak
t}(\bar Q \restriction L^{**})},
\Vdash_{Q_{t,{\underset\tilde {}\to \eta_\delta}}} ``
{\underset\tilde {}\to \nu_t}$ is a $R$-cover of the old reals" 
with ${\underset\tilde {}\to \nu_t}$ some $Q_{t,{\underset\tilde {}\to
\nu_t}}$-name as in (1); (usually ${\underset\tilde {}\to \nu_t}$ is
the generic real of $Q_{t,{\underset\tilde {}\to \nu_1}}$ (this we
assume absolutely).
\ermn
\ub{Then} $\Vdash_P ``(\forall X \in [{}^\omega \omega]^{< \mu})(\exists \nu
\in {}^\omega \omega)(\dsize \bigwedge_{\rho \in X} \rho R \nu)"$ \nl
(so for $R = <^*$ this means ${\frak b} \ge \mu$). \nl
3) Assume
\mr
\item "{$(a)$}"  $R$ is a Borel 
two-place relation \footnote{so $R$ is defined in $\bold V$; if $R$ is
from $\bold V^{\text{Lim}_{\frak t}(\bar Q \restriction K)}$, we need
partial isomorphism (see below) of $({\frak t},\bar Q)$ extending id$_K$}
on ${}^\omega \omega$ (we use $R = \{(\rho,\nu):\rho,\nu
\in {}^\omega 2$ and $\rho^{-1}\{1\},\nu^{-1}\{1\}$ are infinite with finite
intersection)
\sn
\item "{$(b)$}"  $\kappa$ a cardinality, cf$(\kappa) > 2^{\aleph_0}$ and $\kappa <
\lambda$
\sn
\item "{$(c)$}"  if $t_{i,n} \in L^{\frak t}$ for $i < i(*),n < \omega$ and
$\kappa \le i(*) < \lambda$ and each: $\{t_{i,n}:n < \omega\}$ is
$\bar K$-closed, \ub{then} we can find $t_n \in L^{\frak t}$ for
$n < \omega$ such that $\{t_n:n < \omega\} \subseteq L^{\frak t}$ is
$\bar K$-closed and:
{\roster
\itemitem{ $(*)$ }  for every $i < i(*)$ for some $j < \kappa,j \ne i$ and the mapping
$t_{i,n} \mapsto t_{i,n},t_{j,n} \mapsto t_n$ is a partial 
isomorphism of $({\frak t},\bar Q)$ (see \scite{ad.8} below).
\endroster}
\ermn
\ub{Then} in $\bold V^P$ we have
\mr
\item "{$\boxtimes^R_\mu$}"  if $\rho_i,\nu_i \in {}^\omega \omega$ 
for $i < i(*)$ and $\mu \le i(*) < \lambda$ and $i \ne j \Rightarrow \nu_i
R \rho_j$, \ub{then} we can find $\rho \in {}^\omega \omega$ such that
$i < i(*) \Rightarrow \nu_i R \rho$.
\endroster
\endproclaim
\bigskip

\demo{Proof}  Straight, but being requested: \nl
1) Let $\underset\tilde {}\to \rho$ be a $P$-name of a member of
$({}^\omega \omega)^{\bold V^P}$, so as $P$ satisfies (see
\scite{ad.4}(B)(j)), for each $n$ there is a maximal antichain
$\{p_{n,i}:i < i_n\}$ such that $p_{n,i}$ forces a value of
$\underset\tilde {}\to \rho(n)$ and, of course, $i_n$ is countable.
Let $M = \{a:a$ is a countable $\bar K$-closed subset of $L^{\frak
t}\}$, so obviously $M$ is closed under countable unions and
$\cup\{a:a \in M\} = L^{\frak t}$; and let $L_a=a$ for $a \in M$ so by
\scite{ad.4}(B)(g) we have Lim$_{\frak t}(\bar Q) =
\cup\{\text{Lim}_{\frak t}(\bar Q \restriction L_a):a \in M\}$ but $P
= \text{ Lim}_{\frak t}(\bar Q)$, hence for $n < \omega,i < i_n$ for some
$a_{n,i} \in M$ we have $p_{n,i} \in \text{ Lim}_{\frak t}(\bar Q
\restriction L_a)$.  But $M$ is $\aleph_1$-directed so for some $a \in
M$ we have $\{a_{n,i}:n < \omega,i < i_n\} \subseteq \text{ Lim}_{\frak
t}(\bar Q \restriction L_a)$.  Also by \scite{ad.4}(B)(e) we know
Lim$_{\frak t}(\bar Q \restriction L_a) \lessdot \text{ Lim}_{\frak
t}(\bar Q) = P$, so $\underset\tilde {}\to \rho$ is a Lim$_{\frak
t}(\bar Q \restriction L_a)$-name.  Now by assumption (c) of what we
are proving, as $L_a \subseteq L$ is countable, we can find $t \in L^*
\subseteq L^{\frak t}$ such that $L_a \in I^{\frak t}_t$.  Also we
know that $K_t \in I^{\frak t}_t$ (see Definition \scite{ad.1}(2)(c)
hence $A =: K_t \cup L_a$ belongs to $I^{\frak t}_t$ and is $\bar
K$-closed; and easily also $B = A \cup \{t\}$ is $\bar K$-closed.

So $A \subseteq B \subseteq L^{\frak t}$ are $\bar K$-closed so as
above Lim$_{\frak t}(\bar Q \restriction A) \lessdot \text{
Lim}_{\frak t}(\bar Q \restriction B) \lessdot \text{ Lim}_{\frak
t}(\bar Q) = P$ and $\underset\tilde {}\to \rho$ is a Lim$_{\frak
t}(\bar Q \restriction A)$-name (hence also a Lim$_{\frak t}(\bar Q
\restriction B))$ of a member of ${}^\omega \omega$.

Now by assumption (d) in $\bold V^{\text{Lim}_{\frak t}(\bar Q
\restriction A)}$ we have $\Vdash_{Q_{t,n_i}} ``\underset\tilde {}\to
\rho R {\underset\tilde {}\to \nu_t}"$, hence by \scite{ad.4}(B)(h) we
know that Lim$_{\frak t}(\bar Q \restriction B) = \text{ Lim}_{\frak
t}(\bar Q \restriction A) * {\underset\tilde {}\to
Q_{t,{\underset\tilde {}\to n_t}}}$, so  together
$\Vdash_{\text{Lim}_{\frak t}(B)} ``\underset\tilde {}\to \rho R
{\underset\tilde {}\to \nu_t}"$ hence the previous sentence and
obvious absoluteness we have $\Vdash_P 
``\underset\tilde {}\to \rho R {\underset\tilde {}\to \nu_t}"$.  So as
$\underset\tilde {}\to \rho$ was any $P$-name of a member of
$({}^\omega \omega)^{\bold V^P}$ we are done. \nl
1A)  Same proof. \nl
2) So assume $p \Vdash_P ``\underset\tilde {}\to X \subseteq {}^\omega
\omega$ has cardinality $< \mu$".  As we can increase $p$ \wilog \,
for some $\theta < \mu$ we have $p \Vdash_P ``|\underset\tilde {}\to
X| = \theta"$ so we can find a sequence $\langle {\underset\tilde
{}\to \rho_\alpha}:\alpha < \theta \rangle$ of $P$-names of members of
$({}^\omega \omega)^{\bold V^P}$ such that $p \Vdash_P
``\underset\tilde {}\to X = \{{\underset\tilde {}\to
\rho_\alpha}:\alpha < \theta\}"$.  Let $\{p_{\alpha,n,i}:i <
i_{\alpha,n}\}$ be a maximal antichain of $P$, with $p_{\alpha,n,i}$
forcing a value to ${\underset\tilde {}\to \rho_\alpha(n)}$ and
$i_{\alpha,n}$ countable.

Define $M = \{a \subseteq L^{\frak t}$: a countable $\bar
K$-closed$\}$, so for each $\alpha < \theta,n < \omega,i <
i_{\alpha,n}$ for some $a_{\alpha,n,i} \in M$ we have $p_{\alpha,n,i}
\in \text{ Lim}_t(\bar Q \restriction L_a)$.  So for some $\bar
K$-closed $L^{**} \subseteq L^{\frak t}$ and $t \in L^{\frak t}$ we
have $L^{**} \in I^{\frak t}_t$ and $a_{\alpha,n,i} \subseteq L^{**}$
for $\alpha < \theta,n < \omega,i < i_{\alpha,n}$.  We now continue
as in part (1). \nl
3) So assume $i(*) \in [\kappa,\lambda)$ and $\Vdash_P
``{\underset\tilde {}\to \nu_i},{\underset\tilde {}\to \rho_i} \in
{}^\omega \omega$ and $i \ne j \Rightarrow {\underset\tilde {}\to
\nu_i} R {\underset\tilde {}\to \rho_j}"$.  So as above we can find
countable $\bar K$-closed $K^*_i \subseteq L^{\frak t}$ such that
${\underset\tilde {}\to \nu_i},{\underset\tilde {}\to \rho_i}$ are
Lim$_{\frak t}(\bar Q \restriction K^*_i)$-names; \wilog \, $K^*_i \ne
\emptyset$ and even $|K^*_i| = \aleph_0$; this is impossible only if $L^{\frak t}$
is finite and then all is trivial.  Let $\langle t_{i,n}:n < \omega
\rangle$ be a list of the members of $K^*_i$ with no repetitions.  Let $f_{i,j}$ be
the mapping from $K^*_j$ to $K^*_i$ defined by $f_{i,j}(t_{j,n}) =
t_{i,n}$.

We define a two place relation $E_1,E_2$ on $i(*)$ and on $i(*) \times
i(*)$ respectively

$$
\align  
i E_1 j \text{ \ub{iff} } &f_{i,j} \text{ is a partial isomorphism of }
({\frak t},\bar Q) \\
  &\text{such that } \hat f_{i,j} \text{ maps } ({\underset\tilde
{}\to \rho_j},{\underset\tilde {}\to \nu_j}) \text{ to }
({\underset\tilde {}\to \rho_i},{\underset\tilde {}\to \nu_i})
\endalign
$$

$$
\align
(i_1,i_2) E_2 (j_1,j_2) \text{ \ub{iff} } &i_1 E_1 j_1,i_2 E_2 j_2 \text{
and} \\
  &f_{i_1,j_1} \cup f_{i_2,j_2} \text{ is a partial isomorphism of }
({\frak t},\bar Q).
\endalign
$$
\mn
Easily
\mr
\item "{$\boxtimes(i)$}"  $E_1,E_2$ are equivalence classes over their
domain
\sn
\item "{$(ii)$}"  $E_1,E_2$ has $\le 2^{\aleph_0}$ equivalence classes
\sn
\item "{$(iii)$}"  $f_{j,i} = f^{-1}_{i,j}$.
\ermn
As $|i(*)/E_1| \le 2^{\aleph_0} < \text{ cf}(\kappa)$ (by $(*)(ii)$
and assumption (b) respectively) and we can replace $i(*)$ by $i(*) +
\kappa$, \wilog \, $i < \kappa \Rightarrow 0E_1i$.  Now we apply
assumption (c), and get $\langle t_n:n < \omega \rangle$.  By $(*)$ of
clause (c) for any $i,j$ clearly $K^*_i \cup K^*_j$ and
$K^*_i \cup \{t_n:n < \omega\}$ are $\bar K$-closed (see the definition below).  For any $i < i(*)$ let $j_i <
\kappa$ be as in $(*)$ of clause (c) which means: $j_i \ne i$ and the
following mapping $g_i$ is a partial isomorphism of $({\frak t},\bar
Q) = \text{ Dom}(g_i) = \{t_{i,n},t_{j_i,n}:n < \omega\},g_i(t_{i,n})
= t_{i,n},g_i(t_{j,n}) = t_n$.

Let $\underset\tilde {}\to \nu,\underset\tilde {}\to \rho$ be
Lim$_{\frak t}(\bar Q \restriction K^*)$-names such that for some
equivalently any $i,\hat g_i$ maps ${\underset\tilde {}\to \nu_{j_i}},
{\underset\tilde {}\to \rho_{j_i}}$ to $\underset\tilde {}\to
\nu,\underset\tilde {}\to \rho$ respectively (this is O.K. as for any
$i_1,i_2$ we have $j_{i_1} E_1 j_{i_2}$ because $j_{i_1},j_{i_2}$
hence $g_{i_2} \circ f_{j_2,j_{i_1}} = g_{i_1} \restriction
K^*_{j_{i_1}}$).  Now for any $i < \mu$, as $j_i \ne i$, we know
$\Vdash_{\text{Lim}_{\frak t}(\bar Q \restriction (K^*_i \cup
K^*_{j_i}))} ``{\underset\tilde {}\to \nu_i} R {\underset\tilde {}\to
\rho_{j_i}}"$, so applying $g_i$ we have $\Vdash_{\text{Lim}_{\frak
t}(K^*_i \cup K^*)} ``{\underset\tilde {}\to \nu_i} R \underset\tilde
{}\to \rho"$.  So we have proved $\boxtimes^R_\mu$.
\hfill$\square_{\scite{ad.7}}$\margincite{ad.7}
\enddemo
\bn
In \scite{ad.8} we note that isomorphisms (or embeddings) of ${\frak t}$'s
tend to induce isomorphisms (or embeddings) of Lim$_{\frak t}(\bar Q)$, and
deal (in \scite{ad.9},\scite{ad.10}) with some natural operation.  In
\scite{ad.8} we could use two ${\frak t}$'s, but this can trivially be
reduced to one.
\definition{\stag{ad.8} Definition/Claim}  Assume that 
${\frak t},\bar K$ and $\bar Q = \langle
{\underset\tilde {}\to Q_{t,{\underset\tilde {}\to \eta_t}}}:t \in
L^{\frak t} \rangle$ are as in \scite{ad.5}.  
By induction on $\zeta$ we define
and prove \footnote{if $K_t = \emptyset$ and all
${\underset\tilde {}\to Q_{t,\eta}}$ have the same
definition of forcing notion, as in our main case, we can separate the
definition and claim}
\mr
\item "{$(A)$}"  [Def] $\quad$ we say $f$ is a partial isomorphism of
$({\frak t},\bar Q)$ if: (writing ${\frak t}$ instead of \nl

$\qquad ({\frak t},\bar Q)$ means we assume
$Q_{t,{\underset\tilde {}\to \eta_t}} = Q$, i.e. constant)
{\roster
\itemitem{ $(a)$ }  $f$ is a partial one-to-one function from
$L^{\frak t}$ to $L^{\frak t}$
\sn
\itemitem{ $(b)$ }  Dom$(f)$, Rang$(f)$ are $({\frak t},\bar K)$-closed sets
of depth $\le \zeta$
\sn
\itemitem{ $(c)$ }  for $t \in \text{ Dom}(f)$ and $A \subseteq
\text{ Dom}(f)$ we have $A \in I^{\frak t}_t \Leftrightarrow f''(A) \in
I^{\frak t}_{f(t)}$
\sn
\itemitem{ $(d)$ }  for $t \in \text{ Dom}(f)$, we have: $f$ maps $K_t$
onto $K_{f(t)}$ and $f \restriction K_t$ maps ${\underset\tilde {}\to \eta_t}$
to ${\underset\tilde {}\to \eta_{f(t)}}$, more exactly the isomorphism
which $f$ induces  from 
Lim$_{\frak t}(\bar Q \restriction K_t)$ onto Lim$_{\frak t}(\bar Q
\restriction K_{f(t)})$ does this.
\endroster}
\item "{$(B)$}"  [Claim] $\quad f$ induces naturally an isomorphism
which we call $\hat f$ from Lim$(\bar Q \restriction \text{ Dom}(f))$ \nl

$\qquad \quad$ onto Lim$_{\frak t}(\bar Q
\restriction \text{ Rang}(f))$.
\endroster
\enddefinition
\bigskip

\demo{Proof}  Straightforward.
\enddemo
\bigskip

\definition{\stag{ad.9} Definition}  1) We say ${\frak t} = {\frak t}^1 +
{\frak t}^2$ if
\mr
\item "{$(a)$}"  $L^{\frak t} = L^{{\frak t}^1} + L^{{\frak t}^2}$ (as
linear orders)
\sn
\item "{$(b)$}"  for $t \in L^{{\frak t}^1},I^{{\frak t}^1}_t =
I^{\frak t}_t$
\sn
\item "{$(c)$}"  for $t \in L^{{\frak t}^2},I^{{\frak t}^2}_t = \{A
\subseteq L^{\frak t}:A \cap L^{{\frak t}^2} \in I^{{\frak t}^2}_t\}$.
\ermn
So ${\frak t}^1 + {\frak t}^2$ is well defined if ${\frak t}^1,{\frak t}^2$
are disjoint, i.e. $L^{{\frak t}^1} \cap L^{{\frak t}^2} = \emptyset$. \nl
2) We say ${\frak t}^1 \le_{wk} {\frak t}^2$ iff
\mr
\item "{$(a)$}"  $L^{{\frak t}^1} \subseteq L^{{\frak t}^2}$ (as linear 
orders)
\sn
\item "{$(b)$}"  for $A \subseteq L^{{\frak t}^1}$ and $t \in L^{{\frak t}^1}$
we have $A \in I^{{\frak t}^1}_t \Leftrightarrow A \in I^{{\frak t}^2}_t$.
\ermn
3) If $\langle {\frak t}^\zeta:\zeta < \xi \rangle$ is 
$\le_{wk}$-increasing, $\xi$ a
limit ordinal, we define ${\frak t}^\xi =: \dbcu_{\zeta < \xi} {\frak t}
^\zeta$ by

$$
L^{{\frak t}^\xi} = \dbcu_{\zeta < \xi} L^{{\frak t}^\zeta} \qquad
\text{ (as linear orders)}
$$

$$
I^{{\frak t}^\xi}_t = \cup \{I^{{\frak t}^\zeta}_t:\zeta < \xi \text{
and } t \in L^{{\frak t}_\zeta}\}
$$
\mn
Clearly $\zeta < \xi \Rightarrow {\frak t}^\zeta \le_{wk} {\frak
t}^\xi$.  Such ${\frak t}^\xi$ is called the limit of $\langle {\frak
t}^\zeta:\zeta < \xi \rangle$; now a $\le_{wk}$-increasing
sequence $\langle {\frak t}^\zeta:\zeta < \xi \rangle$ is called
continuous if for every limit ordinal $\delta < \xi$ we have ${\frak
t}^\delta = \dbcu_{\zeta < \delta} {\frak t}^\zeta$.
\nl
4) If $\langle {\frak t}^\zeta:\zeta < \xi \rangle$ are pairwise disjoint
(that is $\zeta \ne \varepsilon \Rightarrow L^{{\frak t}^\zeta} \cap
L^{{\frak t}^\varepsilon} = \emptyset$) 
we define $\dsize \sum_{\zeta < \xi} {\frak t}^\zeta$ by induction on $\xi$
naturally: for $\xi = 1$ it is ${\frak t}^0$, for $\xi$ limit it is
$\dbcu_{\zeta < \xi} (\dsize \sum_{\zeta < \varepsilon} {\frak t}^\zeta)$ and
for $\xi = \varepsilon + 1$ it is $(\dsize \sum_{\zeta < \varepsilon}
{\frak t}^\zeta) + {\frak t}^\varepsilon$, so $\xi_1 \le \xi_2 \Rightarrow
\dsize \sum_{\zeta < \xi_1} {\frak t}^\zeta \le_{wk} \dsize \sum_{\zeta < \xi_2}
{\frak t}^\zeta$ (even an initial segment). \nl
5) We can replace in 0) - 4) above ${\frak t}^\zeta$ by $({\frak t}^\zeta,\bar K^\zeta)$. 
\enddefinition
\bigskip

\proclaim{\stag{ad.10} Claim}  Let ${\frak t}$ be an FSI-template. \nl
1) If $L^{\frak t} = \emptyset$ or just $L^{\frak t}$ is finite \ub{then}
${\frak t}$ is smooth. \nl
2) If ${\frak t}^1,{\frak t}^2$ are disjoint FSI-templates, \ub{then}
${\frak t}^1 + {\frak t}^2$ is a FSI-template and $\rho \in \{1,2\}
\Rightarrow {\frak t}^\ell \le_{wk} {\frak t}^1 + {\frak t}^2$. \nl
3) If ${\frak t}^1,{\frak t}^2$ are disjoint smooth FSI-templates \ub{then}
${\frak t}^1 + {\frak t}^2$ is a smooth FSI-template. \nl
4) If $\langle {\frak t}^\zeta:\zeta < \xi \rangle$ is an 
$\le_{wk}$-increasing (\scite{ad.9}(2)) sequence of FSI-templates and 
$\xi$ is a limit ordinal, \ub{then} ${\frak t}^\xi =:
\dbcu_{\zeta < \xi} {\frak t}^\zeta$ is an FSI-template and $\zeta < \xi
\Rightarrow {\frak t}^\zeta \le_{wk} {\frak t}^\xi$. \nl
5) If $\langle {\frak t}^\zeta:\zeta < \xi \rangle$ is an increasing
continuous (see Definition \scite{ad.9}(3)) sequence of smooth
FSI-templates and $\xi$ is a limit ordinal, \ub{then} ${\frak t}^\xi =:
\dbcu_{\zeta < \xi} {\frak t}^\zeta$ is a smooth FSI-template and $\zeta < \xi
\Rightarrow {\frak t}^\zeta \le_{wk} {\frak t}^\xi$. \nl
6) If $\langle {\frak t}^\zeta:\zeta < \xi \rangle$ is a sequence of pairwise
disjoint [smooth] FSI-templates, \ub{then} $\dsize \sum_{\zeta < \xi}
t^\zeta$ is a [smooth] FSI-template and $\langle \dsize \sum_{\zeta <
\varepsilon} t^\zeta:\varepsilon \le \zeta \rangle$ is increasing continuous.
\nl
7) We can add $\bar K^\zeta$ to ${\frak t}^\zeta$. \nl
8) We can restrict ourselves to locally countable ${\frak t}$'s (so the sums
are locally countable if the summands are).
\endproclaim
\bigskip

\demo{Proof}  Easy.
\enddemo
\bn
\centerline {$* \qquad * \qquad *$}
\bn
\ub{\stag{ad.10a} Discussion}:
To prove our desired result CON$({\frak a} > {\frak d})$ we need to construct
an FSI-template ${\frak t}$ of the right form.  Now we do it using a 
measurable
cardinal.  The point is that if we are given $\left< \langle t_{i,n}:n <
\omega \rangle:i < i(*) \right>,L^{\frak t},i(*) \ge \kappa$ and $D$ is a
normal ultrafilter on $\kappa$, then in ${\frak t}^\kappa/D,\left< \langle
t_{i,n}:i < \kappa \rangle/D:n < \omega \right>$ is as required in
\scite{ad.7}(3)(c), considering ${\frak t}^\kappa/D$ an extension of
${\frak t}$.  
\bigskip

\definition{\stag{ad.11} Definition}  For a template ${\frak t}$ and an
$(2^{\aleph_0})^+$-complete ultrafilter $D$ on $\kappa$ we define 
${\frak t}^* =: {\frak t}^\kappa/D$ as follows:

$$
L^{{\frak t}^*} = (L^{\frak t})^\kappa/D \text{ as a linear order}
$$
\mn
and if $t^* = \langle t_i:i < \kappa \rangle/D$ where 
$t_i \in L^{\frak t}$ then we let
$I^{{\frak t}^*}_{t^*} = \{A:\text{we can find } A_i \in I^{\frak t}_{t_i}$
for $i < \kappa$ such that $A \subseteq \dsize \prod_{i < \kappa} A_i/D\}$.
We then let $\bold j_{D,{\frak t}}$ be the canonical embedding of ${\frak t}$
into ${\frak t}^\kappa/D$
and ${\frak t}' = \bold j_{D,{\frak t}}({\frak t})$ is defined by
$L^{{\frak t}'} = L^{t^*} \restriction \{\bold j_{d,{\frak t}}(s):s
\in L^{\frak t}\},I^{{\frak t}'}_s = I^{{\frak t}^*}_s \restriction L^{{\frak t}'}$. \nl
[We can deal with $\bar K$, if $D$ if $(\dbcu_{t \in L} |K_t|^+)$-complete
and can deal also with 
$\bar Q$ if we have $<$ com$(D)$ kinds of $\bar \varphi_t$.]
\enddefinition
\bigskip

\proclaim{\stag{ad.12} Claim}  1) In \scite{ad.11}, 
${\frak t}^\kappa/D$ is also a template and 
$\bold j_{D,{\frak t}}({\frak t}) \le_{wk} {\frak t}^\kappa/D$. \nl
2) If ${\frak t}$ is a smooth FSI-template \ub{then} ${\frak t}^\kappa/D$ is a
smooth FSI-template.
\endproclaim
\bigskip

\demo{Proof}  Straight.

Now \scite{ad.13}, \scite{ad.14} below are used only in the short proof of
\scite{ad.15} depending on \S1, so you may ignore them.
\enddemo
\bigskip

\definition{\stag{ad.13} Definition}  Fix $\aleph_0 < \kappa < \mu =
\text{ cf}(\mu) < \lambda = \text{ cf}(\lambda) = \lambda^\kappa$
and $D$ a $\kappa$-complete (or just $(2^{\aleph_0})^+$-complete) uniform
ultrafilter on $\kappa$.  We define by induction on $\zeta \le \lambda$,
smooth FSI-template ${\frak t}_{\gamma,\zeta}$ for $\gamma < \mu$ such that:
\mr
\item "{$(a)$}"  ${\frak t}_{\gamma,\zeta}$ is a locally countable
FSI-template
\sn
\item "{$(b)$}"  if 
$\gamma_1 \ne \gamma_2$ then ${\frak t}_{\gamma_1,\zeta},
{\frak t}_{\gamma_2,\zeta}$ are disjoint, i.e. 
$L^{{\frak t}_{\gamma_1,\zeta}} \cap L^{{\frak t}_{\gamma_2,\zeta}}
= \emptyset$
\sn
\item "{$(c)$}"  for $\xi < \zeta$ we have ${\frak t}_{\gamma,\xi} 
\le_{wk} {\frak t}_{\gamma,\zeta}$
\sn
\item "{$(d)$}"  if $\zeta$ is limit then ${\frak t}_{\gamma,\zeta} =
\dbcu_{\xi < \zeta} {\frak t}_{\gamma,\xi}$, see \scite{ad.9}(3), \scite{ad.10}(6).
\sn
\item "{$(e)$}"  if $\zeta = \xi +1$ and $\xi$ is even, \ub{then} there is an
isomorphism $\bold j_{\gamma,\zeta}$ from $\dsize \sum_{\beta \le \gamma} \,
{\frak t}_{\beta,\xi}$ onto ${\frak t}_{\gamma,\zeta}$ which is the identity
over ${\frak t}_{\gamma,\xi}$
\sn
\item "{$(f)$}"  if $\zeta = \xi +1$ and $\xi$ is odd, \ub{then} there is an
isomorphism $\bold j_{\gamma,\zeta}$ from $({\frak t}_{\gamma,\xi})^\kappa/D$
onto ${\frak t}_{\gamma,\zeta}$ which extends the inverse of
$\bold j_{D,{\frak t}_{\gamma,\xi}}$.
\endroster 
\enddefinition
\bn
\ub{\stag{ad.14} Observation}:  The definition is \scite{ad.13} is legitimate.
\bigskip

\demo{Proof}  By the previous claims.
\enddemo
\bn
\ub{\stag{ad.15} Conclusion}:  Assume: $\kappa$ is measurable, $\kappa < \mu
= \text{ cf}(\mu) < \lambda = \text{ cf}(\lambda) = \lambda^\kappa$.
\ub{Then} for some c.c.c. forcing notion $P$ of cardinality $\lambda$, in
$\bold V^P$ we have ${\frak a} = \lambda,{\frak b} = {\frak d} = \mu$.
\bigskip

\demo{Proof}  \ub{Short Proof} (depending on \S1).  Let ${\frak t}_{\gamma,\zeta}$ 
(for $\gamma < \mu,\zeta \le \lambda)$ 
be as in \scite{ad.13}.  Let ${\frak t} = \dsize \sum_{\gamma < \mu} 
{\frak t}_{\gamma,\lambda}$ and let 
$\bar K = \langle K_t:t \in L^{\frak t} \rangle,
K_t = \emptyset$ and let $\bar Q = \langle
{\underset\tilde {}\to Q_t}:t \in L^{\frak t} \rangle$ with 
${\underset\tilde {}\to Q_t}$ being constantly the
dominating real forcing (= Hechler forcing).  Lastly let 
$P = \text{ Lim}_{\frak t}(\bar Q)$.
\nl
The rest is as in the end of \S1.
\enddemo
\bn
\ub{Alternative presentation, self contained not depending on
\scite{ad.13}, \scite{ad.14}}:  We define an 
FSI-template ${\frak t}^\zeta$
for $\zeta \le \lambda$ by induction on $\zeta$.
\mn
\ub{Case 1}:  For $\zeta = 0$. \nl

Let ${\frak t}^\zeta$ be defined as follows:

$$
L^{{\frak t}^\zeta} = \mu
$$

$$
I^{{\frak t}^\zeta}_\alpha = \{A:A \subseteq \alpha \text{ is countable}\}
$$
\mn
\ub{Case 2}:  For $\zeta = \xi +1$. \nl

We choose ${\frak t}^\zeta$ such that there is an isomorphism $\bold j_\zeta$
from $L^{{\frak t}^\zeta}$ onto $(L^{{\frak t}^\xi})^\kappa/D$, such that
$\bold j_\zeta \restriction L^{{\frak t}^\xi}$ is the canonical
embedding $\bold j_{D,{\frak t}^\xi}$,
and if $x \in L^{{\frak t}^\zeta},\bold j_\zeta(x) = \langle x_\varepsilon:
\varepsilon < \kappa \rangle/D \in (L^{{\frak t}^\xi})^\kappa/D$ then:

$$
\align
A \in I^{{\frak t}^\zeta}_x &\text{ \ub{iff} for some } \bar A = \langle
A_\varepsilon:\varepsilon < \kappa \rangle \text{ we have} \\
  &A_\varepsilon \in I^{{\frak t}^\xi}_{x_\varepsilon} \text{ and} \\
  &\{\bold j_\zeta(y):y \in A\} \subseteq \{\langle y_\varepsilon:
\varepsilon < \kappa \rangle/D:\{\varepsilon < \kappa:y_\varepsilon \in
A_\varepsilon\} \in D\}
\endalign
$$
\mn
\ub{Case 3}:  $\zeta$ limit.

We choose ${\frak t}^\zeta$ as follows:

$$
L^{{\frak t}^\zeta} = \dbcu_{\xi < \zeta} L^{{\frak t}^\xi}
\text{ as linear orders}
$$

$$
I^{{\frak t}^\zeta}_x \text{ is } \{A:A \subseteq \{s:L^{{\frak t}^\zeta}
\models ``s < x"\}\}
$$
\mn
\ub{if} $x \in L^{{\frak t}^0}$ and \ub{is} otherwise \footnote{this
is the ``veteranity privilege", i.e. ``founding father right";
members $t$ of $L^{{\frak t}^0}$ have the maximal $I^{{\frak t}^\zeta}_t$.}

$$
\align
\{A:&\text{for some } \xi < \zeta \text{ we have } x \in L^{{\frak t}^\xi}
\text{ and if } y = \text{ Min}\{y \in L^{{\frak t}^0}:
L^{{\frak t}^\zeta} \models x < y\} \\
  &\text{which is } \in L^{{\frak t}^0} \text{ then} \\
  &A \backslash \{x \in L^{{\frak t}^\zeta}:L^{{\frak t}^\zeta} 
\models x < z \text{ for some } z \text{ such that }
L^{{\frak t}^0} \models z < y\} 
\text{ belongs to } I^{{\frak t}^\xi}_x\}.
\endalign
$$
\mn
We now prove by induction on $\zeta \le \lambda$ that:
\mr
\item "{$(*)(a)$}"  ${\frak t}^\zeta$ is an FSI-template
\sn
\item "{$(b)$}"  $L^{{\frak t}^0}$ is an unbounded subset of
$L^{{\frak t}^\zeta}$
\sn
\item "{$(c)$}"  ${\frak t}^\zeta$ is smooth
\sn
\item "{$(d)$}"  ${\frak t}^\xi \le_{wk} {\frak t}^\zeta$ for $\xi < \zeta$
\sn
\item "{$(e)$}"  if $y \in L^{{\frak t}^\zeta}$
then $\{z:\text{for some } x \in L^{{\frak t}^0}$ we have 
$L^{{\frak t}^\zeta} \models z \le x$ and \nl
$L^{{\frak t}^\zeta} \models x < y\} \in I^{{\frak t}^\zeta}_x$
\sn
\item "{$(f)$}"  $L^{{\frak t}^\zeta}$ has cardinality $\le (\mu + |\zeta|)^\kappa$.
\ermn
Lastly let for $\zeta \le \lambda,P_\zeta = \text{ Lim}_{\frak t}
(\bar Q \restriction L^{{\frak t}^\zeta})$.
Now
\mr
\item "{$(\alpha)$}"  $P_\lambda$ is a c.c.c. forcing notion of
cardinality $\le \lambda^{\aleph_0}$ hence $\bold V^{P_\lambda}
\models 2^{\aleph_0} \le \lambda$ by \scite{ad.4}(B)(j)
\sn
\item "{$(\beta)$}"  in $\bold V^{P_\lambda}$ we have ${\frak d} \le
\mu$, by \scite{ad.7}(1) applied with $R = <^*$ and $L^* = L^{{\frak
t}^0}$ using $(*)$(b)+(e)
\sn
\item "{$(\gamma)$}"  in $\bold V^{P_\lambda}$ we have ${\frak b} \ge
\mu$ by \scite{ad.7}(2) applied with $R = <^*$
\sn
\item "{$(\delta)$}"  ${\frak b} = {\frak d} = \mu$ and ${\frak a} \ge
\mu$ by $(\beta) + (\gamma)$ as it is well known that ${\frak b} \le
{\frak d}$ and ${\frak b} \le {\frak a}$.
\ermn
But why the demand (c) from \scite{ad.7}(3) holds?  So assume $i(*)
\in [\kappa,\lambda)$ and $t_{i,n} \in L^{{\frak t}^\lambda}$ for $i <
i(*),n < \omega$ be given.  As $\lambda$ is regular $> i(*)$,
necessarily for some $\xi < \lambda$ we have $\{t_{i,n}:i < i(*),n <
\omega\} \subseteq L^{{\frak t}^\xi}$.  Now let $t_n \in L^{t^{\xi
+1}}$ be such that $\bold j_{\xi +1}(t_n) = \langle t_{i,n}:i < \kappa
\rangle/D$; so
\mr
\item "{$(\varepsilon)$}"  in $\bold V^{P_\lambda}$ we have ${\frak a}
\ge \kappa \Rightarrow {\frak a} \ge \lambda$ by \scite{ad.7}(3), see
there.
\ermn
Together we are done.  \hfill$\square_{\scite{ad.15}}$\margincite{ad.15}
\newpage

\head {\S3 eliminating the measurable} \endhead  \resetall \sectno=3
\bigskip

Without a measurable cardinal our problem is to verify condition (c) in
\scite{ad.7}(3).  Toward this it is helpful to show that for some $\aleph_1$-complete
filter $D$ on $\kappa$, for any $i(*) \in [\kappa,\lambda)$ and
$t_{i,n} \in L^{\frak t}$, for $i < i(*),n < \omega$, we have: for every $j <
i(*)$ for some $A \in D^+$ we have for any $i_0,i_1 \in A$, the mapping
$t_{j,n} \mapsto t_{j,n}$; $t_{i_0,n} \mapsto t_{i_1,n}$ is a partial
isomorphism of ${\frak t}$.  So $D$ behaves as an $\aleph_1$-complete
ultrafilter for our purpose. \nl
[If you know enough model theory, this is the problem of finding
convergent sequences, see \cite[Ch.II]{Sh:300}; 
for stable first order $T$ with $\kappa = \kappa_r (T)$ any 
indiscernible sequence (equivalently set) 
$\langle \bar a_\alpha:\alpha < \alpha^* \rangle$ of cardinality 
$\ge \kappa$, is convergent; why? as for any 
$\bar{\bold b} \in {}^{\kappa >} {\frak C}$, for all but $< \kappa$ ordinals 
$\alpha < \alpha^*,\bar{\bold b} \char 94 \bar a
_\alpha$ has a fixed type so average is definable.  In \cite[Ch.II]{Sh:300}, we
deal with it in general, (so harder to prove existence
which we do there under the relevant assumptions).]
\bigskip

\proclaim{\stag{abd.1} Lemma}  Assume $2^{\aleph_0} < \mu = \text{ cf}
(\mu) < \lambda = \text{ cf}(\lambda) = \lambda^{\aleph_0}$.  \ub{Then} for
some $P$ we have
\mr
\item "{$(a)$}"  $P$ is a c.c.c. forcing notion of cardinality $\lambda$
\sn
\item "{$(b)$}"  in $\bold V^P$ we have ${\frak b} = {\frak d} = \mu$ and
${\frak a} = 2^{\aleph_0} = \lambda$.
\endroster
\endproclaim
\bigskip

\demo{Proof}  We rely on \scite{ad.5} + \scite{ad.7}.  Let $L^+_0$ be a
linear order isomorphic to $\lambda$, let $L^-_0$ be a linear order
anti-isomorphic to $\lambda$ (and $L^-_0 \cap L^+_0 = \emptyset$) and let
$L_0 = L^-_0 + L^+_0$.

Let $\bold J$ be the following linear order:
\mr
\item "{$(a)$}"  its set of elements is ${}^{\omega >}(L_0)$
\sn
\item "{$(b)$}"  the order is: $\eta <_{\bold J} \nu$ \ub{iff} for some $n <
\omega$ we have $\eta \restriction n = \nu \restriction n$ and \nl
$\ell g(\eta) = n \and \nu(n) \in L^+_0$ or $\ell g(\nu) = n \and \eta(n) \in
L^-_0$ or we have
$\ell g(\eta) > n \and \ell g(\nu) > n \and L_0 \models \eta(n) <
\nu(n)$.
\ermn
[See more on such orders \cite{Lv} and \cite[AP]{Sh:220}, but we are
self contained.] \nl
Note that
\mr
\widestnumber\item{$\boxtimes^+$}
\item "{$\boxtimes$}"  every interval of $\bold J$ has cardinality
$\lambda$
\sn
\item "{$\boxtimes^+$}"  if $\aleph_0 < \theta = \text{ cf}(\theta) <
\lambda$ or $\theta = 1$ and $\langle t_i:i < \theta \rangle$ is a strictly decreasing
sequence in $\bold J$ then $\bold J \restriction \{y \in \bold J:(\forall i
< \theta)(y <_{\bold J} t_i)\}$ has cofinality $\lambda$
\sn
\item "{$\boxtimes^-$}"  the inverse of $\bold J$ satisfies
$\boxtimes^+$, moreover is isomorphic to $\bold J$.
\ermn
We now define by induction on $\zeta < \lambda$ an FSI-templates ${\frak t}_\zeta$ 
such that
\mr
\item "{$(*)^1_\zeta$}"  the set of members of ${\frak t}_\zeta$ is a set
of finite sequences starting with $\zeta$ hence disjoint to 
${\frak t}_\varepsilon$ for $\varepsilon < \zeta$; for $x \in {\frak t}_\zeta$
let $\xi(x) = \zeta$.
\endroster
\enddemo
\bn
\ub{Defining ${\frak t}_\zeta$}: \ub{Case 1}:  $\zeta = 0,\zeta$ successor
or cf$(\zeta) = \aleph_0$.

Let $L^{{\frak t}_\zeta} = \{\langle \zeta \rangle\}$ and 
$I^{{\frak t}_\zeta}_{<\zeta>} = \{\emptyset\}$.
\bn
\ub{Case 2}:  cf$(\zeta) > \aleph_0$.

Let $h_\zeta:\bold J \rightarrow \zeta$ be a function such that:
$\varepsilon < \zeta \Rightarrow h^{-1}_\zeta\{\varepsilon\}$ is a
dense subset of $\bold J$.  
The set of elements of ${\frak t}_\zeta$ is

$$
\{\langle \zeta \rangle\} \cup \{\langle \zeta \rangle \char 94 
\langle \eta \rangle 
\char 94 x:\eta \in \bold J \text{ and } x \in
\dbcu_{\varepsilon \le h_\zeta(\eta)} L^{{\frak t}_\varepsilon}\}
$$
\mn
order $<_{{\frak t}_\zeta}$ defined by:

$$
\langle \zeta \rangle \text{ is maximal}
$$

$$
\align
\langle \zeta \rangle \char 94 \langle \eta_1 \rangle 
\char 94 x_1 <_{{\frak t}_\zeta}
\langle \zeta \rangle \char 94 \langle \eta_2 \rangle 
\char 94 x_2 \text{ \ub{iff} } \eta_1
<_{\bold J} \eta_2 \vee (\eta_1 = \eta_2 &\and \xi(x_1) < \xi(x_2)) 
\vee (\eta_1 = \eta_2 \\
  &\and \xi(x_1) = \xi(x_2) \and x_1 <_{{\frak t}_{\xi(x_1)}} x_2).
\endalign
$$ 
\mn
Lastly, for $y \in {\frak t}_\zeta$ we define the ideal $I = 
I^{{\frak t}_\zeta}_y$:
\mr
\item "{$(\alpha)$}"  if $y = \langle \zeta \rangle$ then \nl
$I = \bigr\{ Y:Y \subseteq L^{{\frak t}_\zeta} \backslash 
\{ \langle \zeta \rangle\}\}$
\sn
\item "{$(\beta)$}"  if $y = \langle \zeta \rangle \char 94 \langle
\nu \rangle \char 94 x$,
\ub{then} $I$ is the family of sets $Y$ satisfying the following conditions: \nl

$\qquad (i) \quad Y \subseteq L^{{\frak t}_\zeta}$ \nl
\smallpagebreak

$\qquad (ii) \quad (\forall z \in Y)(z <_{{\frak t}_\zeta} y)$ \nl
\smallpagebreak

$\qquad (iii) \quad$ we each $\eta \in \bold J$ and $\xi < \zeta$
we have: \nl
\smallpagebreak

$\hskip60pt \{z:\langle \zeta \rangle \char 94 \langle \eta \rangle \char 94 
z \in Y$ and $\xi(z) = \xi$ and $z \ne \langle \xi \rangle\} 
\in I^{{\frak t}_\xi}_{\langle \xi \rangle}$ \nl
\smallpagebreak

$\qquad (iv) \quad$ the set $\{\eta \in \bold J:
(\exists x)(\langle \zeta \rangle \char 94 \langle \eta \rangle 
\char 94 x \in Y)\}$ is finite.
\ermn
Why is ${\frak t}_\zeta$ really a FSI-template?  We prove, of course, by
induction on $\zeta$ that:
\mr
\item "{$(*)^2_\zeta$}"  $(i) \quad L^{{\frak t}_\zeta}$ is a linear order
\sn
\item "{${{}}$}"  $(ii) \quad I^{{\frak t}_\zeta}_t$ is an ideal of subsets of
$\{s \in I^{{\frak t}_\zeta}_t:s < t\}$
\sn
\item "{${{}}$}"  $(iii) \quad {\frak t}_\zeta$ is an FSI-template,
\sn
\item "{${{}}$}"  $(iv) \quad {\frak t}_\zeta$ is disjoint to ${\frak
t}_\varepsilon$ for $\varepsilon < \zeta$.   
\ermn
[Why?  By \scite{ad.10}.] \nl
Next we prove by induction on $\zeta$, that ${\frak t}_\zeta$ is a smooth
FSI-template. Arriving to $\zeta$
\mr
\item "{$(*)^3_\zeta$}"  for $\eta \in \bold J$ and $\varepsilon \le h_\zeta(\eta)+1$, we
have ${\frak t}_\zeta \restriction \{\langle \zeta \rangle \char 94
\langle \eta \rangle \char 94 \rho:\rho \in \dbcu_{\xi < \varepsilon} 
{\frak t}_\xi\}$ is a
smooth FSI-template. \nl
[Why?  We prove by induction on $\varepsilon$; for $\varepsilon =0$ by 
\scite{ad.10}(1), for $\varepsilon$ successor by \scite{ad.10}(3) 
for $\varepsilon$ limit by \scite{ad.10}(5) and \scite{ad.10}(6)]
\sn
\item "{$(*)^4_\zeta$}"  for $Z \subseteq \bold J$ we have ${\frak t}_\zeta
\restriction (\dbcu_{\eta \in Z} \{\langle \zeta,\eta \rangle \char 94
\rho:\rho \in \dbcu_{\xi < h_\zeta(\eta)} t_\xi\})$ is a smooth FSI-template. \nl
[Why?  By induction on $|Z|$, for $|Z| = 0,|Z| = n+1$ by \scite{ad.10}(3),
for $|Z| \ge \aleph_0$ by \scite{ad.10}(5)]
\sn
\item "{$(*)^5_\zeta$}"  ${\frak t}_\zeta \restriction (L^{{\frak t}_\zeta}
\backslash \{ \langle \zeta \rangle\})$ is a smooth FSI-template. \nl
[Why?  By $(*)^4_\zeta$ for $Z = \bold J$.]
\sn
\item "{$(*)^6_\zeta$}"  ${\frak t}_\zeta$ is a smooth FSI-template \nl
[Why? by \scite{ad.10}(3)]
\sn
\item "{$(*)^7_\zeta$}"  if $K \subseteq L^{{\frak t}_\zeta}$ and $t
\in L^{{\frak t}_\zeta}$ then the ideal $I^{{\frak t}_\zeta}_t \cap
{\Cal P}(K)$ is generated by a countable family of subsets of $\kappa$
\nl
[Why?  Check by induction on $\zeta$.]
\ermn
Now for $\zeta \le \lambda$ let
\mr
\item "{$\boxtimes$}"  ${\frak s}_\zeta =: \dsize \sum_{\varepsilon < \zeta}
{\frak t}_\varepsilon$, i.e.
{\roster
\itemitem{ $(i)$ }  the set of elements of ${\frak s}_\zeta$ is
$\dbcu_{\varepsilon < \zeta} L^{{\frak t}_\varepsilon}$
\sn
\itemitem{ $(ii)$ }  for $x,y \in {\frak s}_\zeta$ we have 
$x <_{{\frak s}_\zeta} y$ iff $\xi(x) < \xi(y) \vee (\xi(x) = \xi(y) \and
x <_{{\frak t}_\zeta} y)$
\sn
\itemitem{ $(iii)$ }  $I^{{\frak s}_\zeta}_y = \{Y \subseteq {\frak s}_\zeta:
(\forall z \in Y)(z <_{{\frak s}_\zeta} y)$ and $\{z \in {\frak s}_\zeta:
\xi(z) = \xi(y)$ and \nl

$\hskip80pt  z \in Y\} \in I^{{\frak t}_{\xi(z)}}_y\}$
\endroster}
\item "{$(*)^8_\zeta$}"  ${\frak s}_\zeta$ is a smooth FSI-template. \nl
[Why?  Just easier than the proof above.]
\sn
\item "{$(*)^9_\zeta$}"  if $K \subseteq L^{{\frak s}_\zeta}$ is
countable and $t \in L^{{\frak s}_\zeta}$, \ub{then} the ideal
$I^{{\frak s}_\zeta}_t \cap {\Cal P}(K)$ of subsets of $K$ is generated
by a countable family of subsets of $K$ \nl
[Why?  By $(*)^7_\zeta$ and the definition of ${\frak s}_\zeta$.]
\endroster
\bn
Let \footnote{but if you like to avoid using $(*)^7_\zeta,(*)^9_\zeta$
and ${\Cal W}$ below just use $\theta = \beth^+_2$.  In fact even
without $(*)^7_\zeta + (*)^9_\zeta$ above, countable ${\Cal W}$ suffice
but then we have to weaken the notion of isomorphisms, and no point}
$\theta = (2^{\aleph_0})^+$, we shall prove below by induction on 
$\zeta$ that ${\frak s}_\zeta,{\frak t}_\zeta$ are $(\lambda,\theta)$-good 
(see definition below and Subclaim \scite{abd.2a}) then we can finish the proof as in
\scite{ad.15} using ${\frak s}_\mu$ (and $(*)^7_\zeta + (*)^9_\zeta$) where 
\definition{\stag{abd.2} Definition}  1) We say ${\frak t}$ is 
$(\lambda,\theta,\tau)$-good if:
\mr
\item "{$\oplus$}"  assume that $t_{\alpha,n} \in L^{\frak t}$ for $\alpha <
\theta,n < \omega,\{t_{\alpha,n}:n < \omega\}$ is $\bar K$-closed 
and ${\Cal W}$ a family of subsets of $\omega$ such that
$2^{|{\Cal W}|} < \theta$, \ub{then} we can find a club $C$ of
$\theta$ and a pressing down function $h$ on $C$ such that:
\sn
\item "{$\oplus'$}"  if $S \subseteq C$ is stationary in $\theta,(\forall 
\delta \in S)[\text{cf}(\delta) > \aleph_0]$ and 
$h \restriction S$ is constant \ub{then}:
{\roster
\itemitem{ $\boxtimes^1_S$ }  for 
every $\alpha < \beta$ in $S$ the truth value
of the following statements does not depend on $(\alpha,\beta)$: \nl
(but may depend on $n,m$ and $w \in {\Cal W}$) \nl
\smallskip

$\quad (i) \qquad t_{\alpha,n} = t_{\beta,m}$
\medskip

$\quad (ii) \qquad t_{\alpha,n} <_{L^{\frak t}} t_{\beta,m}$
\medskip

$\quad (iii) \qquad \{t_{\alpha,\ell}:
\ell \in w\} \in I^{\frak t}_{t_{\alpha,m}}$ 
\medskip

$\quad (iv) \qquad \{t_{\beta,\ell}:\ell \in w\} \in I^{\frak
t}_{t_{\alpha,n}}$
\medskip

$\quad (v) \qquad \{t_{\alpha,\ell}:\ell \in w\} \in I^{\frak
t}_{t_{\beta,n}}$
\sn
\itemitem{ $\boxtimes^2_S$ }  let $\delta^* \le \theta$, cf$(\delta^*)
= \tau,sup(S \cap \delta^*) = \delta^*$;
if $\theta \le \beta^* < \lambda$ and
$s_{\beta,n} \in L^{\frak t}$ for $\beta < \beta^* <
\lambda,n < \omega$ \ub{then} we can find $t_n \in L^{\frak t}$ for $n < 
\omega$ such that $\{t_n:n < \omega\}$ is $\bar K$-closed and
for every $\beta < \beta^*$, for every large enough 
$\alpha \in S \cap \delta^*$ for some ${\frak t}$-partial isomorphism $f$ we have 
$f(t_n) = t_{\alpha,n},f(s_{\beta,n}) = s_{\beta,n}$.
\endroster}
\ermn
2) We say ${\frak t}$ is strongly $(\lambda,\theta,\tau)$-good if
above we allow ${\Cal W} = {\Cal P}(\omega)$ (so if $\theta > \beth_2$
this is the same).  In both cases we may omit $\tau$ if $\tau = \theta$.
\enddefinition
\bn
\stag{abd.02.0}\ub{Observation}:  Instead ``$h$ regressive" it is enough to demand:
for some sequence $\langle X_\alpha:\alpha < \theta \rangle$ of sets, 
increasing continuous, $|X_\alpha| < \theta$ and for every (or club
of) $\delta < \theta$, if cf$(\delta) > \aleph_0$ then $h(\delta) \in 
{\Cal H}_{< \aleph_0}(X_\delta)$.
\bigskip

\proclaim{\stag{abd.2a} Subclaim}  In the proof of \scite{abd.1};
\mr
\widestnumber\item{$(iii)$}
\item "{$(i)$}"  ${\frak t}_\zeta$ is 
$(\lambda,\theta)$-good
\sn
\item "{$(ii)$}"  ${\frak s}_\zeta$ is strongly
$(\lambda,\theta,\aleph_1)$-good
\sn
\item "{$(iii)$}"   if cf$(\zeta) \ne \theta$ then also 
${\frak s}_\zeta$ is strongly $(\lambda,\theta)$-good.
\endroster
\endproclaim
\bigskip

\demo{Proof}  Recall that $\theta = (2^{\aleph_0})^+$, and let ${\Cal
W}$ be given ($2^{|{\Cal W}|} < \theta$ for the first version; ${\Cal
W} = {\Cal P}(\omega)$ for the second, using $(*)^7_\zeta +
(*)^9_\zeta$ from the proof of \scite{abd.1}).  We prove this by induction on $\zeta$.
\bn
\ub{For ${\frak s}_\zeta$}:

If $\zeta=0$ it is empty.  Otherwise given $t_{\alpha,n} \in {\frak s}_\zeta
= \dsize \sum_{\varepsilon < \zeta} {\frak t}_\varepsilon$ for $\alpha <
\theta,n < \omega$ let $h^*_0(\delta)$ be the sequence consisting of:
$\xi_{\alpha,n} =: \text{ Min}\{\xi:\xi \in \{\xi(t_{\beta,m}):\beta <
\delta,m < \omega\} \cup \{\infty\}$ and $\xi \ge \xi(t_{\alpha,n})\}$ 
for $n < \omega$ and $u_\alpha = \{(n,m):\xi(t_{\alpha,n}) = 
\xi_{\alpha,m})\}$ and $\bold w_\alpha = \{(n,w):n < \omega,w \in {\Cal W} \text{ and }
\{t_{\alpha,m}:m \in w\} \in I^{\frak t}_{t_{\alpha,n}}\}$; that is
$h^*_0(\delta) = \langle u_\alpha,\langle \xi_{\alpha,n}:n < \omega
\rangle,\bold w_\alpha \rangle$.  If $S_y = \{\delta:
\text{cf}(\delta) \ge \aleph_1,h^*_0(\delta) = y\}$ is stationary
we define $h^*_1 \restriction S_y$ such that it codes $h^*_0(\delta)$ and if
$n(*) < \omega,\alpha \in S_y \Rightarrow \xi(t_{\alpha,n(*)}) =
\xi_{\alpha,n(*)}$ call it $\xi_{y,n(*)}$ let $u_{y,n(*)} = \{n:
\xi_{\alpha,n} = \xi_{y,n(*)}\}$, \ub{then} $h_1 \restriction S_y$ codes a
function witnessing the $(\lambda,\theta)$-goodness of 
${\frak t}_{\xi_{y,n(*)}}$ for $\langle t_{\alpha,n}:n \in u_{y,n(*)},
\alpha \in S_y \rangle$.
\bn
It is easy to check that this shows $\boxtimes^1_S$ even if cf$(\zeta) = \theta$.
But assume cf$(\zeta) \ne \theta \and \delta^* = \theta$ or $\delta^*
< \theta$, cf$(\delta^*) = \aleph_1$ (or just $\aleph_0 < \text{
cf}(\delta^*) < \theta$), $\delta^* = \sup(S \cap \delta^*)$; 
we shall prove also the statement from $\boxtimes^2_S$.  Let
$w_1 = \{n:\langle \xi(t_{\beta,n}):\beta \in S \rangle$ is strictly
increasing$\}$, $w_0 = \{n:\langle \xi(t_{\beta,n}):\beta \in S \rangle$
is constant$\}$, let $\xi(S,n) = \xi_{S,n} = \cup \{\xi(t_{\beta,n}):\beta \in
S\}$ as cf$(\zeta) \ne \theta$ it is $< \zeta$.

Given $\bar s = \langle s_{\beta,n}:n < \omega \rangle$ we have to find
$\langle t_n:n < \omega \rangle$ as required in $\boxtimes^2_S$.  If
$n \in w_0,w'_{0,n} = \alpha\{m \in w_0:\xi(t_{\alpha,n}) =
\xi(t_{\alpha,m})$ for $\alpha \in S\}$ and to choose $\langle t_m:m
\in w'_{0,n} \rangle$ we use the induction hypothesis on ${\frak t}_{\xi(S,n)}$.  If
$n \in w_1$ then we can find $t^*_n \in {\frak t}_{\xi_{S,n}}$ such that
$\{t:t \in {\frak t}_{\xi_{S,n}},t \le_{{\frak t}_{\xi(S,n)}} t^*\}$ is
disjoint to $\{t_{\beta,m}:\beta < \delta^*,m < \omega\} \cup \{s_{\beta,m}:
\beta < \beta^*$ and $m < \omega\}$ because the lower cofinality of
$L^{{\frak t}_{\xi(S,n)}}$ is the same as that of $L_0$ and is $\lambda >
\theta + |\beta^*|$.  We choose $\eta^* \in \bold J$ such that $(\forall x)
(\langle \zeta,\eta^* \rangle \char 94 \langle x \rangle \in
{\frak t}_{\xi(S,n)} \Rightarrow t <_{{\frak t}_{\xi(S,n)}} t^*)$ and we
choose together $\langle t_{n'}:n' \in w_1,\xi_{S,n'} = \xi_{S,n} \rangle$
taking care of ${\Cal W}$, (inside and automatically for others,
i.e. considering $t_{n_1},t_{n_2}$ such that $\xi_{S,n_1} \ne
\xi_{S,n_2})$, this is immediate.
\enddemo
\hfill$\square_{\scite{abd.1}}$\margincite{abd.1}
\bn
\ub{For ${\frak t}_\zeta$}:

Similar.  \hfill$\square_{\scite{abd.2a}}$\margincite{abd.2a}
\bn
\centerline {$* \qquad * \qquad *$}
\bn
We may like to have ``$2^{\aleph_0} = \lambda$ is singular", ${\frak a} =
\lambda,{\frak b} = {\frak d} = \mu$.  Toward this we would like to have
a linear order $\bold J$ such that if $\bar x = \langle x_\alpha:\alpha <
\kappa \rangle$ is monotonic, say decreasing then for any $\sigma < \lambda$
for some limit $\delta < \kappa$ of uncountable cofinality the linear order
$\{y \in \bold J:\alpha < \delta \Rightarrow y <_{\bold J} x_\alpha\}$ has
cofinality $> \sigma$.  Moreover, $\delta$ can be chosen to suit $\omega$
such sequences $\bar x$ simultaneously.  So every set of $\omega$-tuples from
$\bold J$ of cardinality $\ge \kappa$ but $< \lambda$ can be ``inflated".
\bigskip

\proclaim{\stag{abd.3} Lemma}  Assume
\mr
\item "{$(a)$}"  $(2^{\aleph_0})^+ < \mu = \text{ cf}(\mu) \le \tau
< \lambda = \lambda^{\aleph_0},\lambda$ singular
\sn
\item "{$(b)$}"  $(\forall \alpha < \tau)[|\alpha|^{\aleph_0} < \mu =
\text{ cf}(\tau)]$
\sn
\item "{$(c)$}"  $\tau \ge \aleph_{\text{cf}(\lambda)}$ or at least
\sn
\item "{$(c)^-$}"  there is $f:\lambda \rightarrow \text{ cf}(\lambda)$
such that if $\langle \alpha_\varepsilon:\varepsilon < \tau \rangle$ is
strictly increasing continuous, $\alpha_\varepsilon < \lambda$ and $\gamma
< \text{ cf}(\lambda)$ \ub{then} for some $\varepsilon < \tau$ we have
$f(\alpha_\varepsilon) \ge \gamma$.
\ermn
\ub{Then} for some c.c.c. forcing notion of cardinality $\lambda$ we have
$\Vdash_P ``2^{\aleph_0} = \lambda,{\frak b} = {\frak d} =
\kappa,{\frak a} = \lambda"$.
\endproclaim
\bigskip

\demo{Proof}  Note that $(c) \Rightarrow (c)^-$, just let $\alpha < \lambda
\and \text{ cf}(\alpha) = \aleph_\varepsilon \and \varepsilon <
\text{ cf}(\lambda) \Rightarrow f(\alpha) = \varepsilon$, clearly there is
such a function and it satisfies clause $(c)^-$.  So we can assume $(c)^-$.  Let
$\theta = (2^{\aleph_0})^+$ and $\sigma = \text{ cf}(\lambda)$ and
$\langle \lambda_\varepsilon:\varepsilon < \sigma \rangle$ be a strictly
increasing sequence of regular cardinals $> \tau$ with limit $\lambda$.  Let
$\langle L_{0,\gamma}:\gamma < \text{ cf}(\lambda) \rangle$ be increasing
with $\gamma,L_{0,\gamma}$ like $L_0$ in the proof of \scite{abd.1} with
$\lambda_\varepsilon$ instead of $\lambda$, such that $\beta < \gamma
\Rightarrow L_{0,\beta}$ is an interval of $L_{0,\gamma}$.  Let
$L_0 = \dbcu_{\gamma < \text{cf}(\lambda)} L_{0,\gamma}$ define $g:L_0
\rightarrow \text{ cf}(\lambda)$ by $g(x) = \text{
Min}\{\gamma:x \in L_{0,\gamma}\}$ and let

$$
\align
\bold J^* = \bigl\{ \eta \in {}^{\omega >}(L_0):\eta(0) \in L_{0,0}
\text{ and } \eta(n+1) \in L_{0,g(\eta(n))} \text{ for } n < \omega \bigr\}
\endalign
$$
\mn
ordered as in the proof of \scite{abd.3}.

We define ${\frak s}_\zeta,{\frak t}_\zeta$ as there.  We then prove that
${\frak s}_\zeta,{\frak t}_\zeta$ are $(\tau,\theta)$-good and 
$(\lambda,\tau)$-good as there and this suffices repeating the proof of
\scite{abd.1}. \hfill$\square_{\scite{abd.4}}$\margincite{abd.4}
\enddemo
\bn
\ub{\stag{abd.3a} Discussion}:  We may 
like to separate ${\frak b}$ and ${\frak d}$.  So
below we adapt the proof of \scite{abd.1} to do this (can do it also for
\scite{abd.3}).

A way to do this is to look at 
the forcing in \scite{abd.1} as the limit of the FS iteration
$\langle P^*_i,{\underset\tilde {}\to Q^*_j}:i \le \mu,j < \mu \rangle$, so
the memory of $Q^*_j$ is $\{i:i < j\}$ where
${\underset\tilde {}\to Q^*_j}$ is Lim$_{\frak t}[\langle Q_t:t \in
L^{{\frak t}_j} \rangle]$.  Below we will use the limit of FS iteration
$\langle P^*_i,{\underset\tilde {}\to Q^*_j}:j < \mu \times \mu_1 \rangle,
Q^*_\zeta$ has memory $w_\zeta \subseteq \zeta$ where e.g. for $\zeta
= \mu \,\alpha +i,w_\zeta = \{\kappa \beta +j:\beta \le \alpha,j \le i,(\beta,j) \ne
(\alpha,i)\}$.  Let $P^* = P^*_{\mu \times \mu_i}$ be $\cup\{P_i:i <
\mu \times \mu_1\}$.

Of course, $Q_\zeta$ will be defined as Lim$_{{\frak t}_\zeta}(\bar Q)$,
the ${\frak t}_\zeta$ defined as above and ${\frak b} = \mu,{\frak d} =
\mu_1$.  Should be easy.  If $\langle {\underset\tilde {}\to A_\varepsilon}:
\varepsilon < \varepsilon^{\bar x} \rangle$ exemplifies ${\frak a}$ in
$\bold V^{P^*}$, so $\varepsilon^* \ge \mu$ then for some
$(\alpha^*,\beta^*) \in \mu \times \mu_1$ for $\kappa (= \theta)$ of the
names they involve $\{{\underset\tilde {}\to Q_{\mu \alpha + \beta}}:
\alpha \le \alpha^*,\beta \le \beta^*\}$ only.

Using indiscernibility on the pairs $(\alpha,\beta)$ to making them increase
we can finish.
\bigskip

\proclaim{\stag{abd.4} Lemma}  1) In Lemma \scite{abd.1}, if $\mu =
\text{ cf}(\mu) \le \text{ cf}(\mu_1),\mu_1 < \lambda$, \ub{then} we can 
change in the conclusion ${\frak b} = {\frak d} = 
\mu$ to ${\frak b} = \mu,{\frak d} = \mu_1$. \nl
2) Similarly for \scite{abd.3}.
\endproclaim
\bigskip

\demo{Proof} \ub{First Proof}:  If $\mu_1$ regular, let $\mu_0 =
\mu$.  The proof of \scite{abd.1} for $\ell \in \{0,1\}$ using
$\mu = \mu_\ell$ gives ${\frak s}^\ell_{\mu_\ell}$ and \wilog \,
${\frak s}^0_{\mu_0},{\frak s}^1_{\mu_1}$ are disjoint.  Let 
${\frak s}$ be ${\frak s}_0 +' {\frak s}_1$ meaning
$L[{\frak s}] = L[{\frak s}^0_{\mu_0}] + L[{\frak s}^1_{\mu_1}]$, and
for $t \in L[{\frak s}^\ell_{\mu_\ell}]$ we let $I^{\frak s}_t =:
I^{{\frak s}^\ell_{\mu_\ell}}_t$ (this is not ${\frak s}_0 + {\frak s}_1$
of \scite{ad.10}.  Now the appropriate goodness can be proved. 
\sn
\ub{Second Proof}:  Instead of starting with $\langle Q_i:i < \mu
\rangle$ with full memory we start with 
$\langle {\underset\tilde {}\to Q_\zeta}:\zeta < \mu
\times \mu_1 \rangle,{\underset\tilde {}\to Q_\zeta}$ 
with memory if $\zeta = \mu \alpha +i,
i < \kappa,w_\zeta = \{\mu \beta +j:\beta \le \alpha,j \le i,(\beta,j) \ne
(\alpha,i)\}$.  \hfill$\square_{\scite{abd.4}}$\margincite{abd.4}
\enddemo
\newpage

\head {\S4  On related cardinal invariants} \endhead  \resetall \sectno=4
\bigskip

\proclaim{\stag{au.1} Theorem}  Assume
\mr
\item "{$(a)$}"  $\kappa$ is a measurable cardinal
\sn
\item "{$(b)$}"  $\kappa < \mu = \text{ cf}(\mu) < \lambda =
\text{ cf}(\lambda) = \lambda^\kappa$.
\ermn
\ub{Then} for some c.c.c. forcing notion $P$ of cardinality $\lambda$, in
$\bold V^P$ we have: $2^{\aleph_0} = \lambda,{\frak u} = {\frak d} = {\frak b} =
\mu$ and ${\frak a} = \lambda$.
\endproclaim
\bigskip

\remark{Remark}  Recall ${\frak u} = \text{ Min}\{|{\Cal P}|:
{\Cal P} \subseteq [\omega]^{\aleph_0}$ generates a nonprincipal 
ultrafilter on $\omega\}$.
\endremark
\bigskip

\demo{Proof}  The proof is broken to definitions and claims.
\enddemo
\bigskip

\definition{\stag{au.1a} Definition}  For an ultrafilter $D$ on $\omega$
let $Q(D)$ be: \nl
$\{T:T \subseteq {}^{\omega >} \omega$ is
closed under initial segments, and for some tr$(T)$, the trunk of $T$, we
have:
\mr
\item "{$(i)$}"   $\ell \le \ell g(tr(T)) \Rightarrow T \cap {}^\ell
\omega = \{tr(T) \restriction \ell\}$
\sn
\item "{$(ii)$}"   tr$(T) \trianglelefteq \eta \in {}^{\omega >} \omega
\Rightarrow \{n:\eta \char 94 \langle n \rangle \in T\} \in
{\underset\tilde {}\to D_i}\}$ 
\ermn
ordered by inverse inclusion.
\enddefinition
\bigskip

\definition{\stag{au.1b} Definition}
Let ${\frak K}$ be the family of ${\frak t}$ 
consisting of $\bar Q =
\bar Q^{\frak t} = \langle P_i,{\underset\tilde {}\to Q_i}:i < \mu \rangle
= \langle P^{\frak t}_i,{\underset\tilde {}\to Q^{\frak t}_i}:i < \mu \rangle$
and $\bar D = \bar D^{\frak t} = \langle {\underset\tilde {}\to D_i}:i <
\mu$ and cf$(i) \ne \kappa \rangle$ and $\bar \tau^{\frak t} = \langle
{\underset\tilde {}\to \tau^{\frak t}_i}:i < \mu \rangle$ such that
\mr
\item "{$(a)$}"  $\bar Q$ is a FS-iteration of c.c.c. forcing notions (and
$P_\mu = P^{\frak t}_\mu = \text{ Lim}(\bar Q^{\frak t}) =
\dbcu_{i < \mu} P^{\frak t}_i)$
\sn
\item "{$(b)$}"  if $i < \mu$, cf$(i) \ne \kappa$ then $Q_i = Q
({\underset\tilde {}\to D_i})$, \nl
see Definition \scite{au.1a} above
\sn 
\item "{$(c)$}"  ${\underset\tilde {}\to D_i}$ is a $P_i$-name of a
nonprincipal ultrafilter on $\omega$ when $i < \mu$, cf$(i) \ne \kappa$
\sn
\item "{$(d)$}"  $|P_i| \le \lambda$
\sn
\item "{$(e)$}"  for $i < \mu$, cf$(i) \ne \kappa$ let
${\underset\tilde {}\to \eta_i}$ be the $P_{i+1}$-name of the
${\underset\tilde {}\to Q_i}$-generic real

$$
{\underset\tilde {}\to \eta_i} = \cup\{tr(p(i)):p \in
{\underset\tilde {}\to G_{P_{i+1}}}\}.
$$

Then for $i <j$ from $< \mu$ of cofinality $\ne \kappa$ we have

$$
\Vdash_{P_j} ``\text{Rang}({\underset\tilde {}\to \eta_i}) \in
{\underset\tilde {}\to D_j}"
$$
\sn
\item "{$(f)$}"  ${\underset\tilde {}\to \tau_i}$ is a $P_i$-name of a
function from ${\underset\tilde {}\to Q_i}$ to $\{h:h$ a function from a 
finite set of ordinals to ${\Cal H}(\omega)\}$, such that: \nl
$\Vdash_{P_i}$ ``if $p,q \in {\underset\tilde {}\to Q_i}$ 
are compatible \ub{then} they have a common
upper bound $r$ such that ${\underset\tilde {}\to \tau_i}(r) =
{\underset\tilde {}\to \tau_i}(p) \cup {\underset\tilde {}\to \tau_i}(q)"$
\sn
\item "{$(g)$}"  if cf$(i) \ne \mu$ and $i \in \text{ Dom}(p),p \in
P_j$ and $i < j \le \mu$ then 
${\underset\tilde {}\to \tau_i}(p(i))$ is $\{\langle 0,tr(p) \rangle\}$; i.e. this is
forced to hold 
\sn
\item "{$(h)$}"  we stipulate $P_i = \{p:p$ is a function with domain $a$
finite subset of $i$ for each $j \in \text{ Dom}(p),\emptyset_{P_j}$ forces 
that $p(j) \in {\underset\tilde {}\to Q_j}$ and it forces a value to 
${\underset\tilde {}\to \tau_j}(p(j))\}$
\sn
\item "{$(i)$}"  $\Vdash_{P_i} ``{\underset\tilde {}\to Q_i} \subseteq
{\Cal H}_{< \aleph_0}(\gamma)$ for some ordinal $\gamma"$.
\ermn
Let $\gamma({\frak t})$ be the minimal ordinal $\gamma$ such that $i <
\mu \Rightarrow \Vdash_{P_i}$ ``if $x \in {\underset\tilde {}\to Q_i}$
then dom$({\underset\tilde {}\to \tau_i}(x)) \subseteq \gamma$". \nl
We let $\tau^{{\frak t},i}$ be the function with domain $P_i$ such
that $\tau^{{\frak t},i}(p)$ is a function with domain 
$\{\gamma({\frak t})j + \beta:j \in \text{ Dom}(p)$ and
$p \restriction j \Vdash_{P_j} ``\beta \in \text{
Dom}({\underset\tilde {}\to \tau_j}(p(j))"\}$ and let
$\tau^{{\frak t},i}(\gamma({\frak t})j + \beta)$ be the value 
which $p \restriction j$ forces
on ${\underset\tilde {}\to \tau^{\frak t}_j(\beta)}$.
\enddefinition
\bn
Obviously \nl
\ub{\stag{au.A} Subclaim}:  ${\frak K} \ne \emptyset$.
\bigskip

\demo{Proof}  Should be clear.
\enddemo
\bn
\ub{\stag{au.B} Subclaim}:  if in a universe $\bold V,D$ 
is a nonprincipal ultrafilter on $\omega$
then
\mr
\item "{$(a)$}"  $\Vdash_{Q(\bar D)} ``\{tr(p)(\ell):\ell < \ell g(tr(p)),
p \in {\underset\tilde {}\to Q_{Q(D)}}\}$ is an infinite subset of $\omega$,
almost included in every member of $D$"
\sn
\item "{$(b)$}"  $Q(D)$ is a c.c.c. forcing notion, even $\sigma$-centered
\sn
\item "{$(c)$}"  ${\underset\tilde {}\to \eta_i} =
\cup\{tr(p):p \in {\underset\tilde {}\to G_{Q(D)}}\} \in
{}^\omega \omega$ is forced to dominate $({}^\omega \omega)^{\bold V}$.
\ermn
[Note that this, in particular clause (c), does not depend on
additional properties of $D$; but as we naturally add many Cohen reals
(by the nature of the support) we may add more and then can demand
e.g. ${\underset\tilde {}\to D_i}$ (cf$(i) \ne \kappa$) is a Ramsey ultrafilter.]
\bigskip

\definition{\stag{au.Ba} Definition}  1)  We define 
$\le_{\frak K}$ by: ${\frak t} \le_{\frak K} {\frak s}$ if 
(${\frak t},{\frak s} \in {\frak K}$ and) $i \le \mu \Rightarrow P^{\frak t}_i
\lessdot P^{\frak s}_i$ and $i < \mu \and \text{ cf}(i) \ne \kappa \Rightarrow
\Vdash_{P^{\frak s}_i} ``{\underset\tilde {}\to D^{\frak t}_i} \subseteq
{\underset\tilde {}\to D^{\frak s}_i}"$ and $i < \mu \Rightarrow
\Vdash_{P^{\frak s}_i} ``{\underset\tilde {}\to \tau^{\frak t}_i} \subseteq
{\underset\tilde {}\to \tau^{\frak s}_i}"$.
\nl
2) We say ${\frak t}$ is a canonical
$\le_{\frak K}$-ub of $\langle {\frak t}_\alpha:\alpha < \delta \rangle$ if:
\mr
\widestnumber\item{$(iii)$} 
\item "{$(i)$}"  ${\frak t},{\frak t}_\alpha \in {\frak K}$
\sn
\item "{$(ii)$}"  $\alpha \le \beta < \delta \Rightarrow {\frak t}_\alpha
\le_{\frak K} {\frak t}_\beta \le_{\frak K} {\frak t}$
\sn
\item "{$(iii)$}"  if $i < \mu$ and cf$(i) = \kappa$ then 
$\Vdash_{P^{\frak t}_i} ``{\underset\tilde {}\to Q^{\frak t}_i} =
\dbcu_{\alpha < \delta} {\underset\tilde {}\to Q^{{\frak t}_\alpha}_i}"$.
\ermn
Note that if cf$(\delta) > \aleph_0$ then
$\Vdash_{P^{\frak t}_i} ``{\underset\tilde {}\to Q^{\frak t}_i} =
\dbcu_{\alpha < \delta} 
{\underset\tilde {}\to Q^{{\frak t}_\alpha}_i}"$ for every $i < \mu$,
so ${\frak t}$ is totally determined. \nl
3) We say $\langle {\frak t}_\alpha:\alpha < \alpha^* \rangle$ is
$\le_{\frak K}$-increasing continuous if: $\alpha < \beta < \alpha^*
\Rightarrow {\frak t}_\alpha \le_{\frak K} t_\beta$ and for limit $\delta <
\alpha^*,{\frak t}_\delta$ is a canonical $\le_{\frak K}$-ub of 
$\langle {\frak t}_\alpha:\alpha < \delta \rangle$.  Note that we have
not say ``the canonical $\le_{\frak K}$-u.b." as for $\delta <
\alpha^*$, cf$(\delta) \ne \kappa$ we have some freedom in completing
$\cup\{{\underset\tilde {}\to D^{{\frak t}_\alpha}}:\alpha < \delta\}$
to an ultrafilter (on $\omega$ in $\bold V^{P^{\frak t}_i}$).
\enddefinition
\bn
\ub{\stag{au.C} Subclaim}:  If $P_1 \lessdot P_2$ and ${\underset\tilde {}\to D_\ell}$ 
is a $P_\ell$-name of a nonprincipal ultrafilter on $\omega$ for $\ell =1,2$
and $\Vdash_{P_2} ``{\underset\tilde {}\to D_1} \subseteq 
{\underset\tilde {}\to D_2}"$, \ub{then} $P_1 * Q({\underset\tilde {}\to D_1})
\lessdot P_2 * Q({\underset\tilde {}\to D_2})$. \nl
[Why?  First, we can first force with $P_1$, so \wilog \, $P_1$ is trivial and
$D_1 \in \bold V$ is a nonprincipal ultrafilter on $\omega$.  Now
clearly $p \in Q(D_1) \Rightarrow p \in 
Q({\underset\tilde {}\to D_2})$ and $Q(D_1) \models p \le q
\Rightarrow Q({\underset\tilde {}\to D_2}) \models p \le q$ and if $p,q \in 
Q(D_1)$ are incompatible in $Q(D_1)$ then they are incompatible in 
$Q({\underset\tilde {}\to D_2})$.  Lastly, in $\bold V$, let ${\Cal I} = \{p_n:n <
\omega\} \subseteq Q(D_1)$ be predense in $Q(D_1)$, we shall prove that
${\Cal I}$ is predense in $Q({\underset\tilde {}\to D_2})$.  For this it
suffices to note
\mr
\item "{$\boxtimes$}"  if $D_1$ is a nonprincipal ultrafilter on $\omega,
{\Cal I} \subseteq Q(D_1)$ and $\eta \in {}^{\omega >} \omega$, \ub{then} the
following conditions are equivalent
{\roster
\itemitem{ $(a)_\eta$ }  there is no $p \in Q(D_1)$ incompatible with every
$q \in {\Cal I}$ which satisfies tr$(p)=\eta$
\sn
\itemitem{ $(b)_\eta$ }  there is a set $T$ such that:
\sn
\itemitem{ ${{}}$ }  $(i) \qquad \nu \in T 
\Rightarrow \eta \trianglelefteq \nu \in p$
\sn
\itemitem{ ${{}}$ }  $(ii) \qquad \eta 
\trianglelefteq \nu \trianglelefteq \rho \in T
\Rightarrow \nu \in T$
\sn
\itemitem{ ${{}}$ }  $(iii) \qquad$ if 
$\nu \in T$ then either $\{n:\nu \char 94 \langle n
\rangle \in T\} \in D_1$ or \nl

$\qquad \qquad \qquad (\forall n)(\nu \char 94 \langle n \rangle \notin
T) \and (\exists q \in {\Cal I})(\nu = tr(q))$
\sn
\itemitem{ ${{}}$ }  $(iv) \qquad$ there 
is a strictly decreasing function $h:T \rightarrow \omega_1$
\sn
\itemitem{ ${{}}$ }  $(v) \qquad \,\,\, \eta \in p$.
\endroster}
\endroster
\bigskip

\demo{Proof of $\boxtimes$}  Straight.

So as in $\bold V,{\Cal I} \subseteq Q(D_1)$ is predense, for every
$\eta \in {}^{\omega >} \omega$ we have $(a)_\eta$ for $D_1$ hence by
$\boxtimes$ we have also $\eta \in {}^{\omega >} \omega \Rightarrow
(b)_\eta$, but clearly if $T_\eta$ witness $(b)_\eta$ in $\bold V$
for $D_1$, it witnesses $(b)_\eta$ in $\bold V^{P_2}$ for $D_2$ hence
applying $\boxtimes$ again we get: $\eta \in {}^{\omega >} \omega
\Rightarrow (a)_\eta$ in $\bold V^{P_2}$ for $D_2$, hence ${\Cal I}$
is predense in $Q(D_2)$ in $\bold V^{P_2}$.  So we have proved
Subclaim \scite{au.C}.
\enddemo
\bn
\ub{\stag{au.D} Subclaim}:  If $\langle {\frak t}_\alpha:\alpha < \delta \rangle$ is
$\le_{\frak K}$-increasing continuous and $\delta < \lambda^+$ is a limit
ordinal, \ub{then} it has a canonical $\le_{\frak K}$-ub. \nl
[Why?  By induction on $i < \mu$, we define $P^{\frak t}_i$ and then
${\underset\tilde {}\to Q^{\frak t}_i},{\underset\tilde {}\to \tau_i}$
and ${\underset\tilde {}\to D_i}$ (if cf$(i) \ne \kappa$) such 
that the relevant demands (for ${\frak t} \in {\frak K}$ and for being 
canonical $\le_{\frak K}$-ub of $\bar{\frak t}$) hold.
\sn
Defining $P^{\frak t}_i$ is obvious:
for $i=0$ trivially, if $i = j+1$ it is $P^{\frak t}_j * Q^{\frak
t}_j$ and $i$ is limit it is $\cup\{P^{\frak t}_j:j < i\}$.  As we
have proved the relevant demands on $P^{\frak t}_j,{\underset\tilde
{}\to Q^{\frak t}_j}$ for $j<i$ clearly $P^{\frak t}_i$ is c.c.c. by using $\langle
{\underset\tilde {}\to \tau_j}:j < i \rangle$ and clearly $\langle
P^{\frak t}_\zeta,{\underset\tilde {}\to Q^{\frak t}_\xi}:\zeta \le i,\xi < i
\rangle$ is an FS iteration.  Now we shall prove that $\alpha < \delta \Rightarrow
P^{{\frak t}_\alpha}_i \lessdot P^{\frak t}_i$? \nl
So let ${\Cal I}$ be a predense subset of $P^{{\frak t}_\alpha}_i$ and
$p \in P^{\frak t}_i$ and we should prove that $p$ is compatible with some
$q \in {\Cal I}$ in $P^{\frak t}_i$; we divide the proof to cases.
\bn
\ub{Case 1}:  $i$ is a limit ordinal.

So $p \in P^{\frak t}_j$ for some $j < i$, let ${\Cal I}' = \{q \restriction
j:q \in {\Cal I}\}$, so clearly ${\Cal I}'$ is a predense subset of
$P^{{\frak t}_\alpha}_j$ (as ${\frak t}_\alpha \in {\frak K}$).  By the
induction hypothesis, in $P^{\frak t}_j$ the condition $p$ is compatible with
some $q' \in {\Cal I}'$; so let $r' \in P^{\frak t}_j$ be a common upper
bound of $q',p$ and let $q' = q \restriction j$ where $q \in {\Cal I}$.  So
$r \cup (q \restriction [j,i)) \in P^{\frak t}_i$ is a common upper bound of
$q,p$ as required.
\bn
\ub{Case 2}:  $i = j+1$, cf$(j) = \kappa$.

So \wilog \, for some $\beta < \delta,p(j)$ is a 
$P^{{\frak t}_\beta}_j$-name of a member 
of ${\underset\tilde {}\to Q^{{\frak t}_\beta}_j}$;
and \wilog \, $\alpha \le \beta < \delta$.  By the induction hypothesis
$P^{{\frak t}_\beta}_j \lessdot P^{\frak t}_j$ hence there is $p' \in
P^{{\frak t}_\beta}_j$ such that $[p' \le p'' \in P^{{\frak t}_\beta}_j
\Rightarrow p'',p \restriction j$ are compatible in $P^{\frak t}_j]$. \nl
Let

$$
\align
{\Cal J} = \{q' \restriction j:&q' \in P^{{\frak t}_\beta}_i 
\text{ and } q' \text{ is above some member of } {\Cal I} \\
  &\text{and } q' \restriction j
\Vdash_{P^{{\frak t}_\beta}_j} ``p(j) 
\le^{\underset\tilde {}\to Q^{{\frak t}_\beta}_j} q'(j)"\}.
\endalign
$$
\mn
Now ${\Cal J}$ is a dense subset of $P^{{\frak t}_\beta}_j$ (since if
$q \in P^{{\frak t}_\beta}_j \text{ then } q \cup \{\langle j,p(j) \rangle\}$
belongs to $P^{{\frak t}_\beta}_i$ hence is compatible with some member of
${\Cal I}$). \nl
Hence $p'$ is compatible with some $q'' \in {\Cal J}$, so there is $r$
such that $p' \le r \in
P^{{\frak t}_\beta}_j,q'' \le r$.  As $q'' \in {\Cal J}$ there is $q' \in
P^{{\frak t}_\beta}_i$ such that $q' \restriction j = q'',q'$ 
is above some $q^* \in {\Cal I}$ and $q' \restriction j
\Vdash ``p(j) 
\le^{\underset\tilde {}\to Q^{{\frak t}_\beta}_j} q'(j)"$.

As $P^{{\frak t}_\beta}_j \models ``p' \le r \and q' \restriction j =
q'' \le r"$ and by the
choice of $p'$ there is $p^* \in P^{\frak t}_j$ above $r$ (hence above
$p'$ and above $q'' = q' \restriction j)$,
and above $p \restriction j$.  Now let $r^* = p^* \cup (q'' \restriction \{j\})$,
clearly $r^* \in P^{\frak t}_i$ is above $p \restriction j$ and 
$r^* \restriction j$ forces that $r^*(j)$ is above
$p \restriction \{j\}$.  Clearly $r^* \restriction j$ is above $r$ and
$r^*$ also is above $q^* \in {\Cal I}$ so we are done.
\bn
\ub{Case 3}:  $i = j+1$, cf$(j) \ne \kappa$.

Use Subclaim C above. \nl
So we have dealt with $\alpha < \delta \Rightarrow P^{{\frak
t}_\alpha}_i \lessdot P^{\frak t}_i$.

If $P^{\frak t}_i$ has been defined and cf$(i) = \kappa$ we let
${\underset\tilde {}\to Q^{\frak t}_i} = \dbcu_{\alpha < \delta}
Q^{{\frak t}_\alpha}_i$ and ${\underset\tilde {}\to \tau^{\frak t}_i} =
\dbcu_{\alpha < \delta} {\underset\tilde {}\to \tau^{{\frak t}_\alpha}_i}$,
easy to check that they are as required.  If $P^t_i$ has been defined and
cf$(i) \ne \kappa$, then $\dbcu_{\alpha < \delta} D^t_i$ is a filter
on $\omega$ containing the co-bounded subsets, and we complete it to an
ultrafilter.  Clearly we are done.]
\bn
\ub{\stag{au.E} Subclaim}:  If ${\frak t} \in {\frak K}$ and $E$ is a $\kappa$-complete
nonprincipal ultrafilter on $\kappa$, \ub{then} we can find ${\frak s}$
such that:
\mr
\widestnumber\item{$(iii)$}
\item "{$(i)$}"   ${\frak t} \le_{\frak K} {\frak s} \in {\frak K}$
\sn
\item "{$(ii)$}"  there is $\langle \bold k_i,\bold j_i:
i < \mu$,cf$(i) \ne \kappa)$ such that:
{\roster
\itemitem{ $(\alpha)$ }  $\bold k_i$ is an isomorphism from
$(P^{\frak t}_i)^\kappa/E$ onto $P^{\frak s}_i$
\sn
\itemitem{ $(\beta)$ }  $\bold j_i$ is the canonical embedding of
$P^{\frak t}_i$ into $(P^{\frak t}_i)^\kappa/E$
\sn
\itemitem{ $(\gamma)$ }  $\bold k_i \circ \bold j_i =$ identity on
$P^{\frak t}_i$
\endroster}
\item "{$(iii)$}"  ${\underset\tilde {}\to D^{\frak s}_i}$ is the
image of
$({\underset\tilde {}\to D_i})^\kappa/E$ under $\bold k_i$ and similarly
${\underset\tilde {}\to \tau_i}$ \ub{if} $i < \mu$, cf$(i) \ne \kappa$
\sn
\item "{$(iv)$}"  if $i < \mu$, cf$(i) = \kappa$, then
${\underset\tilde {}\to \tau_i}$ is defined such that, for $j <
\kappa$, cf$(j) \ne \kappa$ we have $\bold k_j$ is an isomorphism from
$(P^{\frak t}_i,\gamma',\tau^{{\frak t},i})^\kappa/D$ onto $(P^{\frak
s}_i,\gamma'',\tau^{{\frak s},i})$ for some ordinals $\gamma',\gamma''$.
\ermn
[Why?  Straight.  Note that if cf$(i) = \kappa,i < \mu$ then 
${\underset\tilde {}\to Q^{\frak s}_i}$ is isomorphic to $P^{\frak s}
_{i+1}/P^{\frak s}_i$ which is c.c.c. as by \L o\'s theorem for the logic
$L_{\kappa,\kappa}$ we have $\dbcu_{j<i} (P^t_j)^\kappa/E \lessdot
(P^t_{i+1})^\kappa/E$, similarly for
${\underset\tilde {}\to \tau_i}$ which guarantees that the quotient 
is c.c.c., too (actually ${\underset\tilde {}\to \tau_i}$ is not needed for
the c.c.c. here).] 
\bn
\ub{\stag{au.F} Subclaim}  If ${\frak t} \in {\frak K}$ then
$\Vdash_{P^{\frak t}_\mu} ``{\frak u} = {\frak b} = {\frak d} = \mu"$. 
\sn
[Why?  In $\bold V^{P^{\frak t}_\mu}$, the family ${\Cal D} = 
\{\text{Rang}({\underset\tilde {}\to \eta_i}):
i < \mu$,cf$(i) \ne \kappa\} \cup \{[n,\omega):n < \omega\}$ generates a 
filter on ${\Cal P}(\omega)^{{\bold V}[P^{\frak t}_\mu]}$, as
Rang$({\underset\tilde {}\to \eta_i}) \in [\omega]^{\aleph_0}$, \nl
$i < j < \mu \and \text{ cf}(i) \ne \kappa \and \text{ cf}(j) \ne 
\kappa \Rightarrow \text{ Rang}({\underset\tilde {}\to \eta_j}) 
\subseteq^* \text{ Rang} ({\underset\tilde {}\to \eta_i})$.  \nl
Also it is an ultrafilter as
${\Cal P}(\omega)^{\bold V[P^{\frak t}_\mu]} = \dbcu_{i < \mu}
{\Cal P}(\omega)^{\bold V[P^{\frak t}_i]}$ and if $i < \mu$, then Rang
$({\underset\tilde {}\to \eta_{i+1}})$ induces an ultrafilter on
${\Cal P}(\omega)^{\bold V[P^{\frak t}_{i+1}]}$.  So ${\frak u} \le \mu$.  Also
$({}^\omega \omega)^{\bold V[P^{\frak t}_\mu]} = \dbcu_{i < \mu}({}^\omega \omega)
^{\bold V[P^{\frak t}_i]},({}^\omega \omega)^{\bold V[P^{\frak t}_i]}$ is increasing with
$i$ and if cf$(i) \ne \kappa$ then ${\underset\tilde {}\to \eta_i} \in
{}^\omega \omega$ dominate $({}^\omega \omega)^{\bold V[P^{\frak t}_i]}$ by
Subclaim \scite{au.B}, so
${\frak b} = {\frak d} = \mu$ as in previous cases.  Lastly, always
${\frak u} \ge {\frak b}$ hence ${\frak u} = \mu$.]
\bn
Now we define ${\frak t}_\alpha \in {\frak K}$ for $\alpha \le \lambda$ 
by induction on $\alpha$ satisfying $\langle {\frak t}_\alpha:\alpha 
\le \lambda \rangle$ is
$\le_{\frak K}$-increasing continuous such that ${\frak t}_{\alpha +1}$ is
gotten from ${\frak t}_\alpha$ as in Subclaim \scite{au.E}.  Let $P = P^{{\frak t}
_\lambda}_\mu$, so $|P| \le \lambda$ hence $(2^{\aleph_0})^{\bold V^P} \le
(\lambda^{\aleph_0})^{\bold V}$ and easily equality holds. \nl
\bn
We finish by \nl
\ub{\stag{au.G} Subclaim}:  $\Vdash_{P_\lambda} ``{\frak a} \ge \text{ cf}(\lambda)"$.
\nl
[Why?  Assume toward a contradiction that $\theta < \text{ cf}(\lambda)$ and
$p \in P$ and $p \Vdash_P ``{\underset\tilde {}\to {\Cal A}} = \{
{\underset\tilde {}\to A_i}:i < \theta\}$ is a MAD family; i.e.
\mr
\widestnumber\item{$(iii)$}
\item "{$(i)$}"  $A_i \in [\omega]^{\aleph_0}$
\sn
\item "{$(ii)$}"  $i \ne j \Rightarrow |{\underset\tilde {}\to A_i} \cap
{\underset\tilde {}\to A_j}| < \aleph_0$
\sn
\item "{$(iii)$}"  under (i) + (ii),
${\underset\tilde {}\to {\Cal A}}$ is maximal".
\ermn
Without loss of generality $\Vdash_P ``{\underset\tilde {}\to A_i} \in
[\omega]^{\aleph_0}"$.  As ${\frak a} \ge {\frak b} = \mu$ by 
Subclaim F, we have $\theta
\ge \mu$.  For each $i < \theta$ and $m < \omega$ there is a maximal
antichain $\langle p_{i,m,n}:n < \omega \rangle$ of $P$ and there is a
sequence $\langle \bold t_{i,m,n}:n < \omega \rangle$ of truth values such
that $p_{i,m,n} \Vdash ``m \in {\underset\tilde {}\to A_i} \equiv
\bold t_{i,m,n}"$.  We can find countable $w_i \subseteq \mu$ such that
$\dbcu_{m,n < \omega} \text{ Dom}(p_{i,m,n}) \subseteq w_i$. Possibly
increasing $w_i$ retaining countability, we can find $\langle R_{i,\gamma}:
\gamma \in w_i \rangle$ such that
\mr
\item "{$(\alpha)$}"  $w_i$ has a maximal element and $\gamma \in w_i
\backslash \{\text{max}(w_i)\} \Rightarrow \gamma +1 \in w_i$
\sn
\item "{$(\beta)$}"  $R_{i,\gamma}$ is a countable subset of
$P^{{\frak t}_\lambda}_\gamma$ and $q \in R_{i,\gamma} \Rightarrow
\text{ Dom}(q) \subseteq w_i$
\sn
\item "{$(\gamma)$}"  for $\gamma_1 < \gamma_2$ in $w_i,q \in R_{i,\gamma_2}
\Rightarrow q \restriction \gamma_1 \in R_{i,\gamma_1}$
\sn
\item "{$(\delta)$}"  for $\gamma_1 \in w_i,\gamma \in \delta_1 \cap
w_i$ and $q \in R_{i,\gamma_1}$ the $P^{\frak t}_\gamma$-name $q(\gamma)$ involves
$\aleph_0$ maximal antichains all included in $R_{i,\gamma}$
\sn
\item "{$(\varepsilon)$}"  $\{p_{i,m,m}:m,n\} \subseteq 
R_{i,\text{max}(w_i)}$.
\ermn
As cf$(\lambda) > \aleph_0$ (otherwise the conclusion is trivial) we have
$P^{\frak t}_\mu = \dbcu_{\alpha < \lambda} P^{{\frak t}_\alpha}_\mu$.
Clearly for some $\alpha < \lambda$ we have $\cup \{R_{i,\gamma}:i < \theta,
\gamma \in w_i\} \subseteq P^{{\frak t}_\alpha}_\mu$.  
But $P^{{\frak t}_\alpha}_\mu \lessdot P^{{\frak t}_\lambda}_\mu$.  So
$\Vdash_{P^{{\frak t}_\alpha}_\mu} ``{\underset\tilde {}\to {\Cal A}} =
\{{\underset\tilde {}\to A_i}:i < \theta\}$ is MAD".

Now, letting $\bold j$ be the canonical elementary embedding of $\bold V$ into
$\bold V^\kappa/D$, we know:
\mr
\item "{$(*)$}"  in $\bold V^\kappa/D,\bold j({\underset\tilde {}\to {\Cal A}})$ is
a $\bold j(P^{{\frak t}_\alpha}_\mu)$-name of a MAD family.
\ermn
As $\bold V^\kappa/D$
is $\kappa$-closed, for c.c.c. forcing notions things are absolute enough but
$\{\bold j(i):i < \mu\}$ is not $\{i:\bold V^\kappa/D \models i < \bold j
(\mu)\}$, so in $\bold V$, it is forced for $\Vdash_{{\bold j}(P^{\frak t}_\mu)}$,
that $\{\bold j({\underset\tilde {}\to A_i}):i < \mu\}$ is not MAD! \nl
Chasing arrows, clearly $\Vdash_{P^{{\frak t}_{\alpha +1}}_\mu} ``\{
{\underset\tilde {}\to A_i}:i < \theta\}$ is not MAD" as required.
\hfill$\square_{\scite{au.1}}$\margincite{au.1}
\bigskip

\proclaim{\stag{au.2} Theorem}  Assume
\mr
\item "{$(a)$}"  $\kappa$ is a measurable cardinal
\sn
\item "{$(b)$}"  $\kappa < \mu = \text{ cf}(\mu) < \lambda = \text{ cf}
(\lambda) = \lambda^\kappa$.
\ermn
\ub{Then} for some 
c.c.c. forcing notion $P$ of cardinality $\lambda$ in $\bold V^P$ we
have: $2^{\aleph_0} = \lambda,{\frak i} = {\frak u} = {\frak d} = {\frak b}
= \mu$ and ${\frak a} = \lambda$.
\endproclaim
\bigskip

\remark{Remark}  Recall ${\frak i} = \text{ Min}\{|{\Cal A}|:{\Cal A}
\subseteq [\omega]^{\aleph_0}$ is a maximal independent$\}$ where
independent means that every nontrivial Boolean combination of finitely many
members of ${\Cal A}$ is not empty and even infinite.
\endremark
\bigskip

\demo{Proof}  Like the proof of \scite{au.1} except for the following changes.

Let $S_0 = \{i < \mu:\text{cf}(i) = \kappa\},S_1 = \{i < \mu:\text{cf}(i) \ne
\kappa,i$ even$\}$, $S_2 = \{i < \mu:\text{cf}(i) \ne \kappa,i$ odd$\}$, so
$\langle S_0,S_1,S_2 \rangle$ is a partition of $\mu$.  For $i \in S_0,i \in
S_1$ we define ${\frak t} \in {\frak K}$ as 
before (so in clauses (c), (e), we restrict ourselves to $i \in S_1$
and $j \in S_1$) but we add
\mr
\item "{$(j)$}"  for $i \in S_2,{\underset\tilde {}\to D_i}$ is a maximal
filter on $\omega$ containing the co-bounded subsets such that
$\{\text{Rang}({\underset\tilde {}\to \eta_j}):j \in i \cap S_2\}$ is an
independent family modulo ${\underset\tilde {}\to D_i}$
\sn
\item "{$(k)$}"  for $i \in S_2$ and $j \in i \cap S_2$ we have
${\underset\tilde {}\to D_j} \subseteq {\underset\tilde {}\to D_i}$.
\ermn
In Subclaim A, define $\bar Q^{\frak t} \restriction i$ by induction
on $i$ such that $\Vdash_{P_i} ``\langle \text{Rang}({\underset\tilde {}\to \eta_j}):
j \in S_2 \cap i \rangle$ is independent modulo 
$\cup\{{\underset\tilde {}\to D_j}:j \in i \cap S_2\}"$. \nl
Note the definition of $Q({\underset\tilde {}\to D_i})$ remains:
tr$(p) \triangleleft \eta \in p \in Q({\underset\tilde {}\to D_i})
\Rightarrow \{n:\eta \char 94 \langle n \rangle \in p\}$ belongs to $D_i$.
Also note
{\roster
\itemitem{ $\boxtimes_1$ }  ``if $A \subseteq \omega,A \ne \emptyset$ mod
$\underset\tilde {}\to D$ then $\Vdash_{Q(D)}" \text{ Rang}
({\underset\tilde {}\to \eta}) \cap A \ne \emptyset \text{ mod } D$. \nl
However
\sn
\itemitem{ $\boxtimes_2$ }  in Subclaim C we assume
$\Vdash_{P_2} ``{\underset\tilde {}\to D_1} = {\underset\tilde {}\to
D_2} \cap {\Cal P}(\omega)^{\bold V^{P_1}}"$ and in the proof we replace
$\boxtimes(b)(iii)$ by: \nl
\smallskip

\hskip20pt  $(iii)'$  if $\nu \in T$ then $\{n:\nu \char 94 \langle n
\rangle \in T\} \ne \emptyset$ mod $D$ or \nl

\hskip40pt $(\forall n)(\nu \char 94
\langle n \rangle \notin T) \and (\exists q \in {\Cal I})(\nu \in q)$.
\endroster}
\enddemo
\bn
Also we should add \nl
\ub{\stag{au.H} Subclaim}:  If ${\frak t} \in {\frak K}$ then
$\Vdash_{P^{\frak t}_\mu} ``\{\text{Rang}({\underset\tilde {}\to \eta_j}):j
\in S_2\}$ is a maximal independent family of subsets of $\omega$". \nl
[Why?  Independence is covered by the previous paragraph. Assume
towards a contradiction that $p \in P^{\frak t}_\mu$ and 
$p \Vdash ``\underset\tilde {}\to A \in [\omega]^{\aleph_0} 
\backslash \{\text{Rang}
({\underset\tilde {}\to \eta_j}):j \in S_2\}$ and $\{\text{Rang}
({\underset\tilde {}\to \eta_j}):j \in S_2\} \cup \{
\underset\tilde {}\to A\}$ is independent". For some $j \in S_2$ we have
$p \in P^{\frak t}_j$ and $\underset\tilde {}\to A$ is a
$P^{\frak t}_j$-name. So in $\bold V^{P^{\frak t}_j},\{\text{Rang}
({\underset\tilde {}\to \eta_i}):i \in S_2 \cap j\} \cup
\{\underset\tilde {}\to A\}$ is an independent family of subsets of $\omega$,
so by the maximality of ${\underset\tilde {}\to D_j}$ we have: for some
$m < n < \omega,i_\ell \in S_2 \cap j$ for $\ell < n$ with no repetitions we
have $\dbca^{m-1}_{\ell=0} \text{ Rang}({\underset\tilde {}\to \eta_{i_\ell}})
\cap \dbca^{n-1}_{\ell =m} (\omega \backslash \text{ Rang}
({\underset\tilde {}\to \eta_{i_\ell}})) \cap \underset\tilde {}\to A \cap
B =0$ for some $B \in D_j$ (maybe interchange $A$ and its complement), so
\wilog \, $p$ forces this (with $B$ replaced by a $P^{\frak t}_j$-name $B$).
But $\Vdash_{P^{\frak t}_{j+1}} ``\text{Rang}({\underset\tilde {}\to \eta_j})
\cap (\omega \backslash \underset\tilde {}\to B)$ is finite" hence 

$$
p \Vdash `` \dbca^{m-1}_{\ell=0} \text{ Rang}
({\underset\tilde {}\to \eta_{i_\ell}} \cap \dbca^{n-1}_{\ell=m} (\omega
\backslash \text{ Rang}({\underset\tilde {}\to \eta_{i_\ell}}))
  \cap \text{ Rang}({\underset\tilde {}\to \eta_j}) \cap A \text{ is
finite}".
$$
\mn
This contradicts the choice of $\underset\tilde {}\to A$ so we are done.
\hfill$\square$ \nl
${{}}$ \hfill$\square_{\scite{au.2}}$\margincite{au.2}
\bn
\ub{\stag{au.3} Discussion}:  1) We can now look at other problems, like what can be the
order and equalities among ${\frak d}, {\frak b}, {\frak a}, {\frak u},
{\frak i}$; have not considered it. 
\newpage
    
REFERENCES.  
\bibliographystyle{lit-plain}
\bibliography{lista,listb,listx,listf,liste}

\def\germ{\frak} \def\scr{\cal}
  \ifx\documentclass\undefinedcs\def\rm{\fam0\tenrm}\fi
  \def\defaultdefine#1#2{\expandafter\ifx\csname#1\endcsname\relax
  \expandafter\def\csname#1\endcsname{#2}\fi} \defaultdefine{Bbb}{\bf}
  \defaultdefine{frak}{\bf} \defaultdefine{mathbb}{\bf}
  \defaultdefine{mathcal}{\cal}
  \defaultdefine{beth}{BETH}\defaultdefine{cal}{\bf} \def\bbfI{{\Bbb I}}
  \def\mbox{\hbox} \def\text{\hbox} \def\om{\omega} \def\Cal#1{{\bf #1}}
  \def\pcf{pcf} \defaultdefine{cf}{cf} \defaultdefine{reals}{{\Bbb R}}
  \defaultdefine{real}{{\Bbb R}} \def\restriction{{|}} \def\club{CLUB}
  \def\w{\omega} \def\exist{\exists} \def\se{{\germ se}} \def\bb{{\bf b}}
  \def\equivalence{\equiv} \let\lt< \let\gt> \def\cite#1{[#1]}
  \def\implies{\Rightarrow}
\begin{thebibliography}{RoSh 670}
\makeatletter \renewcommand{\@biblabel}[1]{[#1]} \makeatother
\par\noindent [References of the form {\tt math.XX/$\cdots$} refer to the
{\tt xxx.lanl.gov} archive]  \par

\bibitem[RoSh 670]{RoSh:670}Andrzej Ros{\l}anowski and Saharon Shelah.
\newblock {Norms on possibilities III: strange subsets of the real line}.
\newblock {\em {in preparation}}.

\bibitem[Sh 592]{Sh:592}Saharon Shelah.
\newblock {Covering of the null ideal may have countable cofinality}.
\newblock {\em {Fundamenta Mathematicae}}, {\bf to appear}.
 {\tt math.LO/9810181}

\bibitem[Sh 630]{Sh:630}Saharon Shelah.
\newblock {Non-elementary proper forcing notions}.
\newblock {\em {Journal of Applied Analysis}}, {\bf accepted}.
 {\tt math.LO/9712283}

\bibitem[Sh 666]{Sh:666}Saharon Shelah.
\newblock {On what I do not understand (and have something to say)}.
\newblock {\em {Fundamenta Mathematicae}}, {\bf to appear}.
 {\tt math.LO/9906113}

\bibitem[Sh 619]{Sh:619}Saharon Shelah.
\newblock {The null ideal restricted to a non-null set may be saturated}.
\newblock {\em {}}.
 {\tt math.LO/9705213}

\bibitem[Sh 220]{Sh:220}Saharon Shelah.
\newblock {Existence of many $L_ {\infty,\lambda}$-equivalent, nonisomorphic
  models of $T$ of power $\lambda$}.
\newblock {\em {Annals of Pure and Applied Logic}}, {\bf 34}:291--310, 1987.
\newblock Proceedings of the Model Theory Conference, Trento, June 1986.

\bibitem[Sh 300]{Sh:300}Saharon Shelah.
\newblock {Universal classes}.
\newblock In {\em {Classification theory (Chicago, IL, 1985)}}, volume 1292 of
  {\em {Lecture Notes in Mathematics}}, pages 264--418. {Springer, Berlin},
  1987.
\newblock {Proceedings of the USA--Israel Conference on Classification Theory,
  Chicago, December 1985; ed. Baldwin, J.T.}

\bibitem[Sh:c]{Sh:c}Saharon Shelah.
\newblock {\em {Classification theory and the number of nonisomorphic models}},
  volume~92 of {\em {Studies in Logic and the Foundations of Mathematics}}.
\newblock {North-Holland Publishing Co., Amsterdam, xxxiv+705 pp}, 1990.

\bibitem[Sh:f]{Sh:f}Saharon Shelah.
\newblock {\em {Proper and improper forcing}}.
\newblock {Perspectives in Mathematical Logic}. {Springer}, 1998.

\bibitem[vD]{vD}Eric~K. van Douwen.
\newblock {The integers and topology}.
\newblock In K.~Kunen and J.~E. Vaughan, editors, {\em Handbook of
  Set-Theoretic Topology}, pages 111--167. Elsevier Science Publishers, 1984.

\end{thebibliography}

\enddocument

\shlhetal   

\bye